\documentclass[12pt]{amsart}
\usepackage{amsfonts, amsmath, latexsym, epsfig}
\usepackage{amssymb}
\usepackage{epsf}
\usepackage{url}
\usepackage{txfonts}

\newtheorem{theorem}{Theorem}

\newtheorem{conjecture}{Conjecture}

\def\QuotS#1#2{\leavevmode\kern-.0em\raise.2ex\hbox{$#1$}\kern-.1em/\kern-.1em\lower.25ex\hbox{$#2$}}

\begin{document}

\author{Mathieu Dutour Sikiri\'c}
\address{Mathieu Dutour Sikiri\'c, Rudjer Boskovi\'c Institute, Bijenicka 54, 10000 Zagreb, Croatia}
\email{mdsikir@irb.hr}

\author{Michel Deza}
\address{Michel Deza, M. Deza, \'Ecole Normale Sup\'erieure, Paris}
\email{Michel.Deza@ens.fr}

\author{Mikhail Shtogrin}
\address{Mikhail Shtogrin, Steklov Mathematical Institute, Moscow}
\email{stogrin@mi.ras.ru}

\thanks{First author has been supported by the Croatian Ministry of Science, Education and Sport under contract 098-0982705-2707. Third author has been supported by Russian Federation grant RFFI 11-01-00633.}

\title{Fullerene-like spheres with faces of negative curvature}
\date{}

\maketitle

\begin{abstract}
Given $R\subset \mathbb{N}$, an  {\em $(R,k)$-sphere} is a $k$-regular 
map on the sphere whose faces have gonalities $i\in R$. 
The most interesting/useful are (geometric) {\em fullerenes}, i.e.,
$(\{5,6\},3)$-spheres.

Call $\kappa_i=1 + \frac{i}{k} - \frac{i}{2}$
the {\em curvature} of $i$-gonal faces. $(R,k)$-spheres  admitting 
$\kappa_i<0$
are much harder to study. We consider the symmetries and construction for
three new instances of such spheres:
$(\{a,b\},k)$-spheres with $p_b\le 3$ (they are listed),
{\em icosahedrites} (i.e., $(\{3,4\},5)$-spheres) and,
for any $c\in \mathbb{N}$, {\em fullerene $c$-disks}, i.e., 
$(\{5,6,c\},3)$-spheres with $p_c=1$. 
\end{abstract}

\section{Introduction}

Given $R\subset \mathbb{N}$, an  {\em $(R,k)$-sphere}
$S$ is
a $k$-regular map on the sphere
whose faces have {\em gonalities} (numbers
of sides) $i\in R$.
Let $v,e$ and $f=\sum_{i}p_i$ be the numbers of
vertices, edges and faces of $S$, where $p_i$ is the number
of $i$-gonal faces.
A graph is called {\em $k$-connected} if after removing and $k-1$
vertices, it remains connected.

Clearly, $k$-regularity implies $k v=2e=\sum_{i}i p_i$ and
the {\em Euler formula} $2=v-e+f$ become Gauss-Bonnet-like one
$2=\sum_{i}\kappa_ip_i$, where $\kappa_{i}=1+\frac{i}{k}-\frac{2}{2}$
is called (dualizing the definition in \cite{Hi01})
the {\em curvature} of the $i$-gonal faces.

Let $a=\min\{i\in R\}$.
Then, besides the cases $k=2$ ($a$-cycle) and 
exotic cases $a=1,2$, it holds
\begin{equation*}
\frac{2k}{k-2}>a>2<k<\frac{2a}{a-2},
\end{equation*}
i.e., $(a,k)$ should belong to the five Platonic parameter
pairs $(3,3)$, $(4,3)$, $(3,4)$, $(5,3)$ or $(3,5)$.

Call an $(R,k)$-sphere {\em standard} if $\min_{i\in R}\kappa_i=0$,   
i.e., $b=\frac{2k}{k-2}$, where $b$ denotes $\max\{i\in R\}$.
Such spheres have $(b,k)=(3,6)$, $(4,4)$, $(6,3)$, i.e., the three
Euclidean parameter pairs.
Exclusion of faces of negative curvature simplifies enumeration,
while the number $p_b$ of faces of curvature zero not being
restricted, there is an infinity of such $(R,k)$-spheres.

An {\em $(\{a,b\},k)$-sphere} is an  $(R,k)$-sphere
with $R=\{a,b\}$, $1\le a, b$. 
Clearly, all possible $(a,b;k)$ for the standard $(\{a,b\},k)$-spheres 
are:
\begin{equation*}
(5,6;3), (4,6;3), (3,6;3), (2,6;3), (3,4;4), (2,4;4), (2,3;6), (1,3;6).
\end{equation*}
Those eight families can be seen as spheric analogs of the
regular plane partitions $\{6,3\}$, $\{4,4\}$, $\{3,6\}$
with $p_a=\frac{2b}{b-a}$  $a$-gonal ``defects'' $\kappa_a$ added to
get the total curvature $2$ of the sphere.
$(\{5,6\},3)$-spheres are
(geometric) {\em fullerenes}, important in Chemistry, while
$(\{a,b\},4)$-spheres are minimal projections of
{\em alternating links}, whose components are their  {\em central circuits}
(those going only ahead) and crossings are the vertices.

We considered above eight families in 
\cite{zig2,oct2,covcent,octa,selfdual,zig4} and the book \cite{book3},
where $i$-faces with $\kappa_i<0$
are allowed, mainly in Chapters 15-19.
Here we consider three new natural instances of $(\{a,b,c\},k)$-spheres, 
each allowing faces of negative curvature.
The first Section describes such $(\{a,b\},k)$-spheres with
$p_b\le 3$.
The second Section concerns the {\em icosahedrites}, i.e.,
$(\{3,4\},5)$-spheres, in which $4$-gonal faces have $\kappa_4=-\frac{1}{5}$.
The third Section treats $(\{a,b,c\},k)$-spheres with $p_c=1$,
in which unique $c$-gonal face can be of negative curvature, especially,
{\em fullerene $c$-disks}, i.e., $(\{5,6,c\},3)$-spheres with $p_c=1$.

Note that all  $(R,k)$-spheres with $1,2\notin R$ and $\kappa_i>0$
for all $i\in R$ have $k=4,5$ or $3$ and, respectively,
$\kappa_i=\frac{1}{4}$, $\frac{1}{10}$ or
$\kappa_i\in \{\frac{1}{6}, \frac{1}{3}, \frac{1}{2}\}$.
So, they are only Octahedron, Icosahedron and
eleven $(\{3,4,5\},3)$-spheres: eight dual deltahedra, Cube and its
truncations on one or two opposite vertices ({\em D\"urer octahedron}).

\begin{figure}
\begin{center}
\begin{minipage}[b]{12.1cm}
\centering
\epsfig{height=3.6cm, file=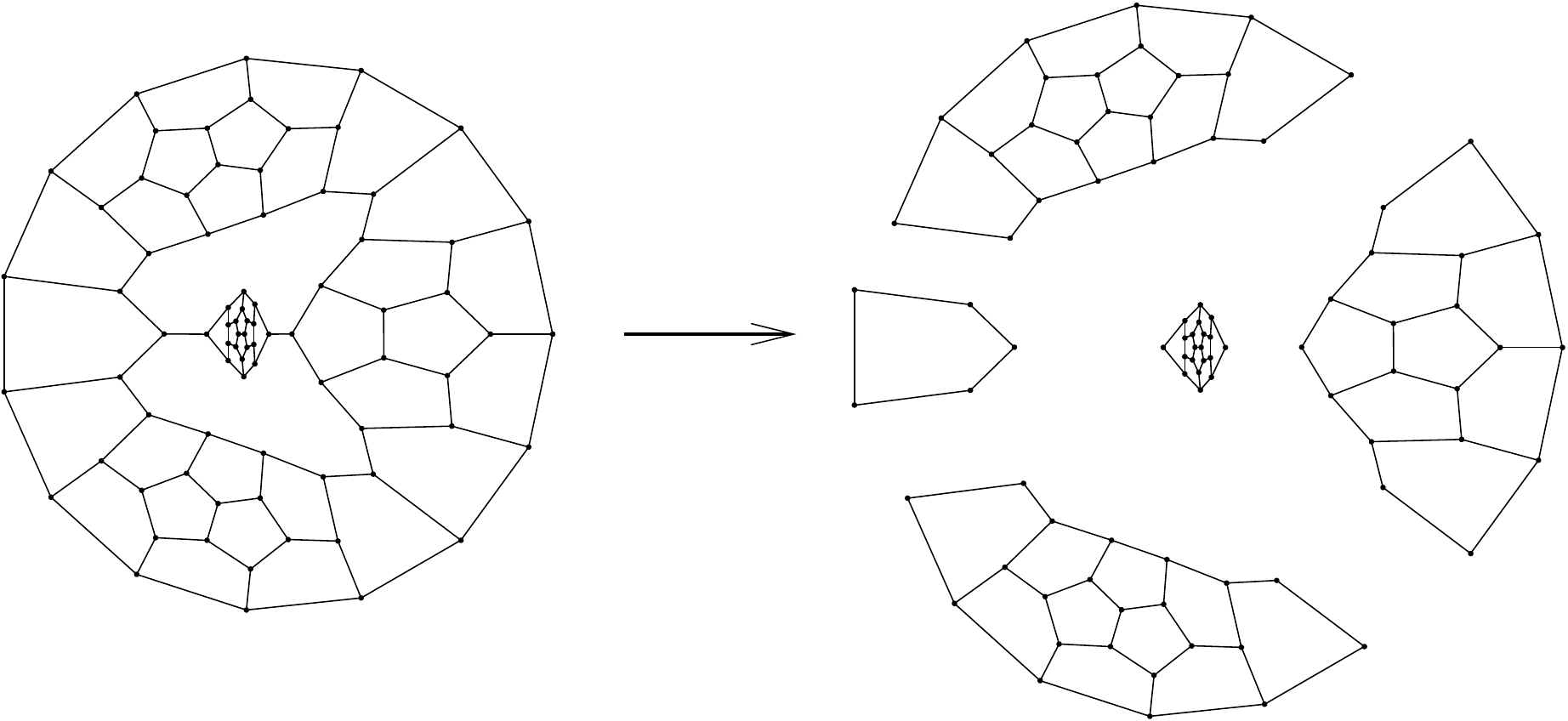}\par
\end{minipage}
\end{center}
\caption{The decomposition of a $(\{5,15\},3)$-sphere 
seen as a $(5,3)$-polycycle with two, $15$- and $18$-gonal, holes,
into elementary polycycles: two 
edge-split $\{5,3\}$'s, $Pen_1$, $C_2$ and $\{5,3\} - e$ in the middle}
\label{ExampleDecompo}
\end{figure}

\begin{figure}
\begin{center}
\begin{minipage}[b]{3.0cm}
\centering
\epsfig{figure=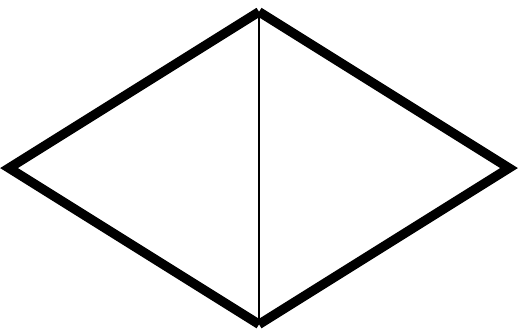,height=1.5cm}\par
$\{3,3\}-e$
\end{minipage}
\begin{minipage}[b]{3.0cm}
\centering
\epsfig{figure=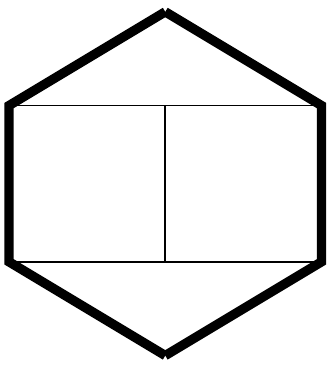,height=1.5cm}\par
$\{4,3\}-e$
\end{minipage}
\begin{minipage}[b]{3.0cm}
\centering
\epsfig{figure=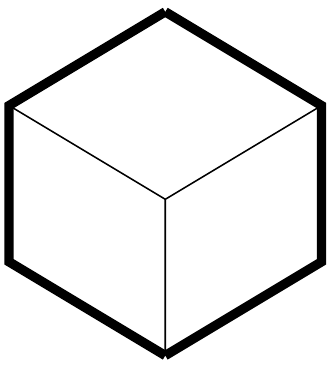,height=1.5cm}\par
$\{4,3\}-v$
\end{minipage}
\begin{minipage}[b]{3.0cm}
\centering
\epsfig{figure=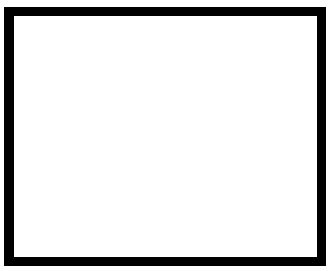,height=1.4cm}\par
$Sq_1$=$\{4,3\}-f$
\end{minipage}
\begin{minipage}[b]{3.0cm}
\centering
\epsfig{figure=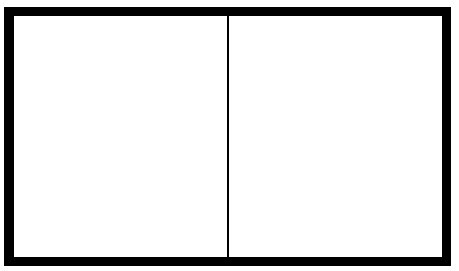,height=1.0cm}\par
$Sq_2$
\end{minipage}
\begin{minipage}[b]{3.0cm}
\centering
\epsfig{figure=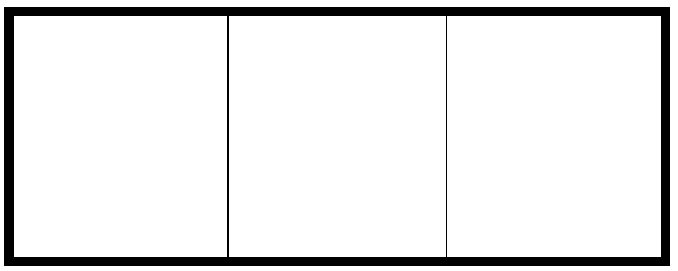,height=1.0cm}\par
$Sq_3$
\end{minipage}
\begin{minipage}[b]{2.3cm}
\centering
\epsfig{file=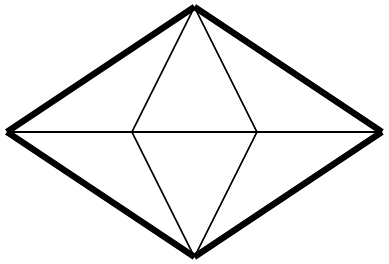, height=1.5cm}\par
$\{3,4\}-e$
\end{minipage}
\begin{minipage}[b]{3.7cm}
\centering
\epsfig{file=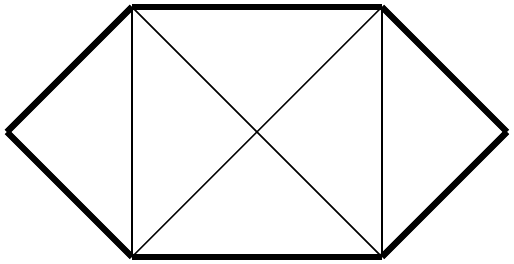, height=1.5cm}\par
vertex-split $\{3,4\}$
\end{minipage}
\begin{minipage}[b]{2.3cm}
\centering
\epsfig{file=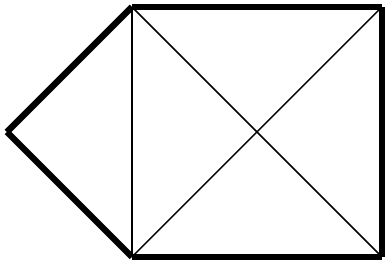, height=1.4cm}\par
$\{3, 4\}-P_3$
\end{minipage}
\begin{minipage}[b]{2.8cm}
\centering
\epsfig{figure=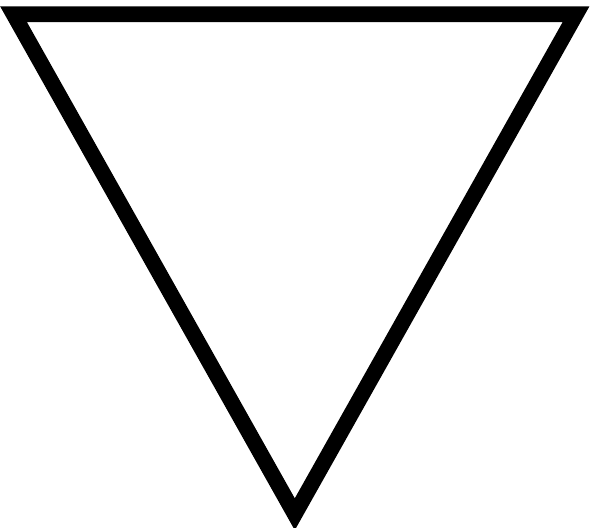,height=1.4cm}\par
$Tr_1=\{3,4\}-f$
\end{minipage}
\begin{minipage}[b]{2.5cm}
\centering
\epsfig{figure=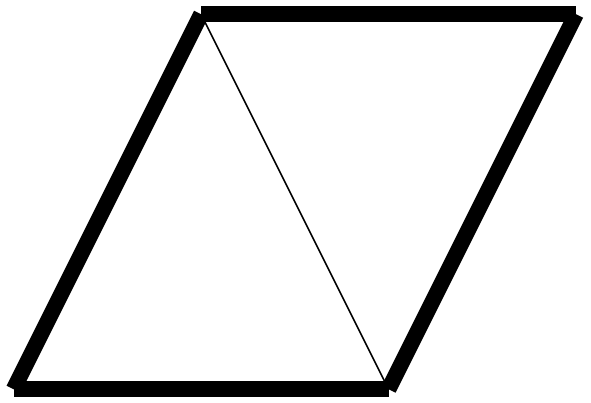,height=1.4cm}\par
$Tr_2$
\end{minipage}
\begin{minipage}[b]{2.8cm}
\centering
\epsfig{figure=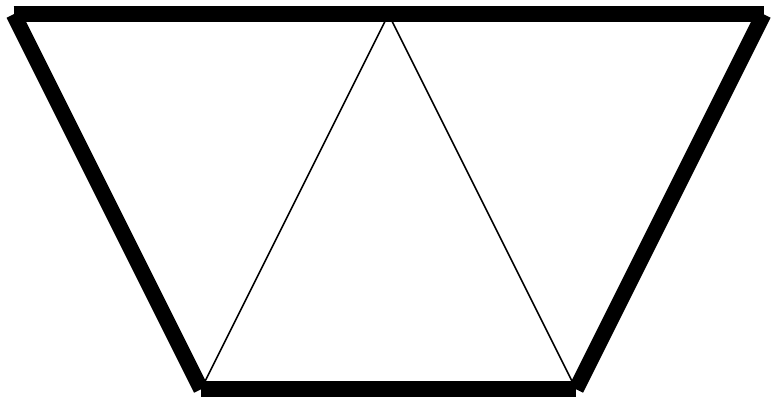,height=1.4cm}\par
$Tr_3$
\end{minipage}
\begin{minipage}[b]{2.8cm}
\centering
\epsfig{figure=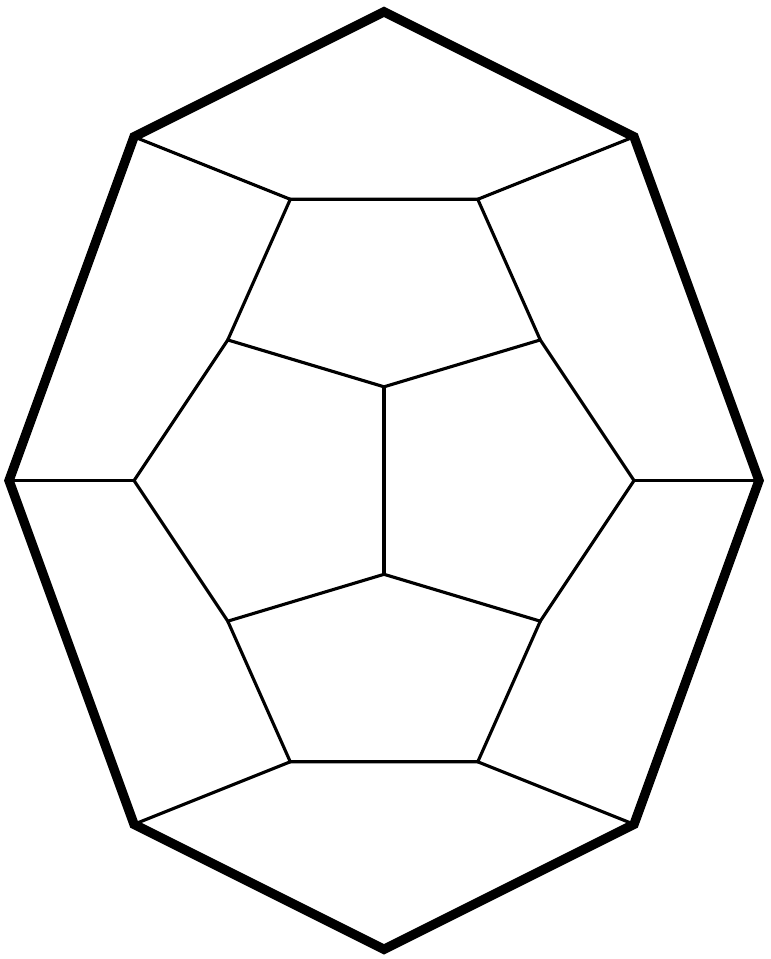,height=1.5cm}\par
$\{5,3\} - e$
\end{minipage}
\begin{minipage}[b]{2.8cm}
\centering
\epsfig{figure=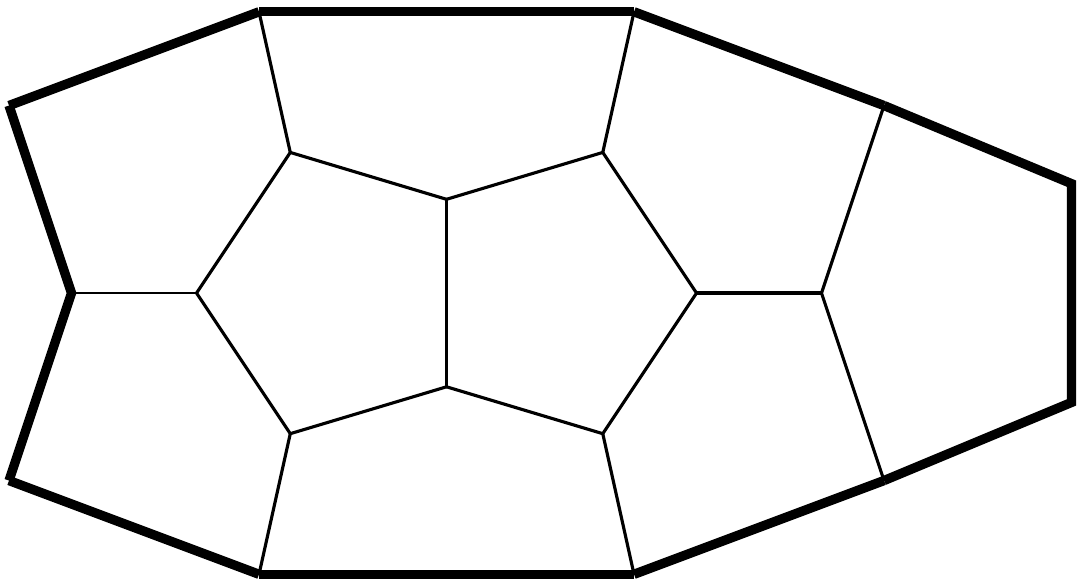,height=1.5cm}\par
$B_2$
\end{minipage}
\begin{minipage}[b]{3.8cm}
\centering
\epsfig{figure=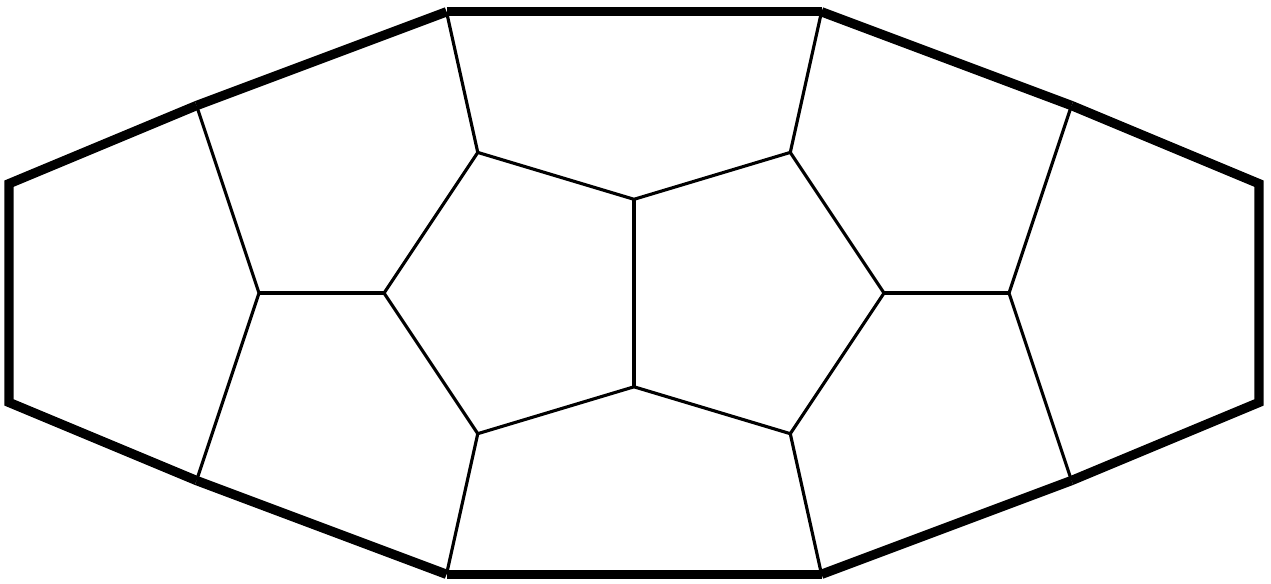,height=1.5cm}\par
edge-split $\{5,3\}$
\end{minipage}
\begin{minipage}[b]{2.8cm}
\centering
\epsfig{figure=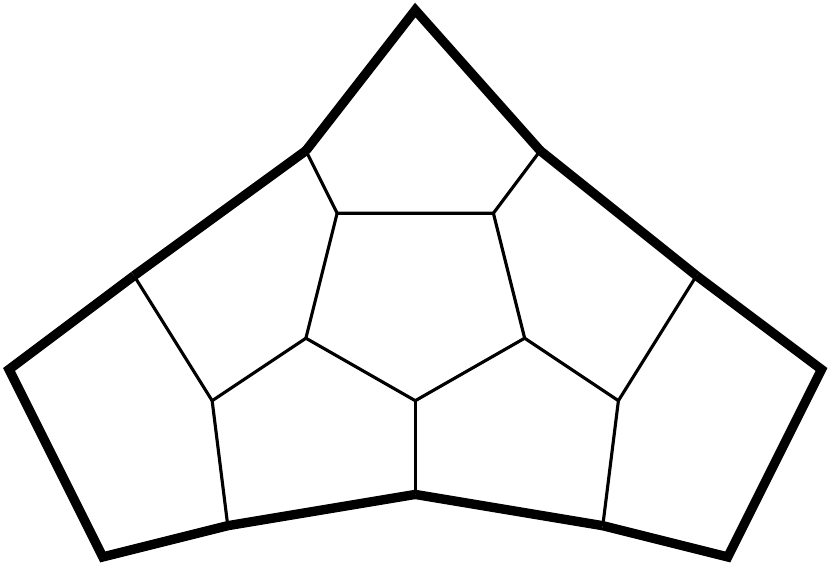,height=1.5cm}\par
$C_2$
\end{minipage}
\begin{minipage}[b]{2.8cm}
\centering
\epsfig{figure=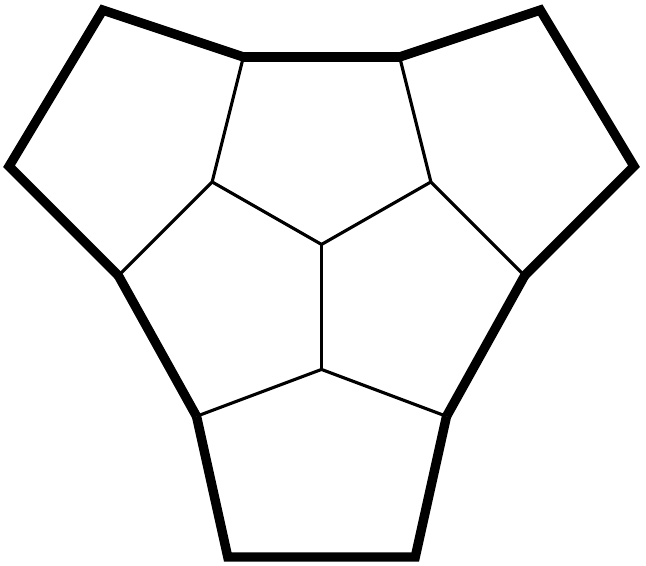,height=1.5cm}\par
$C_3$
\end{minipage}
\begin{minipage}[b]{2.8cm}
\centering
\epsfig{figure=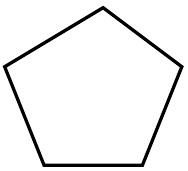,height=1.5cm}\par
$Pen_1$
\end{minipage}
\begin{minipage}[b]{2.8cm}
\centering
\epsfig{figure=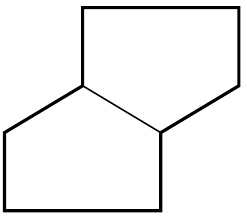,height=1.5cm}\par
$Pen_2$
\end{minipage}
\begin{minipage}[b]{2.8cm}
\centering
\epsfig{figure=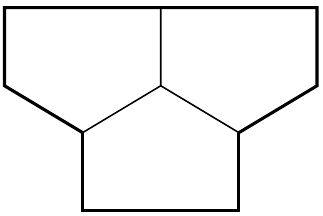,height=1.5cm}\par
$Pen_3$
\end{minipage}
\begin{minipage}[b]{2.8cm}
\centering
\epsfig{figure=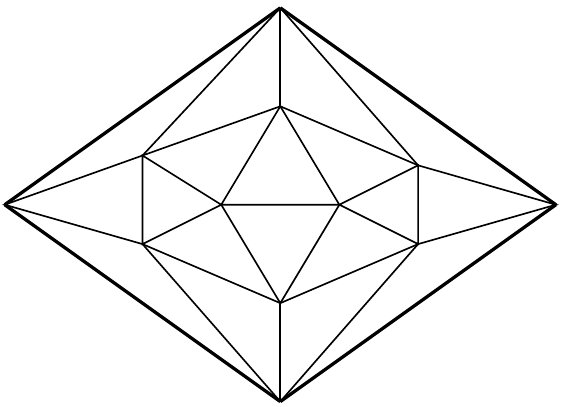,height=1.5cm}\par
$\{3,5\} -e$
\end{minipage}
\begin{minipage}[b]{2.8cm}
\centering
\epsfig{figure=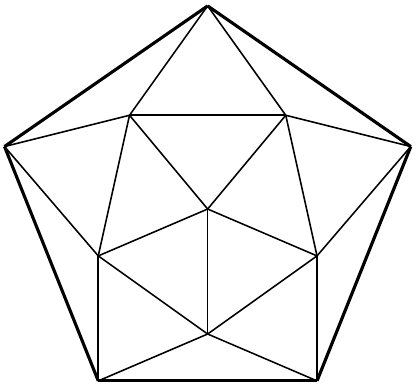,height=1.5cm}\par
$\{3,5\} - v$
\end{minipage}
\begin{minipage}[b]{2.8cm}
\centering
\epsfig{figure=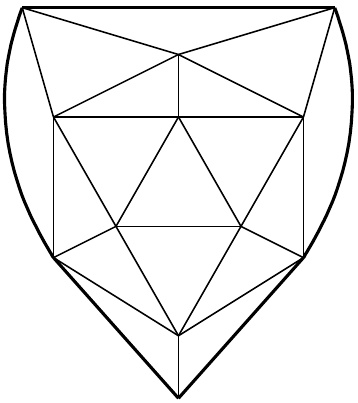,height=1.5cm}\par
$a_3$
\end{minipage}
\begin{minipage}[b]{2.8cm}
\centering
\epsfig{figure=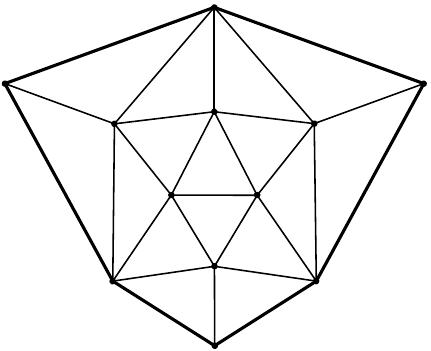,height=1.5cm}\par
$a_4$
\end{minipage}
\begin{minipage}[b]{5.0cm}
\centering
\epsfig{figure=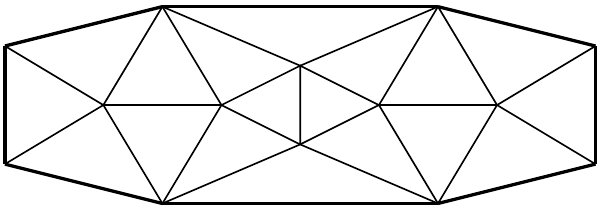,height=1.5cm}\par
edge-split $\{3,5\}$
\end{minipage}
\begin{minipage}[b]{3.7cm}
\centering
\epsfig{figure=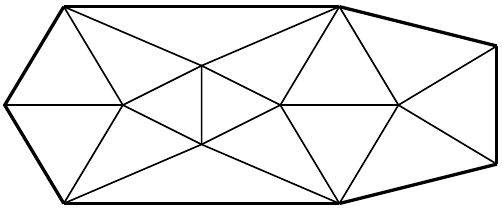,height=1.5cm}\par
$b_3$
\end{minipage}
\begin{minipage}[b]{3.0cm}
\centering
\epsfig{figure=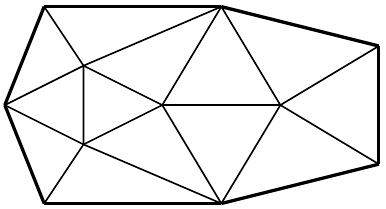,height=1.5cm}\par
$b_4$
\end{minipage}
\begin{minipage}[b]{2.8cm}
\centering
\epsfig{figure=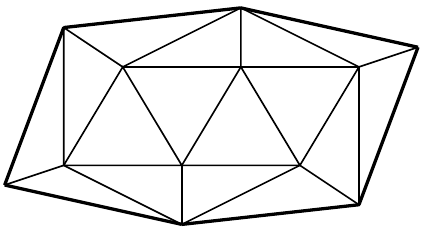,height=1.5cm}\par
$c_1$
\end{minipage}
\begin{minipage}[b]{2.8cm}
\centering
\epsfig{figure=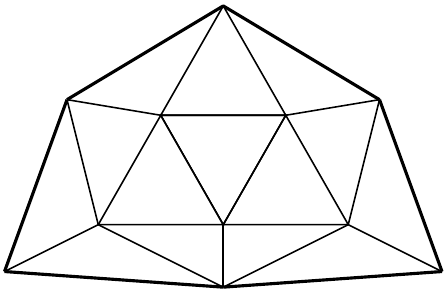,height=1.5cm}\par
$c_2$
\end{minipage}
\begin{minipage}[b]{2.8cm}
\centering
\epsfig{figure=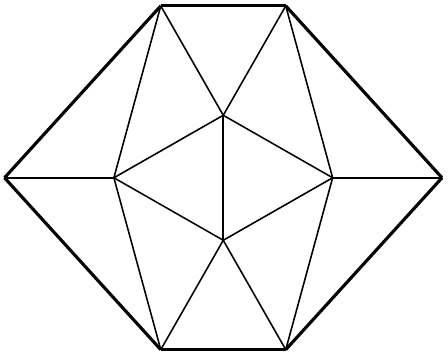,height=1.5cm}\par
$c_3$
\end{minipage}
\begin{minipage}[b]{2.8cm}
\centering
\epsfig{figure=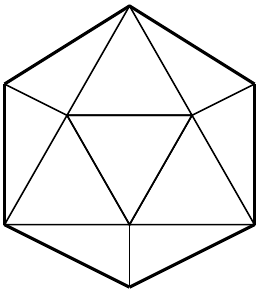,height=1.5cm}\par
$c_4$
\end{minipage}
\begin{minipage}[b]{2.8cm}
\centering
\epsfig{figure=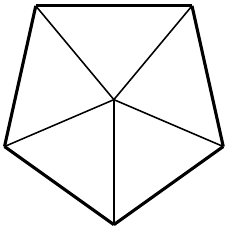,height=1.5cm}\par
$Sun_1$
\end{minipage}
\begin{minipage}[b]{2.8cm}
\centering
\epsfig{figure=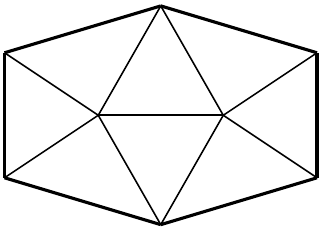,height=1.5cm}\par
$Sun_2$
\end{minipage}
\begin{minipage}[b]{2.8cm}
\centering
\epsfig{figure=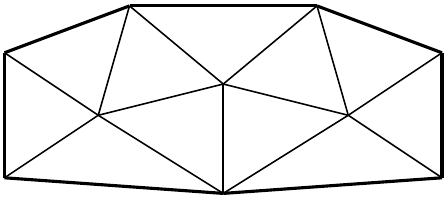,height=1.5cm}\par
$Sun_3$
\end{minipage}

\end{center}
\caption{$(a,k)$-polycycles from Chapters 4,7 of \cite{book3} used in proofs of Theorems \ref{Case_Pb2} and \ref{Case_PB3}}
\label{ExamplePolycycle}
\end{figure}

\section{The list of $(\{a, b\}, k)$-spheres with $p_b\leq 3$}

For $a,k\geq 3$, a {\em $(a,k)$-polycycle} is a plane graph whose faces,
besides some disjoint pairwisely, including exterior one and called
{\em holes}, are $a$-gons, and whose vertices have degree between $2$ and 
$k$ with vertices not on the boundary of holes being $k$-valent.
Let as see any $\{a,b\}),k$-sphere, after merging its adjacent $b$-gons
in larger faces, as a $(a,k)$-polycycle.
For a given $(\{a, b\}, k)$-sphere, we can first remove the edges contained
in two $b$-gons.
If a vertex has a clockwise list of incident faces of the form $b^{x_1} a^{y_1}\dots b^{x_N} a^{y_N}$ with $N\geq 2$ and $x_i$, $y_i\geq 1$, then we split it into $N$ different vertices.
The remaining faces are $a$-gonal and are organized into one or
more $(a,k)$-polycycles
with the pair $(a,k)$ being one of $(3,3)$, $(3,4)$, $(4,3)$, $(3,5)$
and $(5,3)$.

The $(3,3)$-, $(3,4)$- and $(4,3)$-polycycles are easily classified 
(see pages 45,46 of \cite{book3})
and this gives a method for solving
the problem of this section. For the remaining two cases, we have to
introduce another method. 
An $(a,k)$-polycycle is called {\em elementary} if it cannot be cut along
an edge into two $(a,k)$-polycycle.
An $(a,k)$-polycycle admits a unique decomposition into elementary
$(a,k)$-polycycles and the list of elementary $(5,3)$-, $(3,5)$-polycycles
is given in pages 75, 76 of \cite{book3}.

From the set of elementary polycycles, appearing in an $(\{a,b\},k)$-sphere,
we can form a {\em decomposition graph} with the vertices being the
occurring polycycles and two polycycles adjacent if they share an edge,
a vertex or are connected by an edge contained in two $b$-gons.
The connecting vertices and edges are called {\em active}.
See in Figure \ref{ExampleDecompo} an example of such a decomposition.

See the list of polycycles used in this paper on
Figure \ref{ExamplePolycycle}.
Here $Sq_i, Tr_i$ denote horizontal paths of $i$, respectively, squares and  
triangles, while series $Pen_i$, $Sun_i$ are defined similarly.
Note that $Tr_2=\{3,3\}-e$ and $Tr_4=\{3,4\}-P_4$.
We define {\em vertex-split $\{3,5\}$} as unique $(3,5)$-polycycle
obtained from $a_3$ by adjoining a $Tr_1$, and
{\em face-split $\{3,5\}$} as unique $(3,5)$-polycycle
obtained from $c_1$ by adding two $Tr_1$ on two opposite edges.

Note that the $5$-gons of unique minimal fullerene $c$-disk with $c=4$,
$7$, $8$ and $c\geq  13$ given in Figure \ref{MinimalNanodisk} are
organized, respectively, into edge-split $\{5,3\}$, $C_3 + Pen_7$,
$Pen_7 + Pen_7$ and $B_2+Pen_{c-12} + B_2$.

\begin{theorem}\label{Case_PB1}
There is no $(\{a,b\},k)$-sphere with $p_b=1$.
\end{theorem}
\proof The decomposition graph of such spheres is a tree.
If this tree is reduced to a vertex, then the occurring polycycle has
an exterior face being a $b$-gon and an examination of the possibilities
gives $a=b$.
Otherwise, we have at least one polycycle with a unique active vertex
or edge. But an inspection of the list of elementary
polycycles gives that no such one satisfies the condition. \qed

Note that  $(\{6,8\},3)$-maps with $p_6=1$ exist on an oriented  
surface of genus $3$.

Clearly, all $(\{a,b\},k)$-spheres with $a,k>2$ and
$p_b=0$ are five Platonic ones denoted by $\{a,k\}$: Tetrahedron,
Cube ($Prism_4$), Octahedron ($APrism_3$),  Dodecahedron (snub $Prism_5$)
and Icosahedron (snub $APrism_3$).

There exists unique $3$-connected {\em trivial} $(\{a,b\},k)$-sphere
with $p_b=2$ for $(\{4,b\},3)$-, $(\{3,b\},4)$-,  $(\{5,b\},3)$-,
$(\{3,b\},5)$-: $Prism_b$ $D_{bh}$, $APrism_b$ $D_{bd}$,
{\em snub $Prism_b$} $D_{bd}$,
{\em snub $APrism_b$} $D_{bd}$, i.e., respectively,
two $b$-gons separated by $b$-ring of $4$-gons, $2b$-ring of $3$-gons,
two $b$-rings of $5$-gons, two $3b$-rings of $3$-gons.

Clearly, for any $t\geq 2$, there are $2$-connected $b$-vertex $(\{2,b=2t\},3)$- and $(\{2,b\},2t)$-sphere with $p_b=2$: a circle with $t$ disjoint $2$-gons
put on it and a $b$-gon with every edge repeated $t$ times.

\begin{theorem}\label{Case_Pb2}
Let $b > a >2$, $k>2$.
For any non-trivial $(\{a,b\},k)$-sphere with $p_b=2$, the number
$t=\frac{b}{a}$ is an integer.
The list of such spheres consists of following $10$ spheres for each
$t\geq 2$:

(i) For $(a,k)=(3,3)$, $(4,3)$, $(5,3)$, $(3,4)$, $(3,5)$,
the $(\{a,ta\},k)$-sphere $D_{th}$ obtained by putting
on a circle $t$ polycycles $\{a,k\} - e$.
Those polycycles are connected by an edge to their neighbors
and so, only $2$-connected.

(ii) For $(a,k)=(3,4)$, $(5,3)$, $(3,5)$,
the $(\{a,ta\},k)$-sphere $D_{th}$ obtained by 
partitioning of a circle into $t$ polycycles:
respectively, vertex-split $\{3,4\}$, edge-split $\{5,3\}$,
edge-split $\{3,5\}$.

(iii) The $(\{3,3t\}, 5)$-spheres $C_{th}$, 
$D_t$, obtained by partitioning of a circle into $t$ polycycles: 
respectively, vertex-split $\{3,5\}$ ($Tr_1+a_3)$ 
and face-split $\{3,5\}$ ($Tr_1 + c_1 + Tr_1$).

\end{theorem}
\proof If the two $b$-gonal faces are separated by an elementary 
polycycle, then we are in the case of the snub $Prism_b$ or snub $APrism_b$.
Otherwise, the decomposition graph should contain one cycle separating two
$b$-gons. Any nontrivial path connected to this cycle would have a vertex 
of degree $1$ which we have seen to be impossible. 
So, the decomposition graph is reduced to this cycle.
Examination of the list of polycycles
with exactly two connecting edges/vertices and consideration of all
possibilities gives the above list. \qed

Among above spheres, only those coming from edge-split $\{5,3\}$ 
edge-split $\{3,5\}$ 
and face-split $\{3,5\}$ are $3$-connected. 
Those coming from vertex-split $\{3,4\}$ and vertex-split $\{3,5\}$
are $3$-edge connected, but only $2$-(vertex)-connected. 

Let us address now the case $p_b=3$. Note that $3\times K_2$ with $t$ 
disjoint $2$-gons put on each edge is a $(\{2, b=2+4t\}, 3)$-sphere with
$p_b=3$.

\begin{figure}
\begin{center}
\begin{minipage}[b]{3.0cm}
\centering
\epsfig{height=2.4cm, file=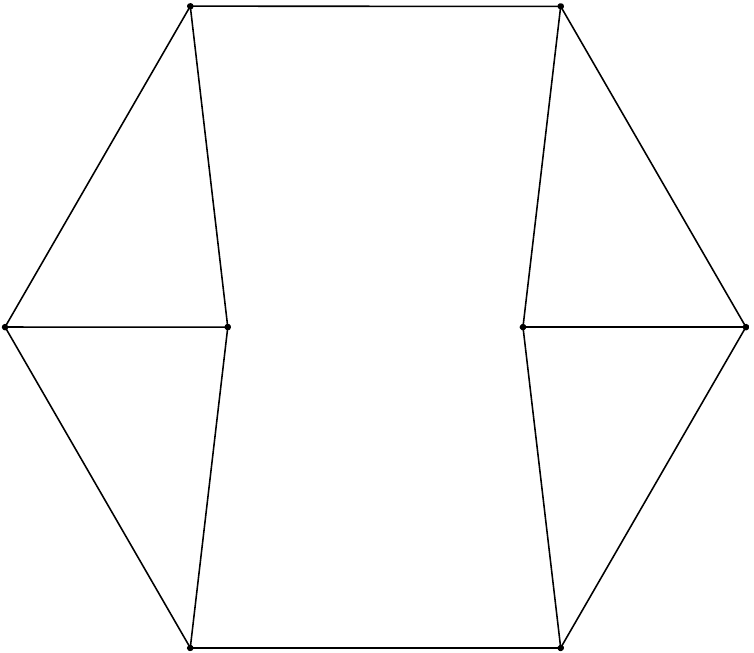}\par
$(\{3,6\},3)$: $D_{2h}$ (i)
\end{minipage}
\begin{minipage}[b]{3.0cm}
\centering
\epsfig{height=2.4cm, file=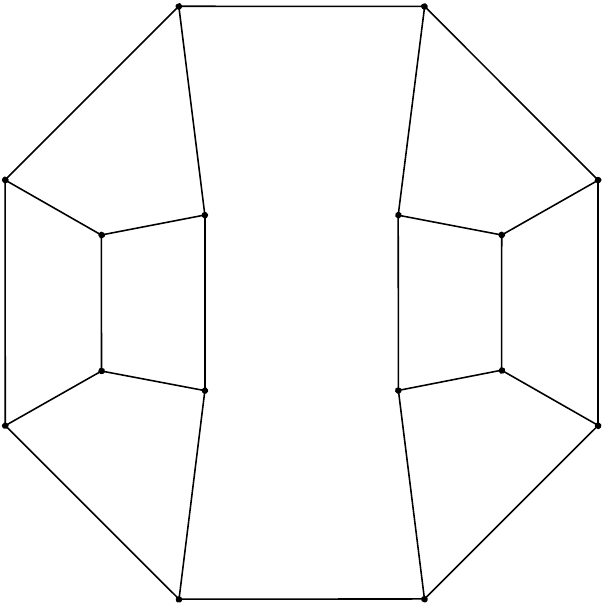}\par
$(\{4,8\},3)$: $D_{2h}$ (i)
\end{minipage}
\begin{minipage}[b]{3.0cm}
\centering
\epsfig{height=2.4cm, file=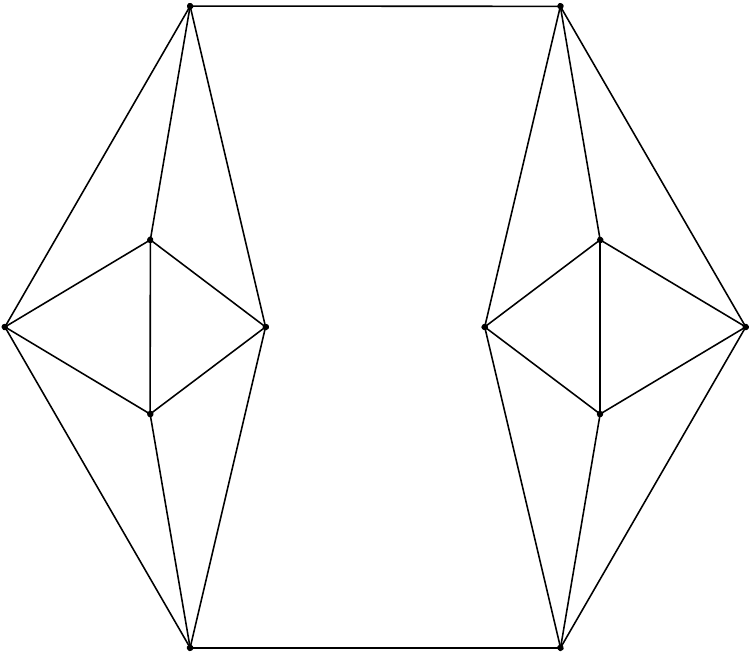}\par
$(\{3,6\},4)$: $D_{2h}$ (i)
\end{minipage}
\begin{minipage}[b]{3.1cm}
\centering
\epsfig{height=2.4cm, file=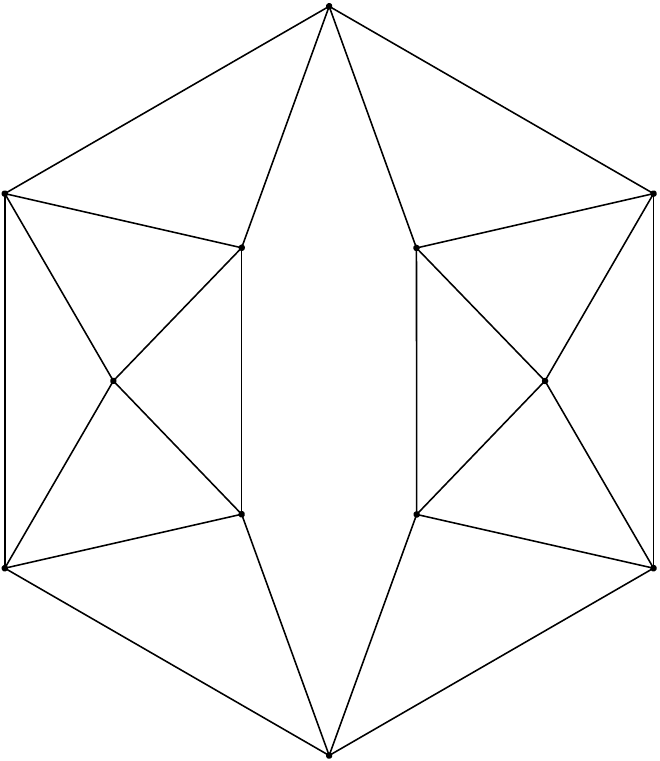}\par
$(\{3,6\},4)$: $D_{2h}$ (ii)
\end{minipage}
\begin{minipage}[b]{3.2cm}
\centering
\epsfig{height=2.4cm, file=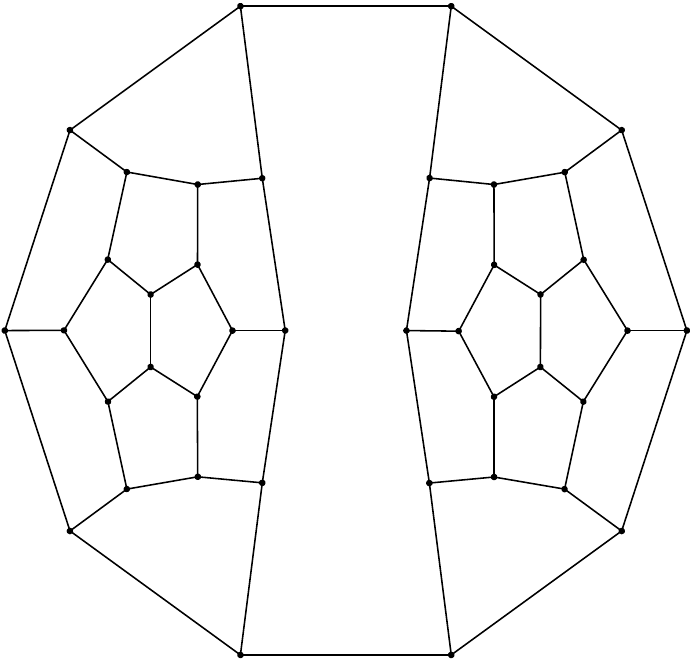}\par
$(\{5,10\},3)$: $D_{2h}$ (i)
\end{minipage}
\begin{minipage}[b]{3.3cm}
\centering
\epsfig{height=2.4cm, file=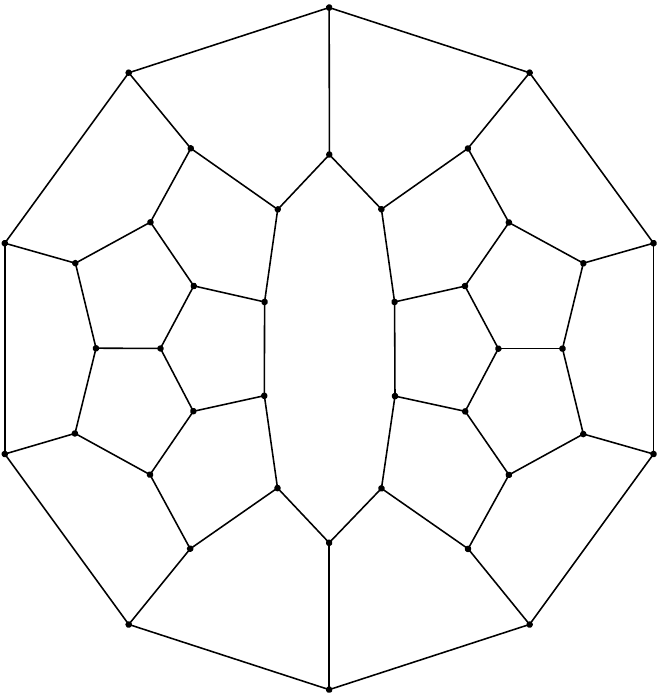}\par
$(\{5,10\},3)$: $D_{2h}$ (ii)
\end{minipage}

\begin{minipage}[b]{3.0cm}
\centering
\epsfig{height=2.4cm, file=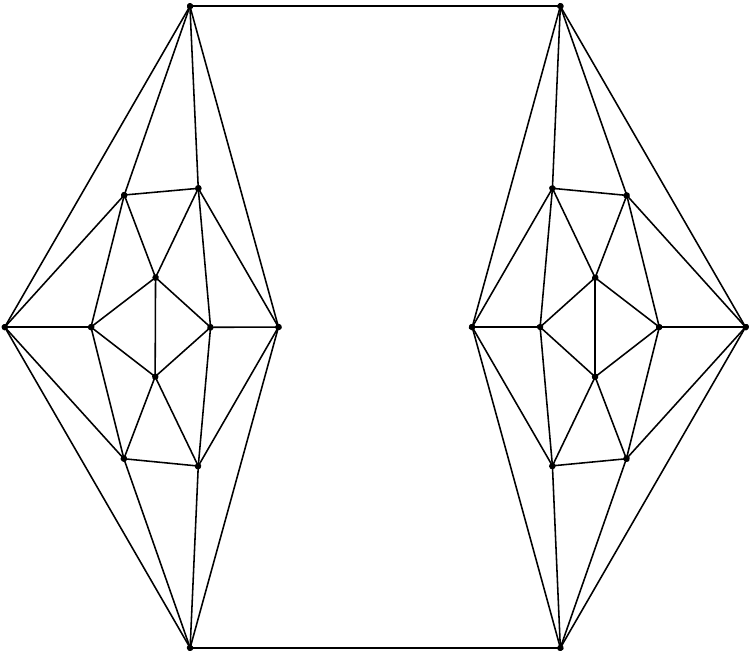}\par
$(\{3,6\},5)$: $D_{2h}$ (i)
\end{minipage}
\begin{minipage}[b]{3.1cm}
\centering
\epsfig{height=2.4cm, file=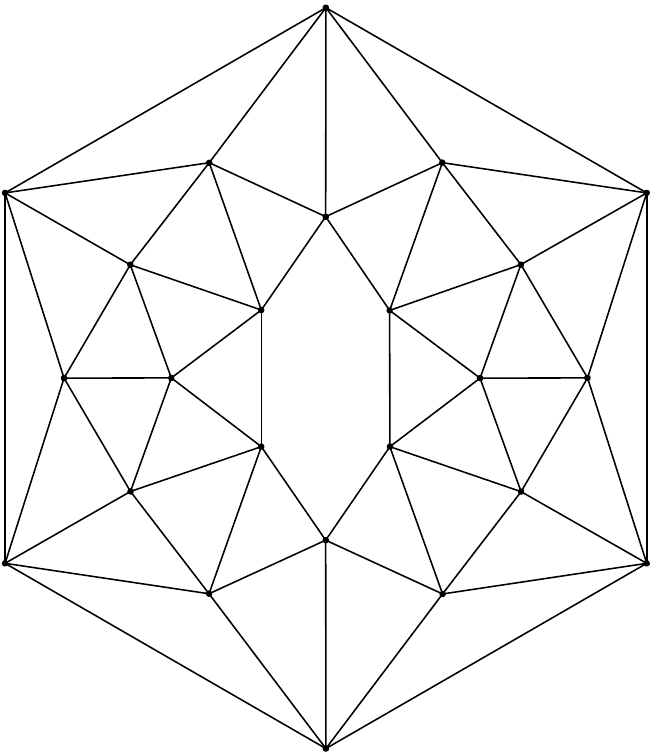}\par
$(\{3,6\},5)$: $D_{2h}$ (ii)
\end{minipage}
\begin{minipage}[b]{3.1cm}
\centering
\epsfig{height=2.4cm, file=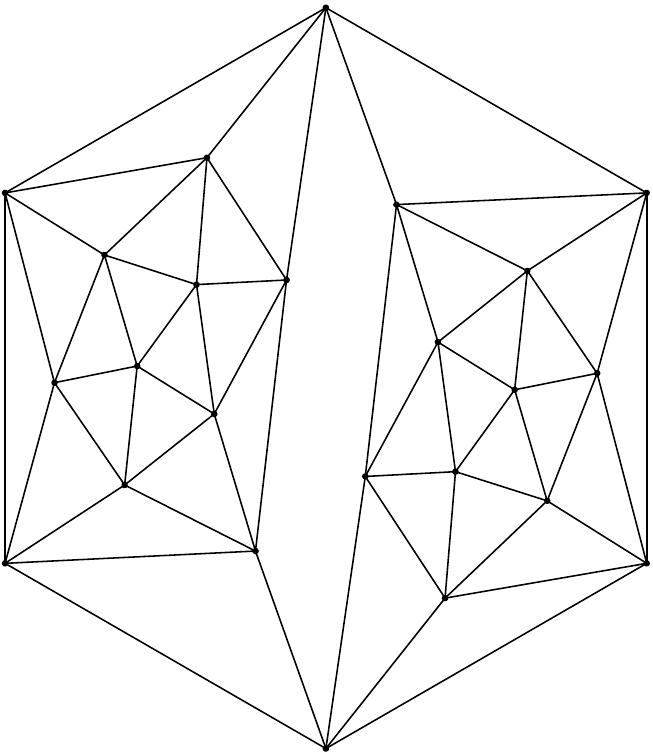}\par
$(\{3,6\},5)$: $C_{2h}$ (iii)
\end{minipage}
\begin{minipage}[b]{3.0cm}
\centering
\epsfig{height=2.4cm, file=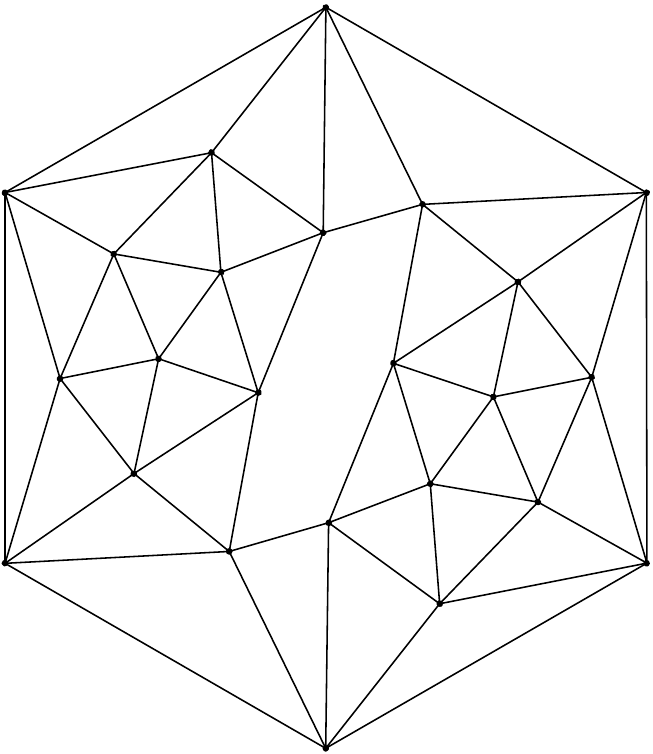}\par
$(\{3,6\},5)$: $D_2$ (iii)
\end{minipage}

\end{center}
\caption{All non-trivial $(\{a,ta\},k)$-spheres with $p_{2a}=t=2$; see Theorem \ref{Case_Pb2}}
\label{SecondTermInfSeries}
\end{figure}

\begin{figure}
\begin{center}
\begin{minipage}[b]{3.0cm}
\centering
\epsfig{height=2.4cm, file=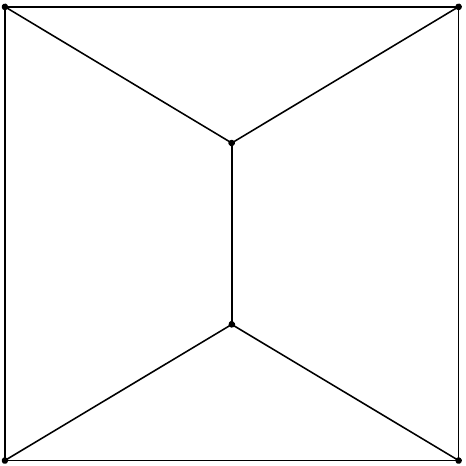}\par
$D_{3h}$, $b=4$
\end{minipage}
\begin{minipage}[b]{3.0cm}
\centering
\epsfig{height=2.4cm, file=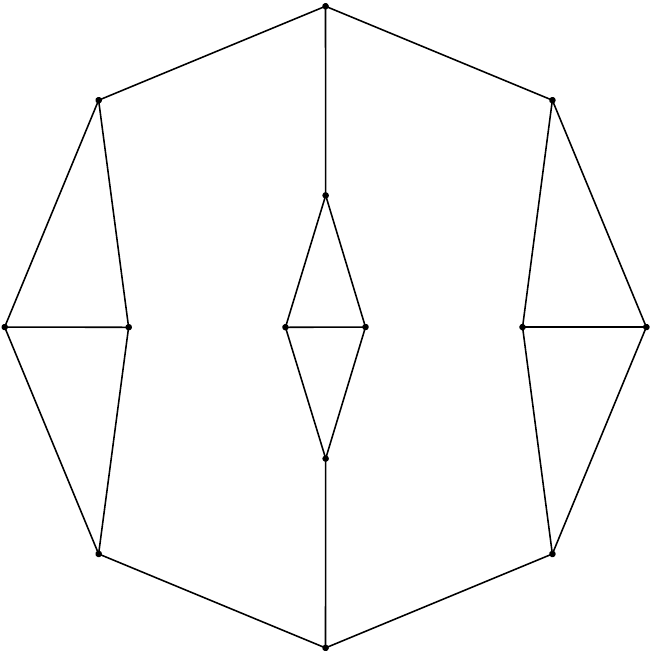}\par
$D_{3h}$, $b=2 + 6$
\end{minipage}
\begin{minipage}[b]{3.0cm}
\centering
\epsfig{height=2.4cm, file=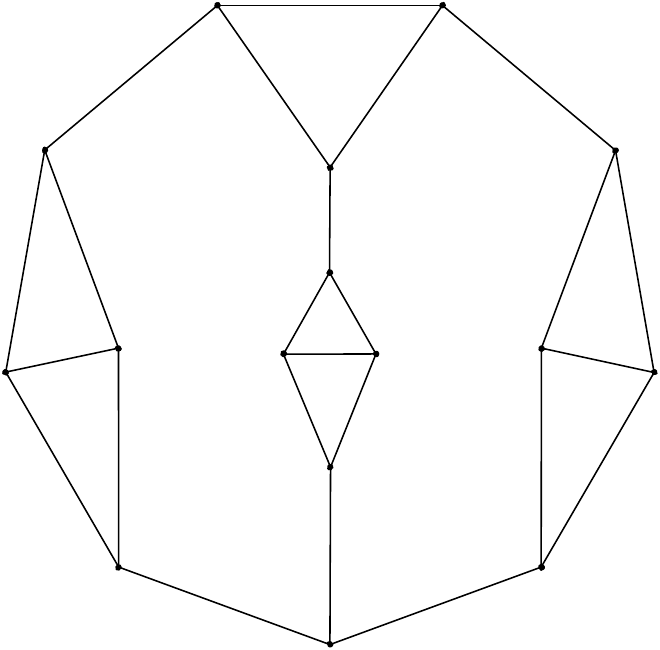}\par
$C_{3v}$, $b=3 + 6$
\end{minipage}
\begin{minipage}[b]{3.0cm}
\centering
\epsfig{height=2.4cm, file=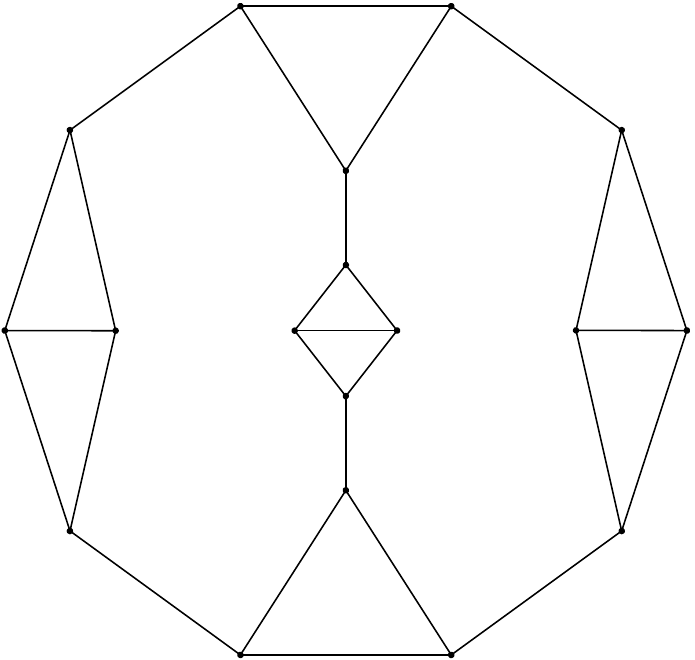}\par
$D_{3h}$, $b=4 + 6$
\end{minipage}
\end{center}
\caption{All $(\{3,b\},3)$-spheres with $p_b=3$
for $2\leq b\leq 10$}
\label{AllKnown3_10to15_3_three}
\end{figure}

\begin{figure}
\begin{center}
\begin{minipage}[b]{3.0cm}
\centering
\epsfig{height=2.4cm, file=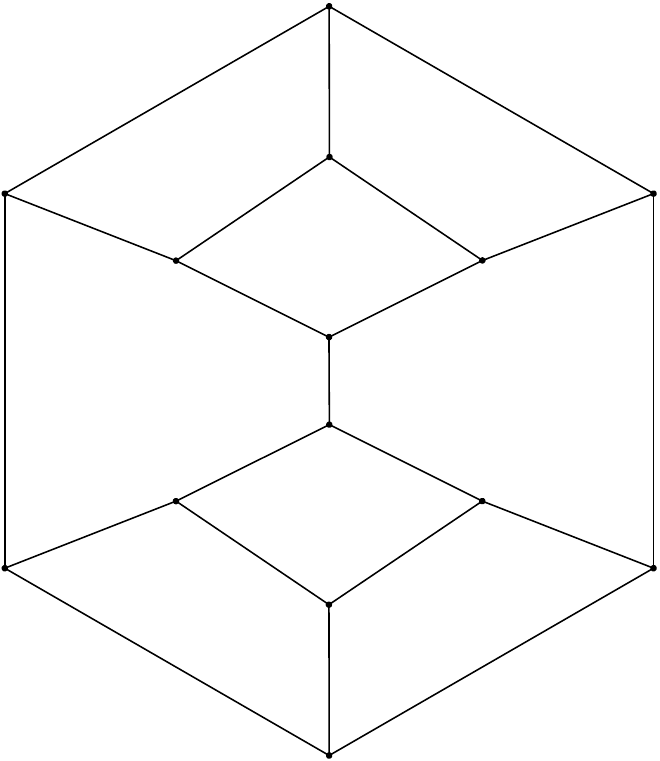}\par
$D_{3h}$, $b=6$
\end{minipage}
\begin{minipage}[b]{3.0cm}
\centering
\epsfig{height=2.8cm, file=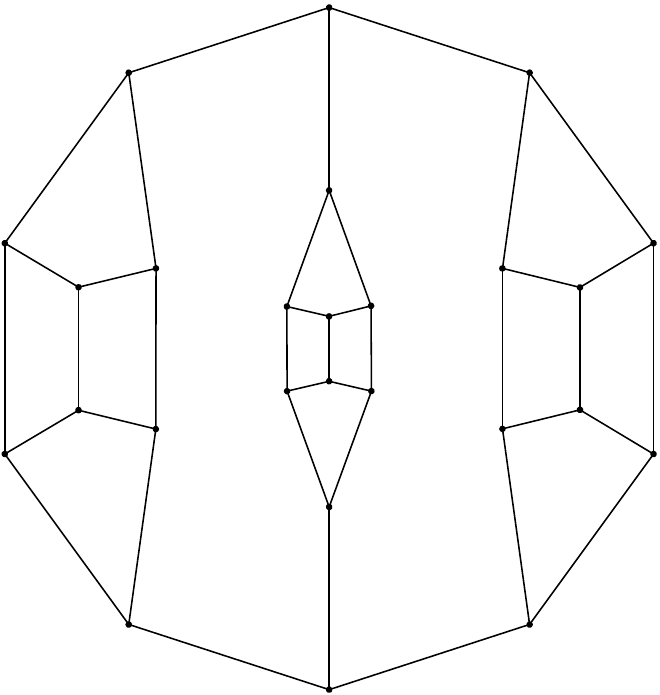}\par
$D_{3h}$, $b=2+8$
\end{minipage}
\begin{minipage}[b]{3.0cm}
\centering
\epsfig{height=2.8cm, file=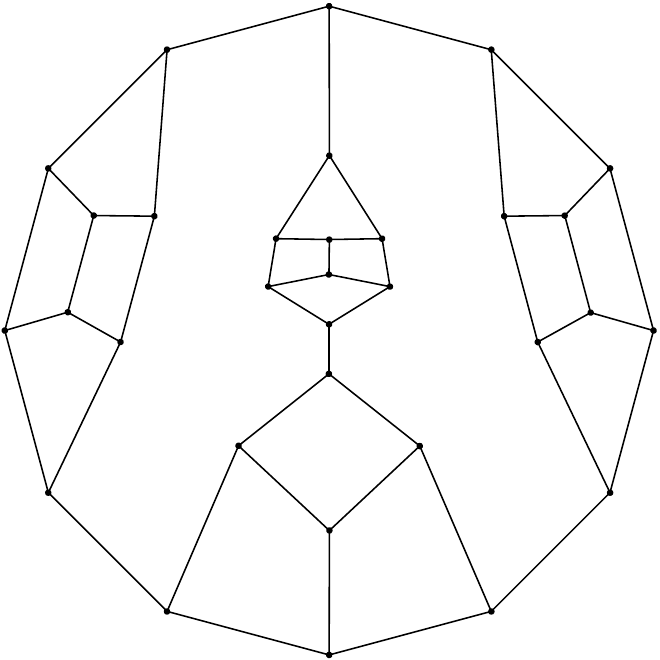}\par
$C_{3v}$, $b=4+8$
\end{minipage}
\begin{minipage}[b]{3.0cm}
\centering
\epsfig{height=2.8cm, file=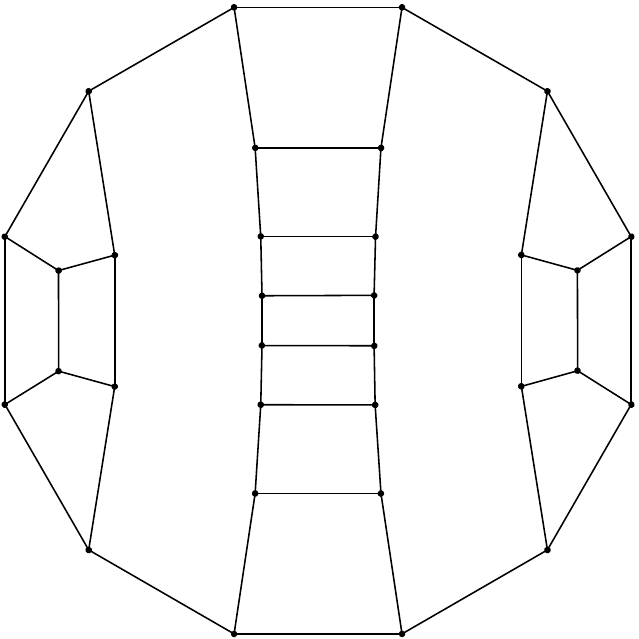}\par
$C_{2v}$, $b=4+8$
\end{minipage}
\end{center}
\caption{All $(\{4,b\},3)$-spheres with $p_b=3$
for $2\leq b\leq 12$}
\label{AllKnown4_6to20_3_three}
\end{figure}

\begin{figure}
\begin{center}
\begin{minipage}[b]{3.0cm}
\centering
\epsfig{height=2.8cm, file=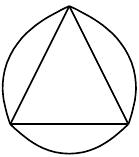}\par
$D_{3h}$, $b=2$
\end{minipage}
\begin{minipage}[b]{3.0cm}
\centering
\epsfig{height=2.8cm, file=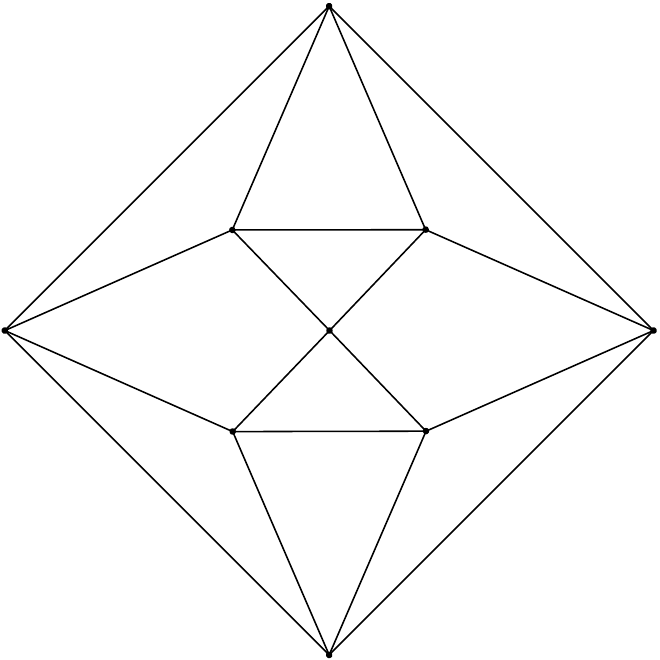}\par
$D_{3h}$, $b=4$
\end{minipage}
\begin{minipage}[b]{3.0cm}
\centering
\epsfig{height=2.8cm, file=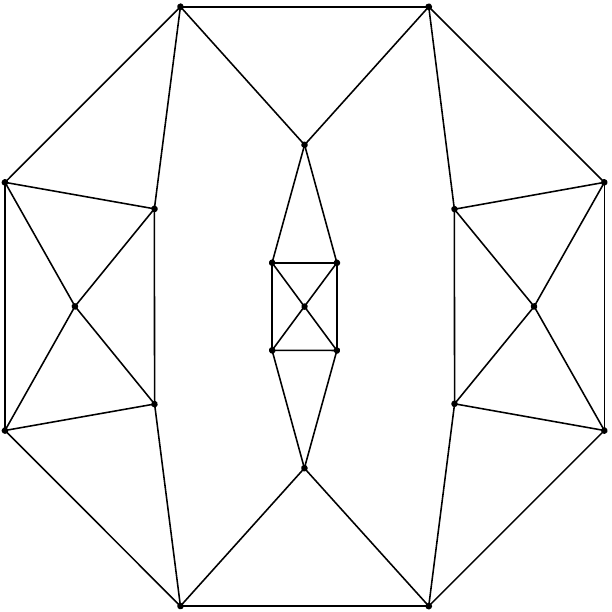}\par
$D_{3h}$, $b=2+6$
\end{minipage}
\begin{minipage}[b]{3.0cm}
\centering
\epsfig{height=2.8cm, file=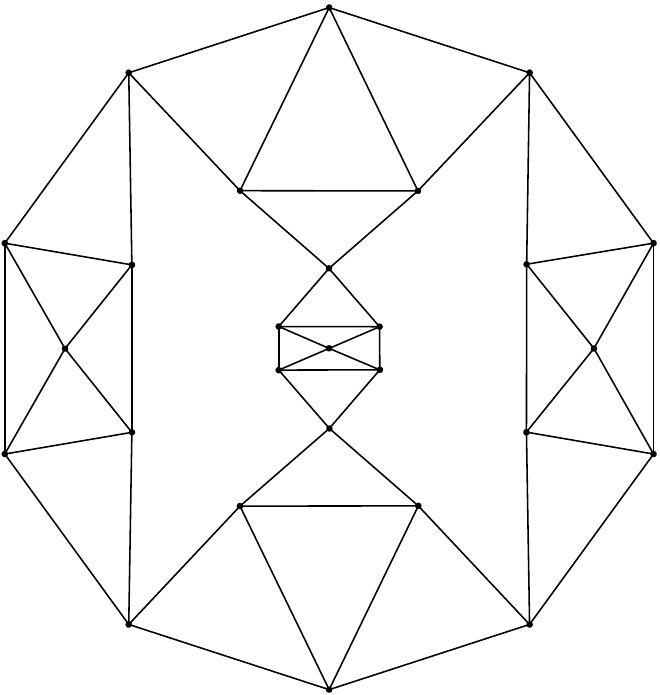}\par
$D_{3h}$, $b=4+6$
\end{minipage}
\begin{minipage}[b]{3.0cm}
\centering
\epsfig{height=2.8cm, file=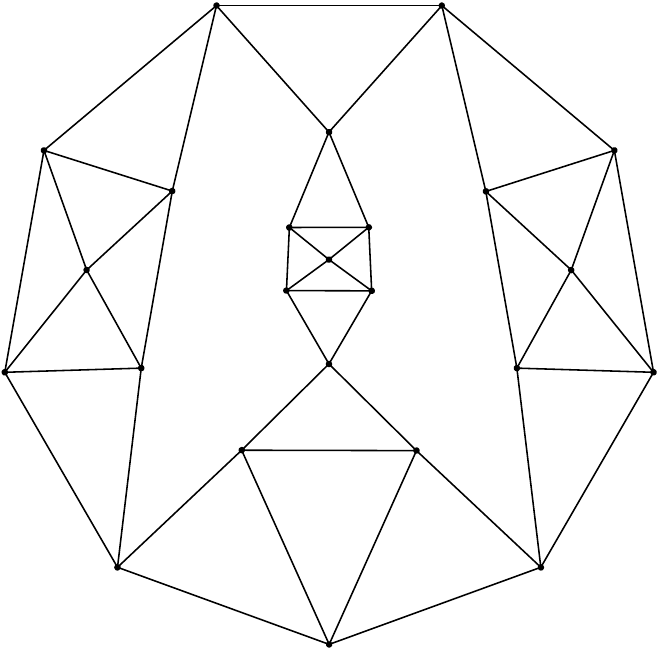}\par
$C_{3v}$, $b=3 + 6$
\end{minipage}
\begin{minipage}[b]{3.0cm}
\centering
\epsfig{height=2.8cm, file=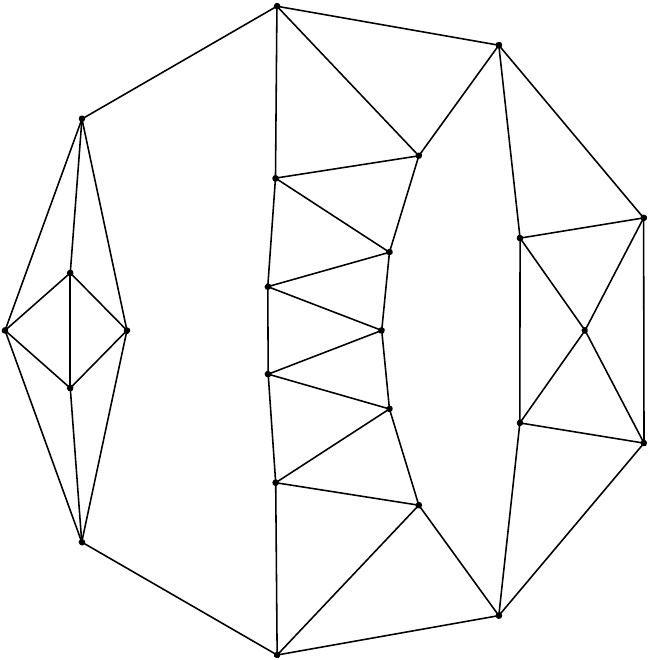}\par
$C_s$, $b=3+6$
\end{minipage}
\begin{minipage}[b]{3.0cm}
\centering
\epsfig{height=2.8cm, file=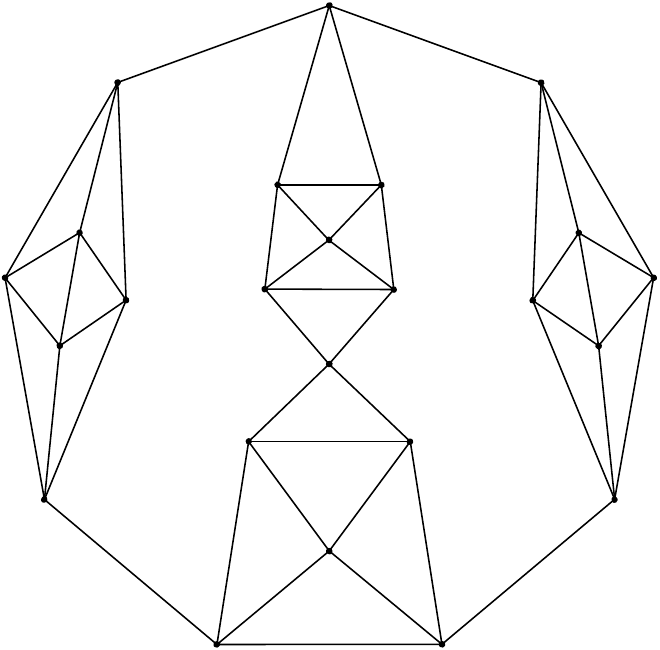}\par
$C_s$, $b=3+6$
\end{minipage}
\end{center}
\caption{All $(\{3,b\},4)$-spheres with $p_b=3$ for $2\leq b\leq 14$}
\label{AllKnown3_4to15_4_three}
\end{figure}

\begin{figure}
\begin{center}
\begin{minipage}[b]{3.0cm}
\centering
\epsfig{height=2.5cm, file=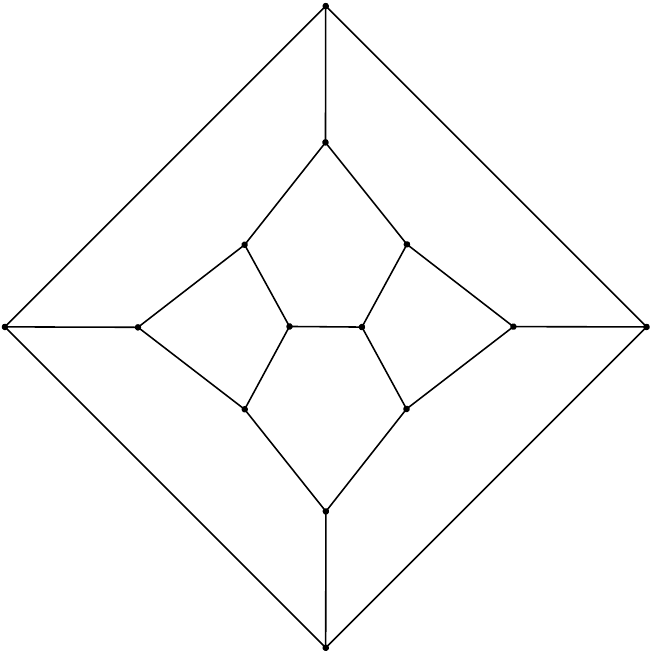}\par
$D_{3h}$, $b=4$\par
(b)
\end{minipage}
\begin{minipage}[b]{3.0cm}
\centering
\epsfig{height=2.5cm, file=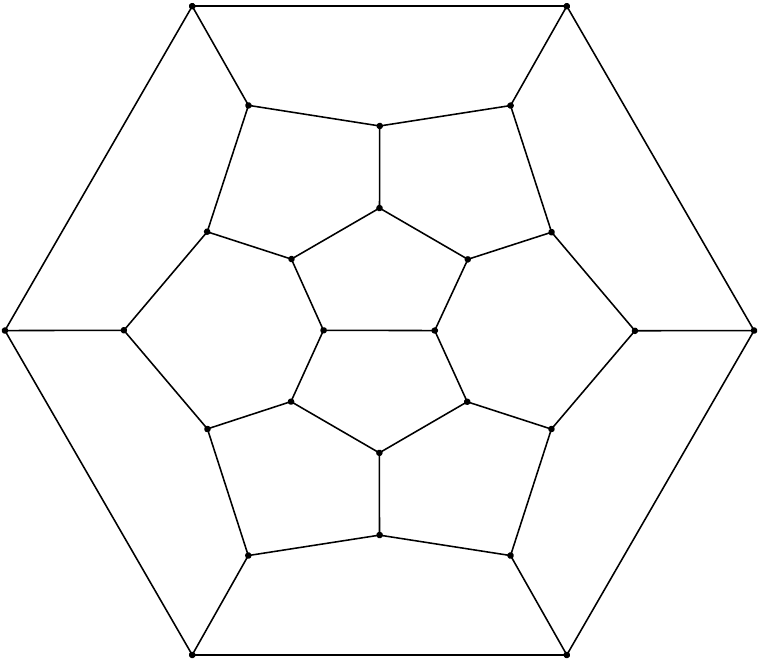}\par
$D_{3h}$, $b=6$\par
(b) 
\end{minipage}
\begin{minipage}[b]{3.0cm}
\centering
\epsfig{height=2.8cm, file=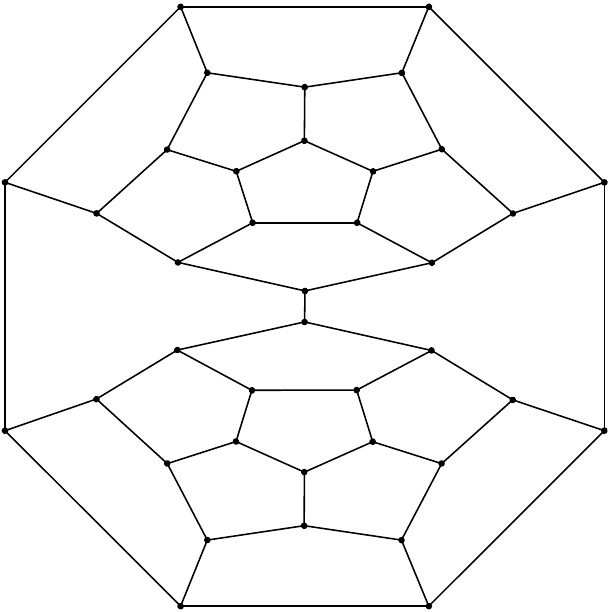}\par
$D_{3h}$, $b=8$\par
(a)
\end{minipage}
\begin{minipage}[b]{3.0cm}
\centering
\epsfig{height=2.8cm, file=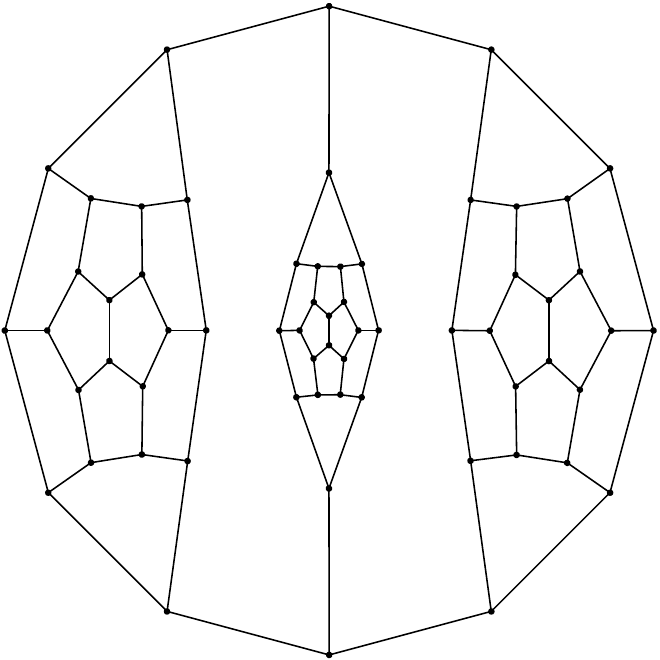}\par
$D_{3h}$, $b=2+10$\par
(a)
\end{minipage}
\begin{minipage}[b]{3.0cm}
\centering
\epsfig{height=2.8cm, file=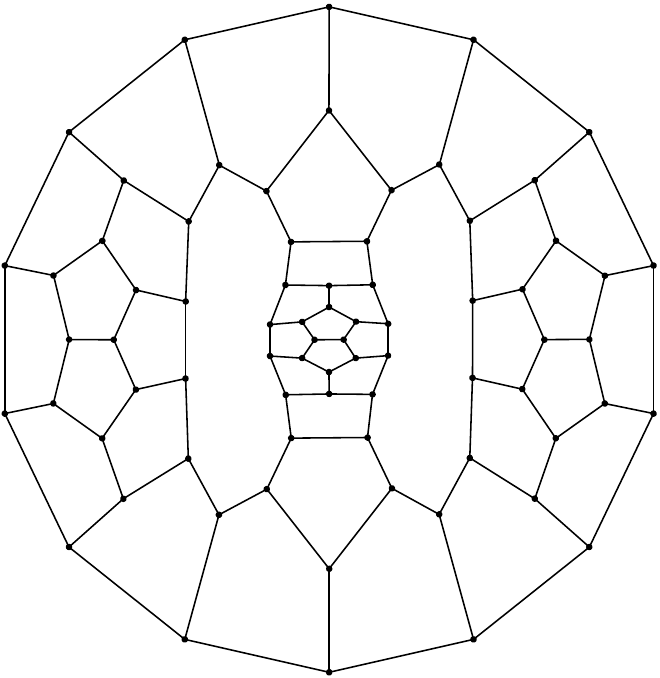}\par
$D_{3h}$, $b=4+10$\par
(b)
\end{minipage}
\begin{minipage}[b]{3.0cm}
\centering
\epsfig{height=2.8cm, file=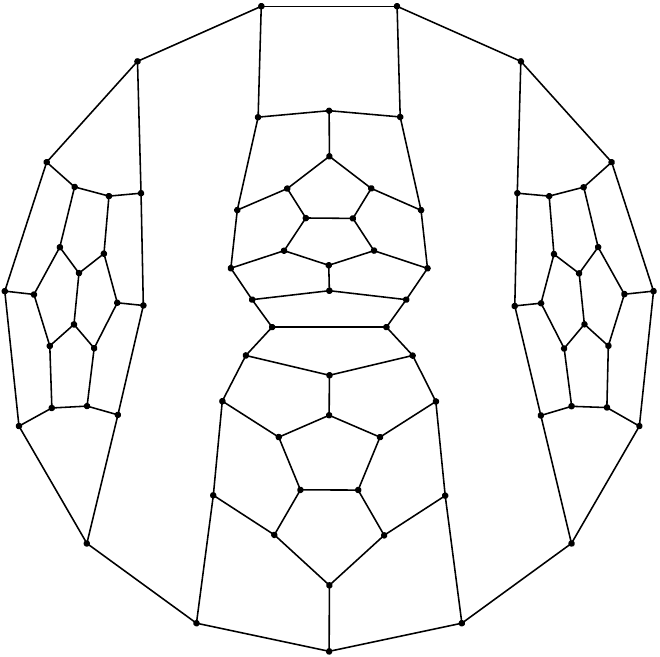}\par
$C_s$, $b=5+15$\par
(c1)
\end{minipage}
\begin{minipage}[b]{3.0cm}
\centering
\epsfig{height=2.8cm, file=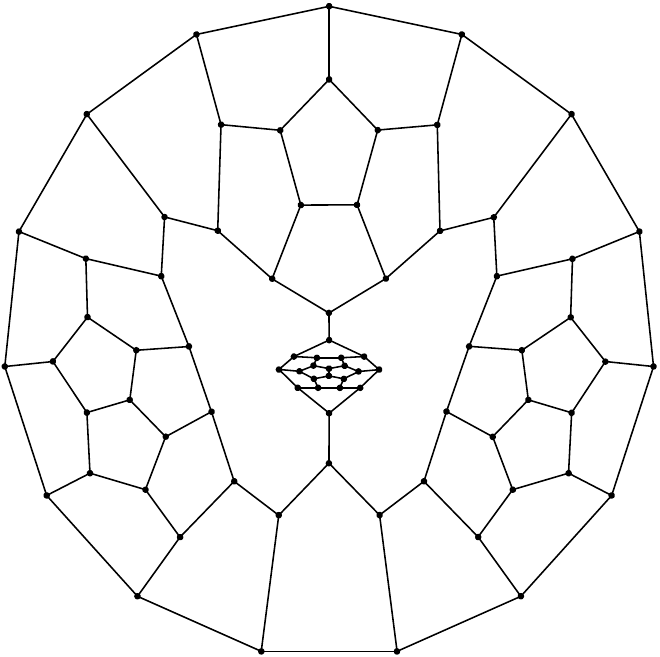}\par
$C_s$, $b=5+15$\par
(c2)
\end{minipage}
\begin{minipage}[b]{3.0cm}
\centering
\epsfig{height=2.8cm, file=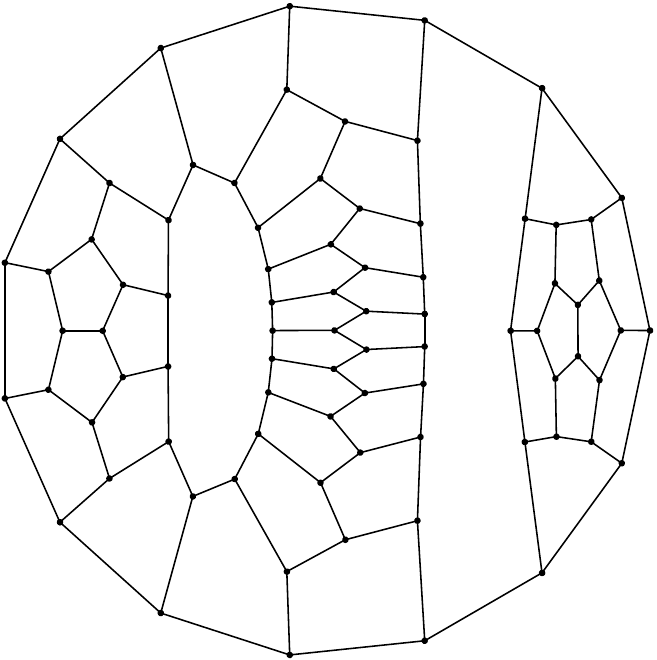}\par
$C_s$, $b=5+15$\par
(c3)
\end{minipage}
\begin{minipage}[b]{3.0cm}
\centering
\epsfig{height=2.8cm, file=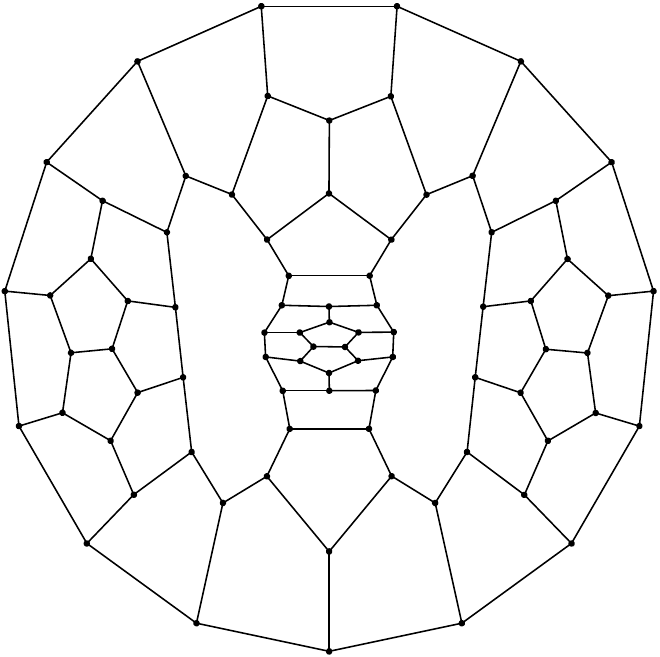}\par
$C_{3v}$, $b=5+15$\par
(b)
\end{minipage}
\begin{minipage}[b]{3.0cm}
\centering
\epsfig{height=2.8cm, file=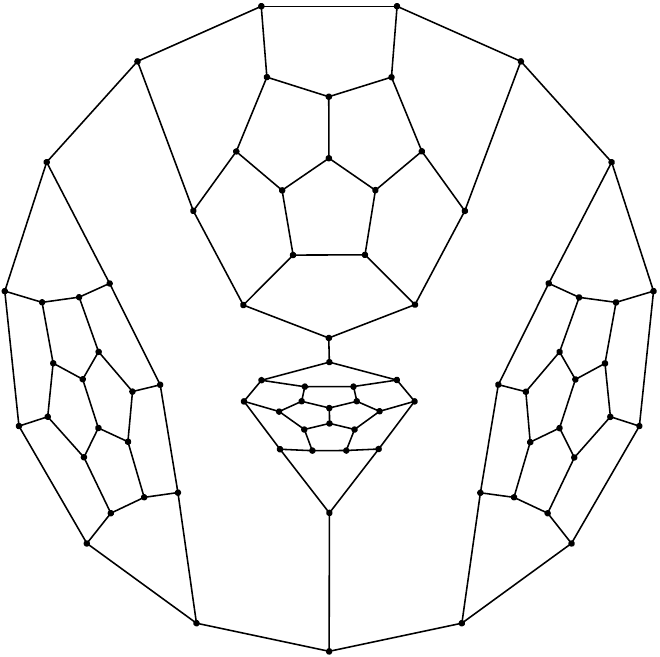}\par
$C_{3v}$, $b=5+15$\par
(a)
\end{minipage}

\end{center}
\caption{All $(\{5,b\},3)$-spheres with $p_b=3$ for $2\leq b\leq 15$}
\label{AllKnown5_15_3_three}
\end{figure}

\begin{theorem}\label{Case_PB3}
Let $b\ge 2, b\neq a$ and let 
$(a,k)=(3,3)$, $(4,3)$, $(3,4)$, $(5,3)$.
Then $(\{a,b\},k)$-spheres with $p_b=3$ exists
if and only if $b\equiv 2, a, 2a-2 \pmod {2a}$
and $b\equiv 4,6\pmod {10}$ for $a=5$.

Such spheres are unique if $b\nequiv a\pmod {2a}$
and their symmetry is $D_{3h}$.
Let $t=\lfloor{\frac{b}{2a}}\rfloor$. There are seven such spheres
with $t=0$ and $3+4+5+9$ of them for any $t\ge 1$; see corresponding
Figures \ref{AllKnown3_10to15_3_three}, \ref{AllKnown4_6to20_3_three},
\ref{AllKnown3_4to15_4_three}, \ref{AllKnown5_15_3_three} and detailed
description below.

(i) The three $(\{3,b\},3)$-sphere with $p_b=3$ and $b=2+6t$, $4+6t$,
$3+6t$ come by putting $t$ polycycles $\{3,3\}-e$ on $3$ edges of, 
respectively, $3\times K_2$, $Prism_3$ (three $4$-$4$ edges) and
Tetrahedron $\{3,3\}$.
Only for $b=2$, $4$, the graph is
$3$-connected. The symmetry is $C_{3v}$ if $b=3+6t$.

(ii) All four but one $(\{4,b\},3)$-spheres with $p_b=3$ and $b=2+8t$,
$6 + 8 t$, $4 + 8 t$ come by putting $t$ polycycles $\{4,3\}-e$
on $3$ edges of, respectively, $3\times K_2$, $14$-vertex $(\{4,6\},3)$-sphere 
(three $6$-$6$ edges) and Cube $\{4,3\}$ ($3$ incident edges).
The remaining sphere is $Sq_{8t-1}$ with its two end edges connected on each of
two sides by a chain of $t$ polycycles $\{4,3\} - e$.
This graph has symmetry $C_{2v}$, while
other graph coming from Cube has symmetry $C_{3v}$.
Only for $b=6$, the graph is $3$-connected.

(iii) All five but two $(\{3,b\},4)$-spheres with $p_b=3$ and $b=2 + 6t$, 
$4+6t$, $3+6t$ come when replacing by $t$ vertex-split $\{3,4\}$'s, 
$3$ vertices of, respectively, edge-doubled triangle,
$9$-vertex $(\{3,4\},4)$-sphere ($3$ vertices common to two $4$-gons)
and Octahedron $\{3,4\}$ ($3$ vertices of a triangle).
First of remaining spheres consists of $Tr_{6t+5}$ connected to its other
end on one side by a chain of $t$ vertex-split $\{3,4\}$'s
and on the other side by a chain of $t$ polycycles $\{3,4\}-e$.
Second remaining sphere consists
of a vertex and $\{3,4\}-P_3$
connected by two chains of $t$ polycycles $\{3,4\}-e$
and one chain of $t$ vertex-split $\{3,4\}$'s.
Those graphs have symmetry $C_s$, while other graph coming from Octahedron
has symmetry $C_{3v}$.
Only for $b=2$ and $b=4$, the graphs are
$3$-connected.

There are nine $(\{5,b\},3)$-spheres with $p_b=3$ 
for each $t=\lfloor{\frac{b}{10}}\rfloor \ge 1$.
\begin{enumerate}
\item[(a)] Three spheres with $b=2 + 10t$, $b=8 + 10t$, $b=5 + 10t$
come by putting $t$ polycycles $\{5,3\}-e$ on $3$ edges of
$3\times K_2$, $22$-vertex $(\{5,8\},3)$-sphere (three $8$-$8$ edges)
and Dodecahedron $\{5,3\}$ ($3$ incident edges).
Symmetry is $C_{3v}$ in last case and $D_{3h}$, otherwise.
Only for $b=8$, the graph is $3$-connected.

\item[(b)] Three $3$-connected spheres 
with $b=4 + 10t$, $6+10t$, $5+10t$ come
by putting $3$ chains of $t$ edge-split $\{5,3\}$'s
between two polycycles, $\{Pen_3, Pen_3\}$, $\{C_3,C_3\}$ and $\{Pen_3, C_3\}$.
Symmetry is $C_{3v}$ in last case and $D_{3h}$, otherwise.

\item[(c)] Three $2$-connected $(\{5,5+10t\},3)$-spheres of symmetry 
$C_s$ come from, respectively:
\begin{enumerate}
\item[c1)] a chain $B_2$ + $t$ times edge-split 
$\{5,3\}$ connected on both ends by $t$ times $\{5,3\}-e$;
\item[c2)] $C_2 + t (\{5,3\}- e) + Pen_1$ connected on both ends by $t$ times 
edge-split $\{5,3\}$;
\item[c3)] $Pen_{10 t + 9}$ connected on one side by $t$ 
times edge-split $\{5,3\}$ and on the other one by $t$ times $\{5,3\}-e$.
\end{enumerate}

\end{enumerate}

\end{theorem}
\proof It is not possible to have 
an elementary polycycle
separating three
$b$-gonal faces. Hence, the decomposition graph has three faces and is either
formed of two vertices of degree $3$ connected by chains of vertices
of degree $2$,
or a vertex of degree $4$ connected by two chains of vertices
of degree $2$ on each side.
An examination of the possibilities along the same lines gives the above
result. \qed

Note that all $(\{a,b\},k)$-spheres with $p_b=3$ and symmetry $\neq C_s$,
are $bR_j$ (i.e., each $b$-gon has exactly $j$ edges of adjacency with
$b$-gons); see {\em face-regularity} in
Section \ref{SectionIcosahedrite}).
$j=2\left\lceil {\frac{b}{2a}} \right\rceil$ for $(a,k)=(3,3)$, $(4,3)$
and $j=0$ for $(a,k)=(3,4)$, $(3,5)$.
For $(a,k)=(5,3)$, we have
$j=0$ or $2\left\lceil {\frac{b}{2a}} \right\rceil$.

\begin{figure}
\begin{center}
\begin{minipage}[b]{3.0cm}
\centering
\epsfig{height=3.0cm, file=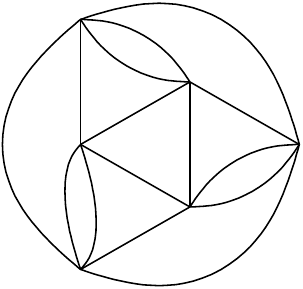}\par
$D_3$, $b=2$ (a)
\end{minipage}
\begin{minipage}[b]{3.0cm}
\centering
\epsfig{height=3.0cm, file=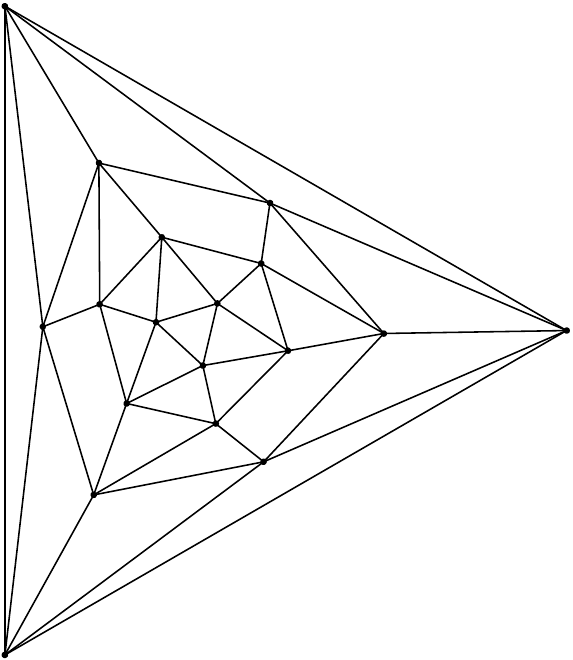}\par
$D_3$, $b=4$ (a)
\end{minipage}
\begin{minipage}[b]{3.0cm}
\centering
\epsfig{height=3.0cm, file=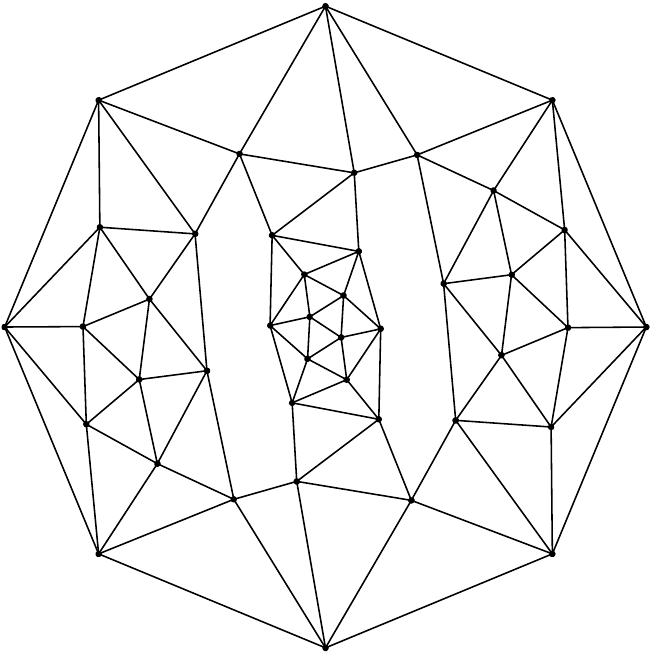}\par
$D_3$, $b=2+6$ (a)
\end{minipage}
\begin{minipage}[b]{3.0cm}
\centering
\epsfig{height=3.0cm, file=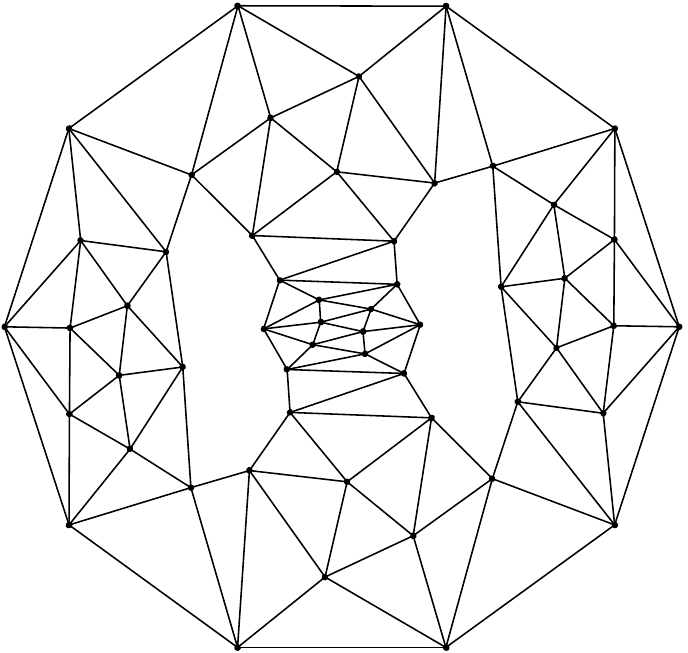}\par
$D_3$, $b=4+6$ (a)
\end{minipage}
\begin{minipage}[b]{4.1cm}
\centering
\epsfig{height=4.0cm, file=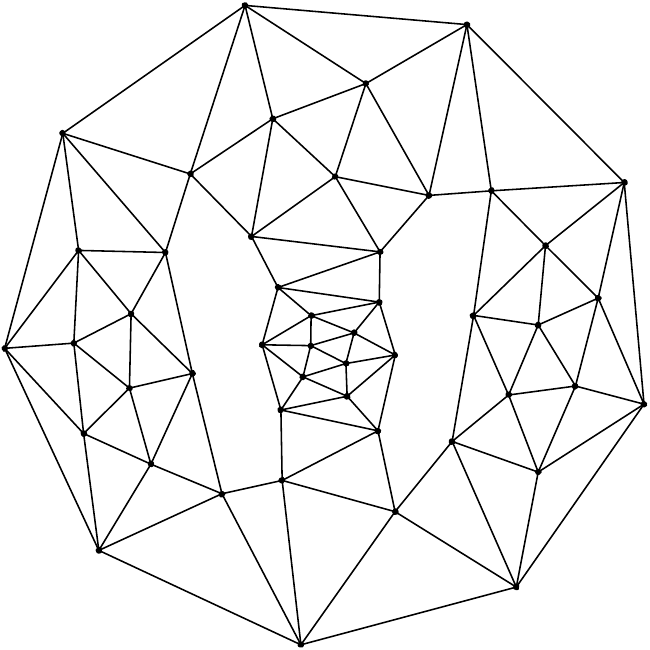}\par
$C_3$, $b=3+6$ (a)
\end{minipage}
\begin{minipage}[b]{4.1cm}
\centering
\epsfig{height=4.0cm, file=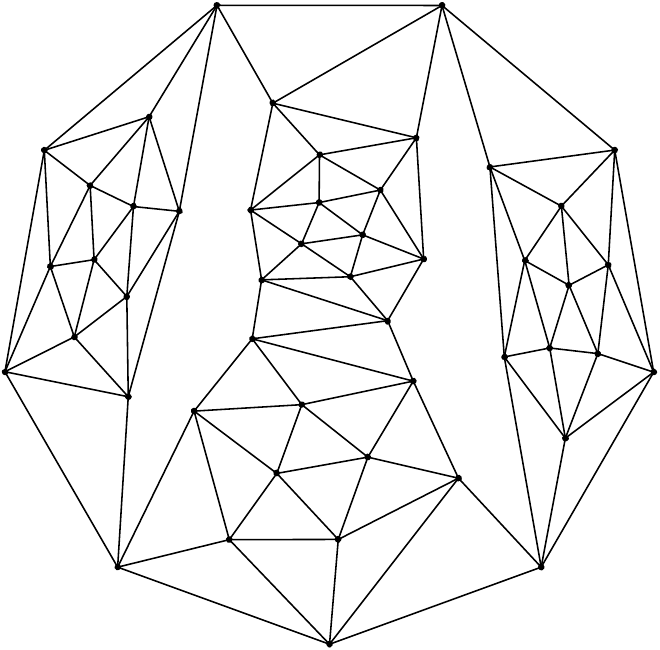}\par
$C_1$, $b=3+6$ ($b_1$)
\end{minipage}
\begin{minipage}[b]{4.1cm}
\centering
\epsfig{height=4.0cm, file=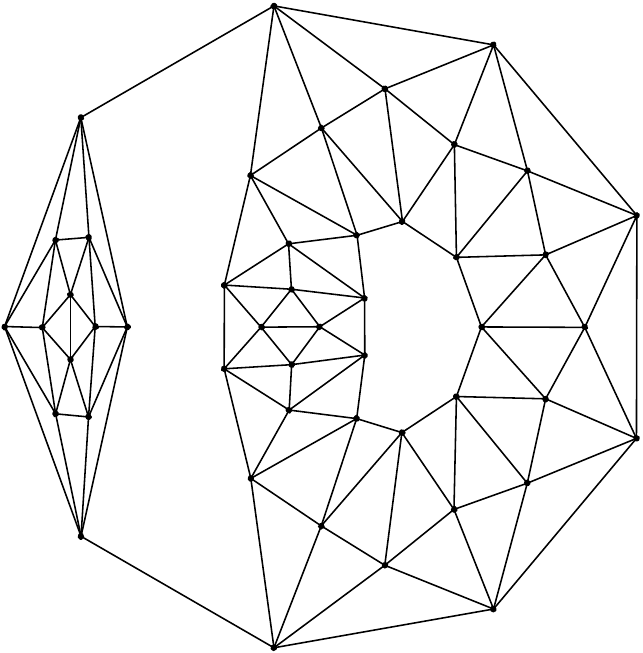}\par
$C_s$, $b=3+6$ ($b_2$)
\end{minipage}
\begin{minipage}[b]{4.1cm}
\centering
\epsfig{height=4.0cm, file=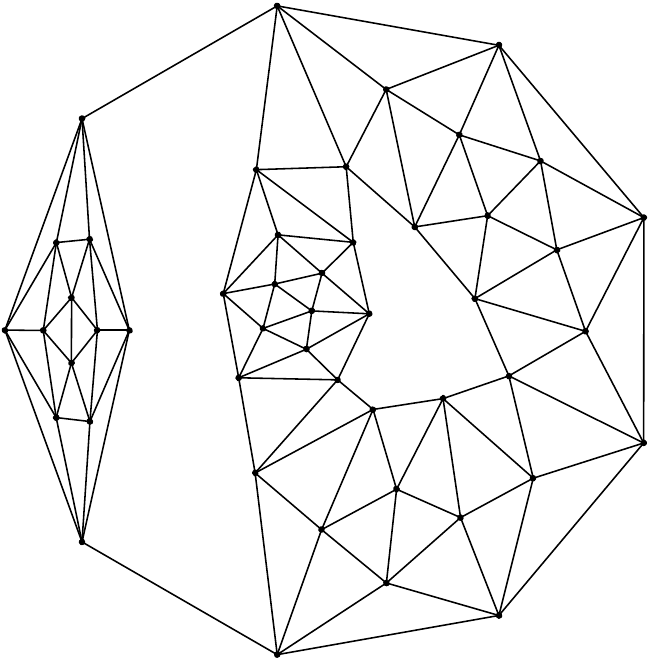}\par
$C_1$, $b=3+6$ ($b_3$)
\end{minipage}
\begin{minipage}[b]{4.1cm}
\centering
\epsfig{height=4.0cm, file=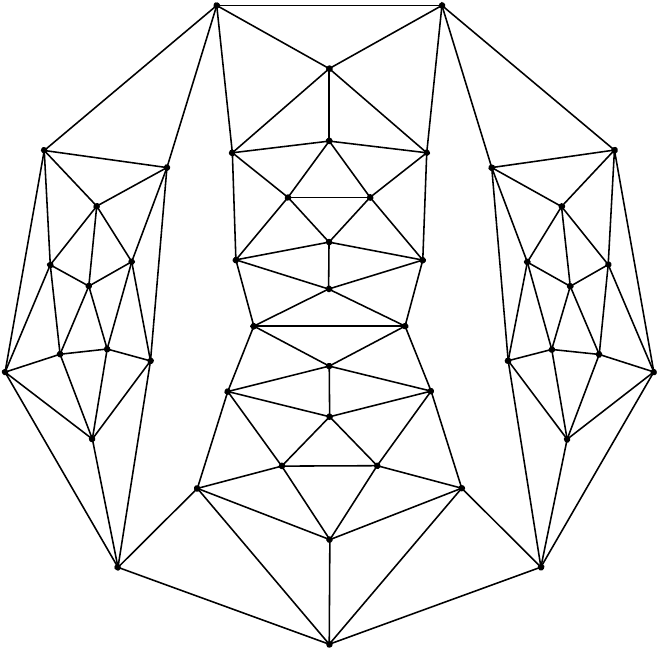}\par
$C_s$, $b=3+6$ ($b_4$)
\end{minipage}
\begin{minipage}[b]{4.1cm}
\centering
\epsfig{height=4.0cm, file=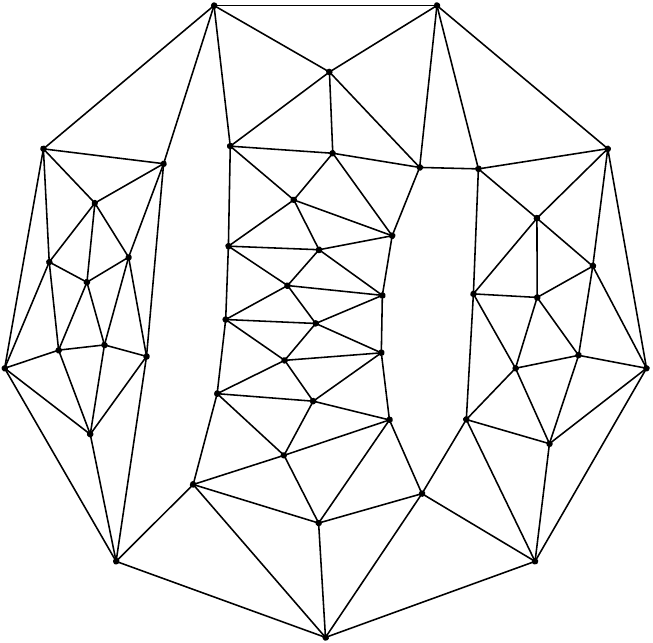}\par
$C_1$, $b=3+6$ ($b_5$)
\end{minipage}
\begin{minipage}[b]{4.1cm}
\centering
\epsfig{height=4.0cm, file=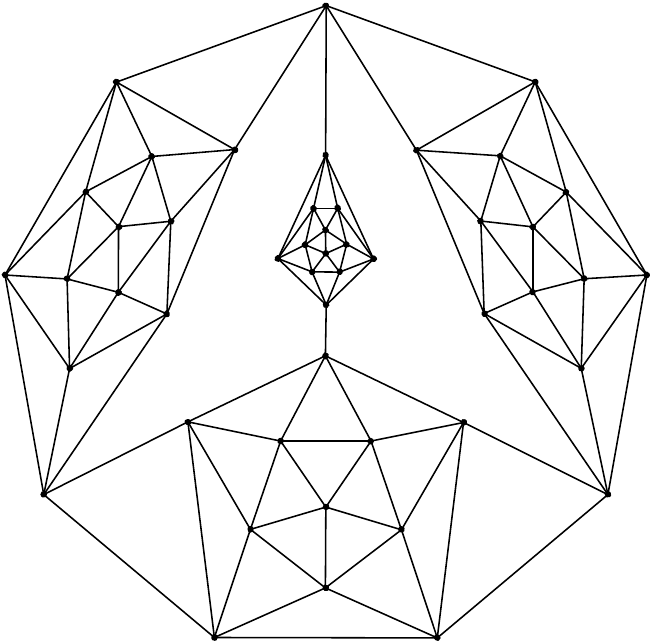}\par
$C_s$, $b=3+6$ ($c_1$)
\end{minipage}
\begin{minipage}[b]{4.1cm}
\centering
\epsfig{height=4.0cm, file=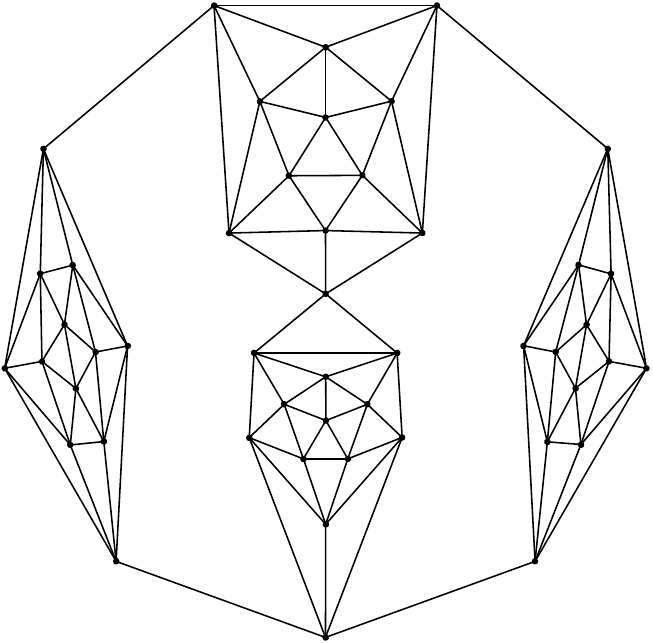}\par
$C_s$, $b=3+6$ ($c_2$)
\end{minipage}
\begin{minipage}[b]{4.1cm}
\centering
\epsfig{height=4.0cm, file=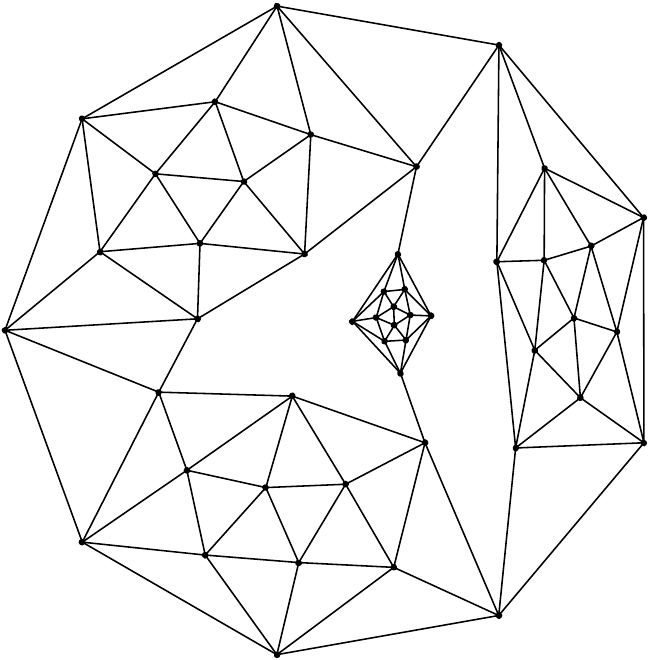}\par
$C_1$, $b=3+6$ ($c_3$)
\end{minipage}

\end{center}
\caption{All known $(\{3,b\},5)$-spheres with $p_b=3$ and $2\leq b\leq 10$ (part 1)}
\label{FirstNonTrivial_53_I}
\end{figure}

\begin{figure}
\begin{center}
\begin{minipage}[b]{4.1cm}
\centering
\epsfig{height=4.0cm, file=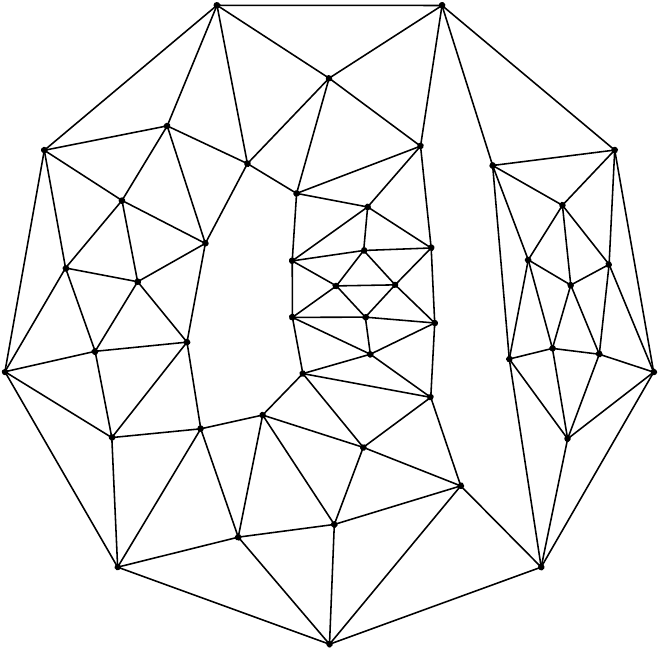}\par
$C_s$, $b=3+6$ ($d_1$)
\end{minipage}
\begin{minipage}[b]{4.1cm}
\centering
\epsfig{height=4.0cm, file=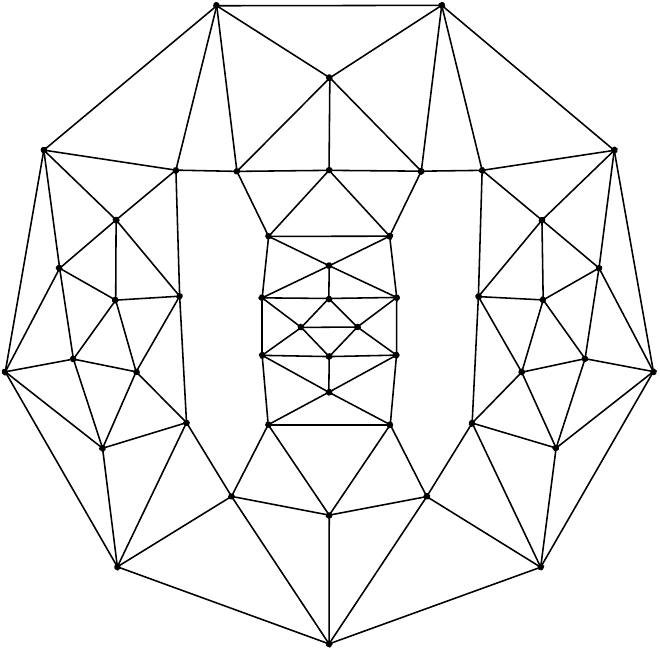}\par
$C_s$, $b=3+6$ ($d_2$)
\end{minipage}
\begin{minipage}[b]{4.1cm}
\centering
\epsfig{height=4.0cm, file=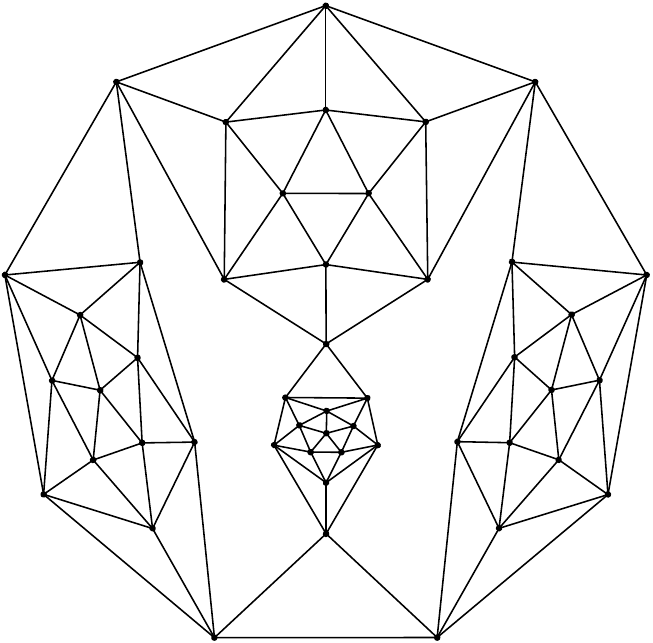}\par
$C_{3v}$, $b=3+6$ ($e_1$)
\end{minipage}
\begin{minipage}[b]{4.1cm}
\centering
\epsfig{height=4.0cm, file=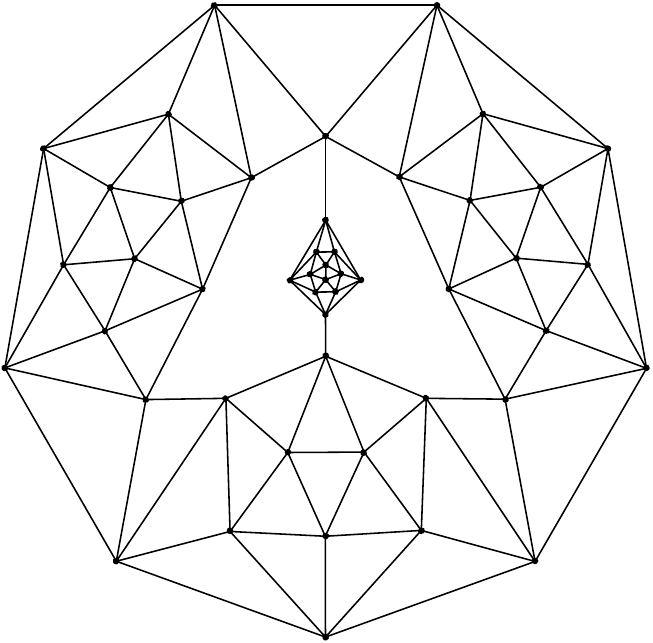}\par
$C_s$, $b=3+6$ ($e_2$)
\end{minipage}
\begin{minipage}[b]{4.1cm}
\centering
\epsfig{height=4.0cm, file=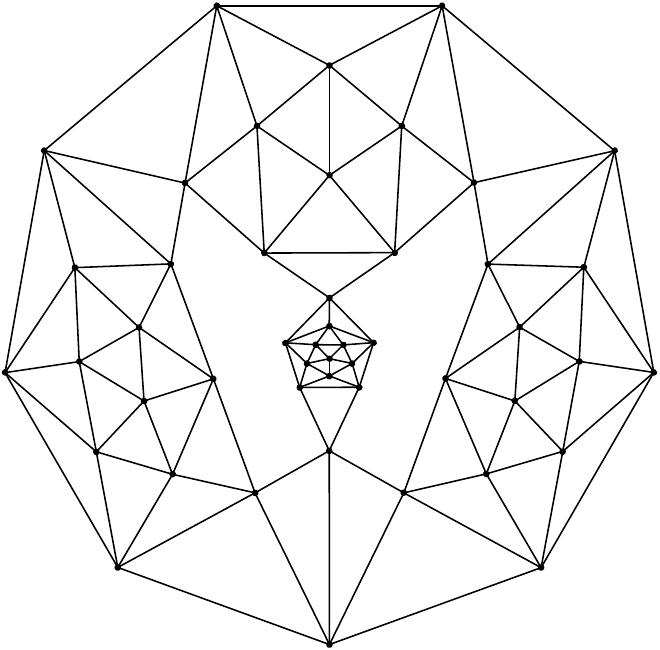}\par
$C_s$, $b=3+6$ ($e_3$)
\end{minipage}
\begin{minipage}[b]{4.1cm}
\centering
\epsfig{height=4.0cm, file=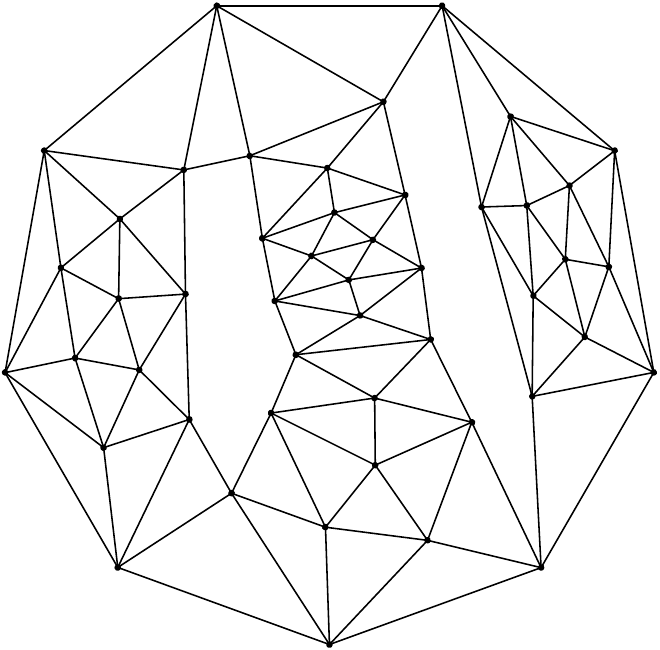}\par
$C_1$, $b=3+6$ ($e_4$)
\end{minipage}

\end{center}
\caption{All known $(\{3,b\},5)$-spheres with $p_b=3$ and $2\leq b\leq 10$ (part 2)}
\label{FirstNonTrivial_53_II}
\end{figure}

In case $(a,k)=(3,5)$, we have $17$ infinite series of spheres but no
proof that the list is complete.
See below the list of $(\{3,b\},5)$-spheres obtained
and in Figure \ref{FirstNonTrivial_53_I}, \ref{FirstNonTrivial_53_II}
their pictures for small $b$. All but (a) have $b=3+6t$.
All but (a) and ($d_2$) are only $2$-connected.
By $R_e$, $V_{sp}$, $E_{sp}$, and $F_{sp}$ we denote a chain of $t$ polycycles $A$
with $A$ being $\{3,5\}-e$, vertex-split $\{3,5\}$, edge-split $\{3,5\}$
and face-split $\{3,5\}$, respectively.

\begin{enumerate}
\item[(a)] Three spheres with $b=2 + 6t$, $4+6t$, $3+6t$
obtained by putting three $E_{sp}$
between two polycycles, 
$\{Tr_1+3Tr_1, Tr_1+3Tr_1\}$, $\{c_4 + 3Tr_1, c_4 + 3Tr_1\}$ and
$\{c_4 + 3Tr_1, Tr_1 + 3Tr_1\}$.
Symmetry is $C_{3}$ in last case and $D_{3}$, otherwise.

\item[($b_1$)] $C_1$: $c_1+Tr_1+F_{sp}$ with ends connected by  $R_{e}$ and $V_{sp}$.
\item[($b_2$)] $C_s$:  $Sun_{6t+5}$ with ends connected by  $R_{e}$ and $E_{sp}$.
\item[($b_3$)] $C_1$: $b_3+Tr_1+E_{sp}$ with ends connected by  $R_{e}$ and $F_{sp}$.
\item[($b_4$)] $C_s$: $b_3+2Tr_1+E_{sp}$  with ends connected by  $V_{sp}$ and $E_{sp}$.
\item[($b_5$)] $C_1$: $Sun_{4+6t}+3Tr_1$  with ends connected by  $V_{sp}$ and $F_{sp}$.

\item[($c_1$)] $C_{s}$: a vertex and $(\{3,5\}-v)+2Tr_1$ connected by  $R_{e}$ and two $V_{sp}$'s.    
\item[($c_2$)] $C_s$: a vertex and $a_3$  connected by two $R_{e}$'s and one $V_{sp}$.
\item[($c_3$)] $C_1$: a vertex and $c_2+2Tr_1$  connected by two $V_{sp}$'s and one $F_{sp}$.

\item[($d_1$)] $C_s$: $Sun_1$ and $Sun_3+Tr_1$  connected by $V_{sp}$ and two $E_{sp}$'s.
\item[($d_2$)] $C_s$: $Sun_2$ and $Sun_1+2Tr_1$ and connected by $E_{sp}$ and two $F_{sp}$'s.

\item[($e_1$)] $C_{3v}$:  $Tr_1$ and $a_4$ connected by  three $V_{sp}$'s.
\item[($e_2$)] $C_s$: $Tr_1$ and $c_2$ connected by $R_{e}$ and two $F_{sp}$'s.
\item[($e_3$)] $C_s$: $2Tr_1$ and $c_3+3Tr_1$ connected by $V_{sp}$ and two $F_{sp}$'s.
\item[($e_4$)] $C_1$: $3Tr_1$ and $b_4+Tr_1$  connected by $V_{sp}$, $E_{sp}$ and $F_{sp}$.

\end{enumerate}
Note that all $15$ orbits of $(\{3,3+6t\},5)$-spheres
constructed were built from the $15$ orbits of triples
of triangles in Icosahedron.
Note also that ${20\choose 3} = \sum_G\frac{120}{|G|}=
5\times 120 + 8\times 60 + 40 + 20$, since
there are  $5$, $8$, $1$, $1$ cases with symmetry
$C_1$, $C_s$, $C_3$, $C_{3v}$, respectively.
There may be other spheres in this case and for other values of $b$.

\section{Icosahedrites}\label{SectionIcosahedrite}
We call {\em icosahedrite} the $(\{3,4\},5)$-spheres.
They are in a sense the simplest $5$-valent plane graphs.
Clearly, for them it holds  $p_3=20+2p_4$ and $v=12+2p_4$.
Note that all $(\{a,3\},5)$-spheres with $a<3$ are: 
$(\{1,3\},5)$-sphere with $(p_1,p_3;v)=(2,6;4)$ and four
$(\{2,3\},5)$-spheres with $(p_2,p_3;v)=(4,4;4)$ (two), $(3,8;6)$, $(2,12;8)$.
Remaining $(\{1,2,3\},5)$-spheres should have 
$(p_1,p_2,p_3;v)=(1,3,1;2)$, $(2,1,2;2)$, $(1,2,5;4)$ or $(1,1,9;6)$; 
only 2nd and 3rd exist.

\begin{figure}
\begin{center}
\begin{minipage}[b]{3.0cm}
\centering
\epsfig{height=2.8cm, file=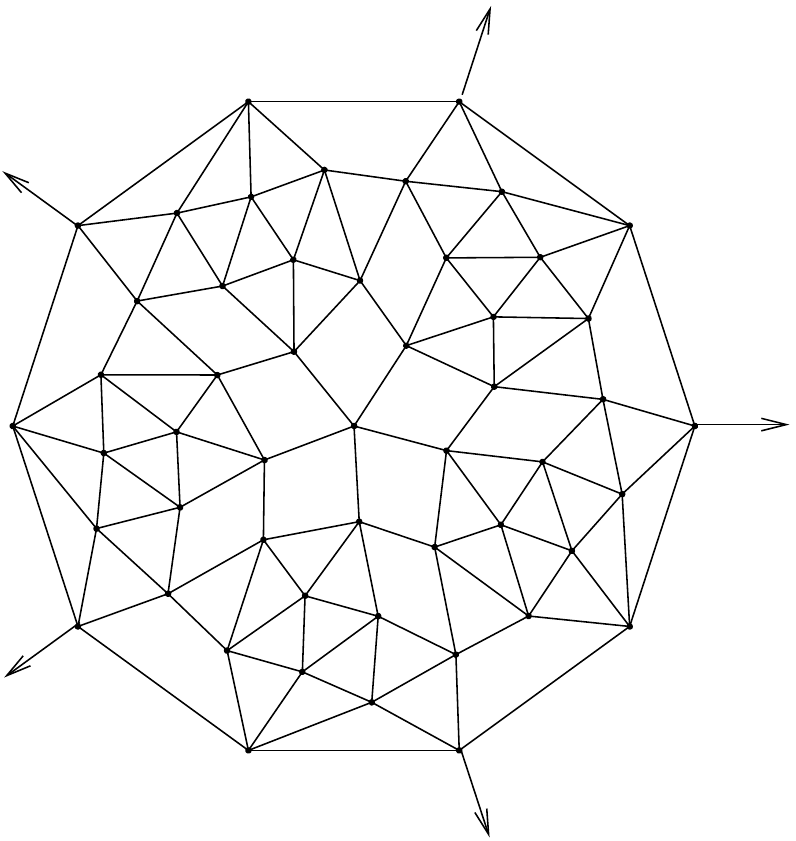}\par
$C_5$, 52
\end{minipage}
\begin{minipage}[b]{3.0cm}
\centering
\epsfig{height=2.8cm, file=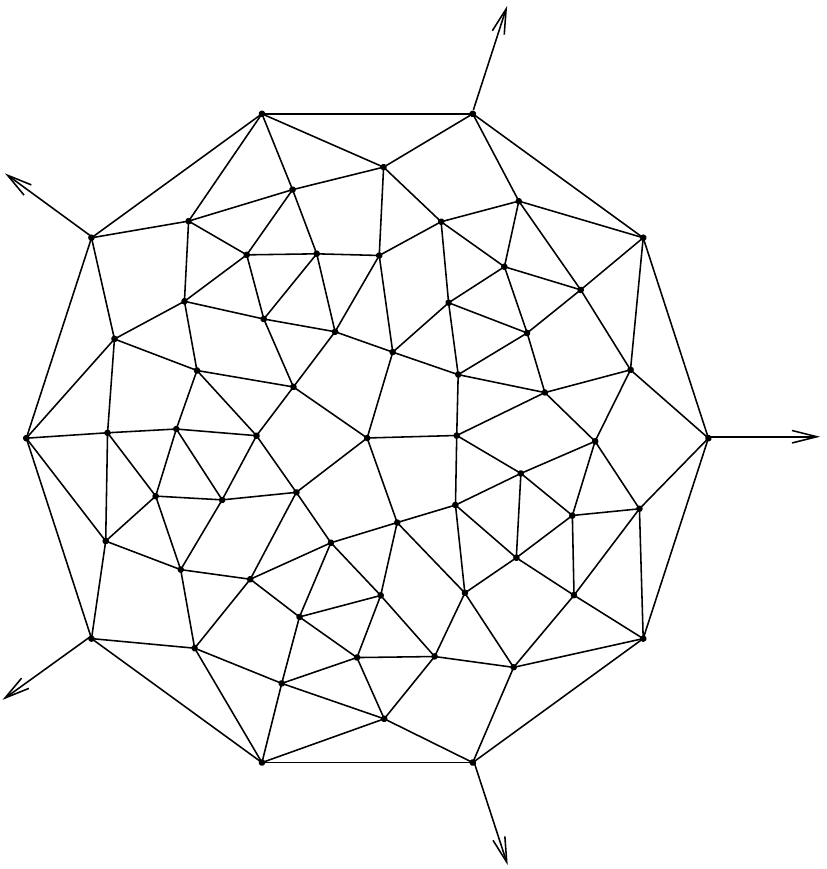}\par
$C_{5h}$, 62
\end{minipage}
\begin{minipage}[b]{3.0cm}
\centering
\epsfig{height=2.8cm, file=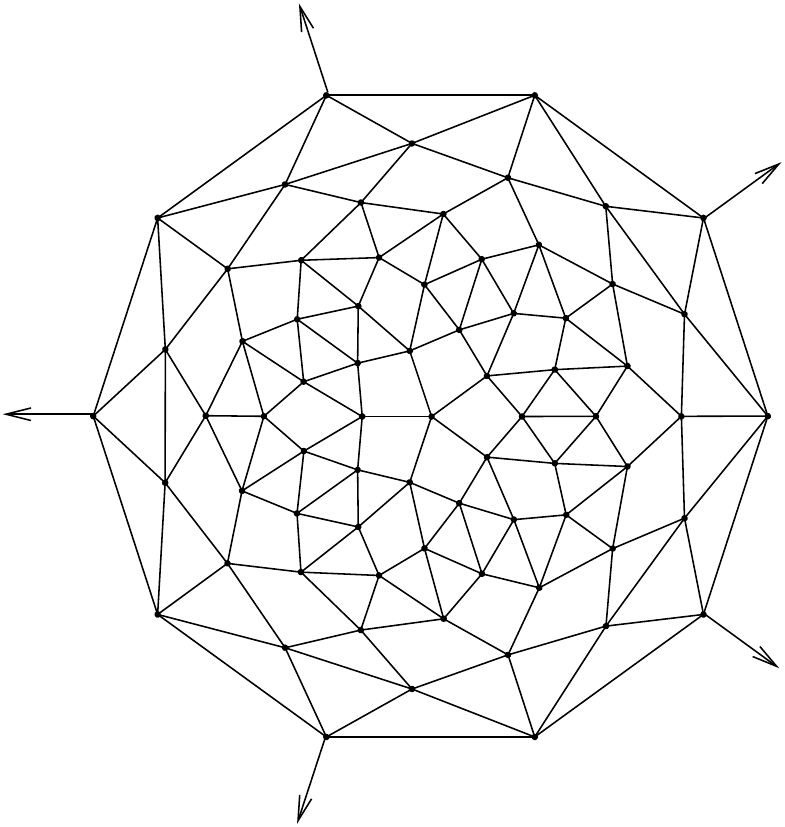}\par
$C_{5v}$, 72
\end{minipage}
\begin{minipage}[b]{3.0cm}
\centering
\epsfig{height=2.8cm, file=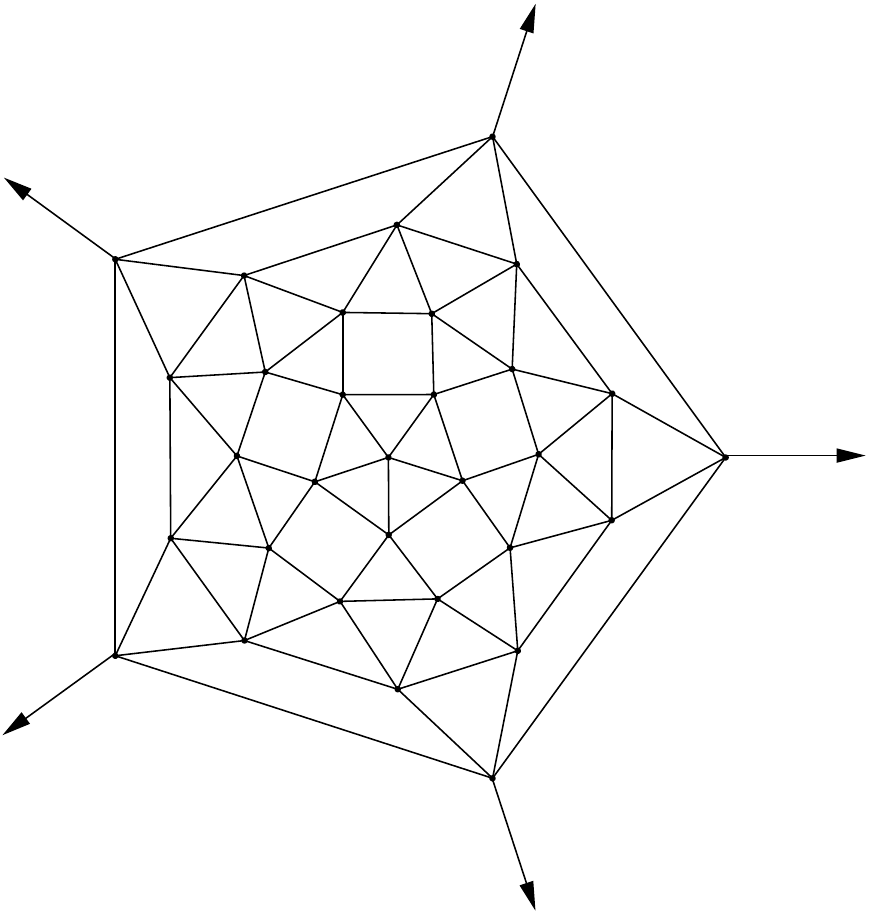}\par
$D_5$, 32, $4R_0$
\end{minipage}
\begin{minipage}[b]{3.0cm}
\centering
\epsfig{height=2.8cm, file=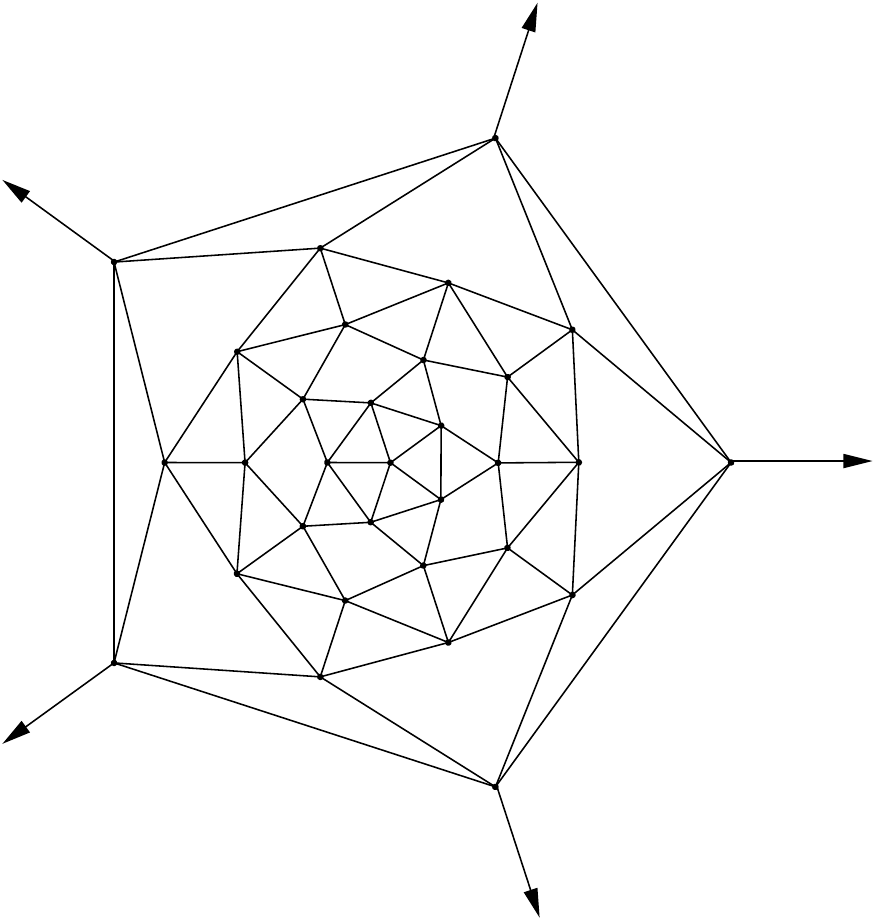}\par
$D_{5d}$, 32, $4R_0$
\end{minipage}
\begin{minipage}[b]{3.0cm}
\centering
\epsfig{height=2.8cm, file=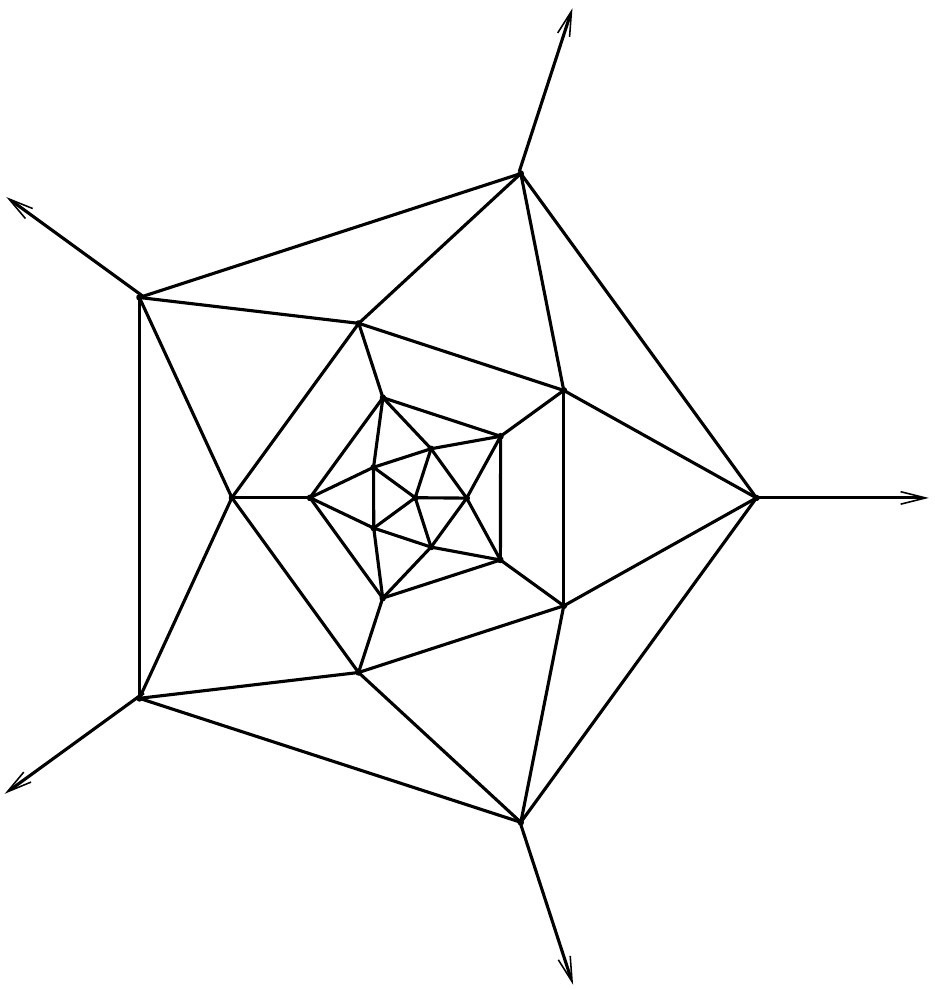}\par
$D_{5h}$, 22, $4R_2$
\end{minipage}
\begin{minipage}[b]{3.0cm}
\centering
\epsfig{height=2.8cm, file=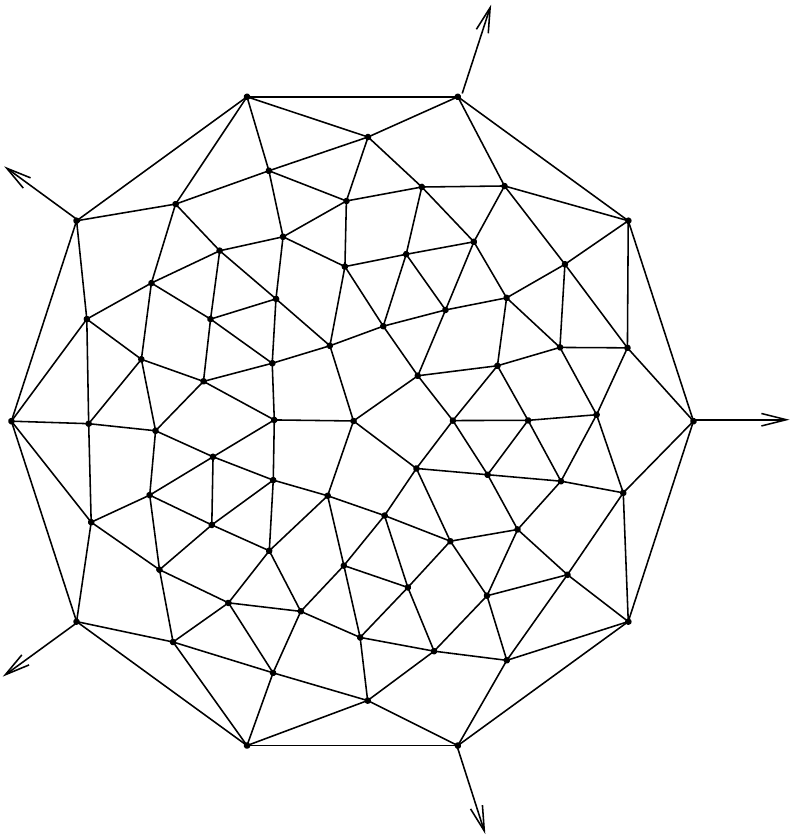}\par
$S_{10}$, 72
\end{minipage}
\begin{minipage}[b]{3.0cm}
\centering
\epsfig{height=2.8cm, file=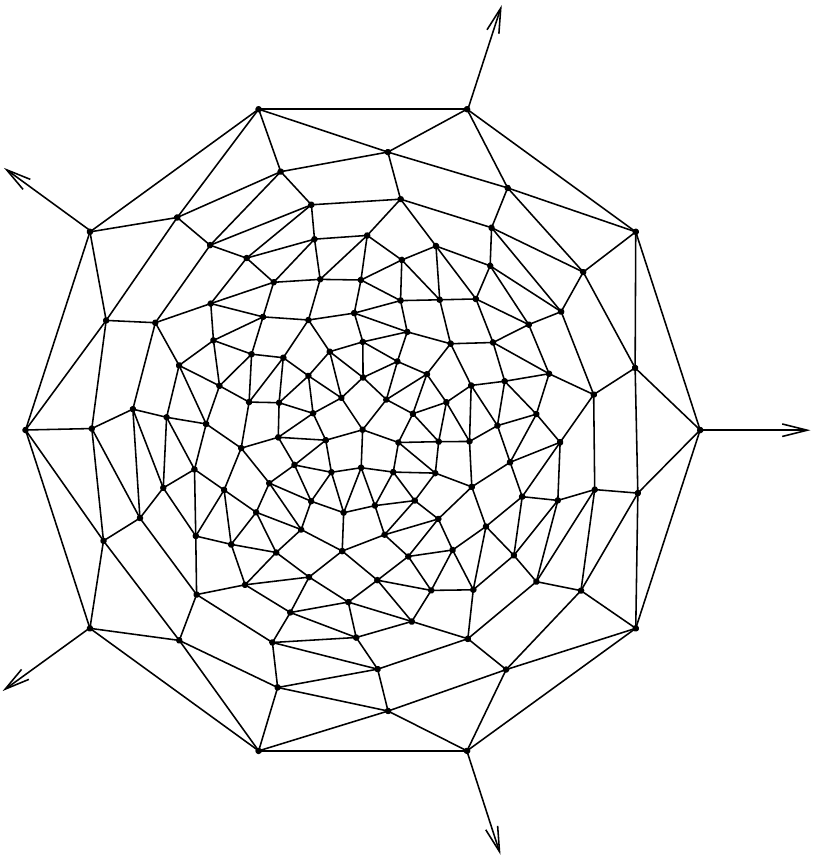}\par
$I$, 132, $4R_2$
\end{minipage}
\begin{minipage}[b]{3.0cm}
\centering
\epsfig{height=2.8cm, file=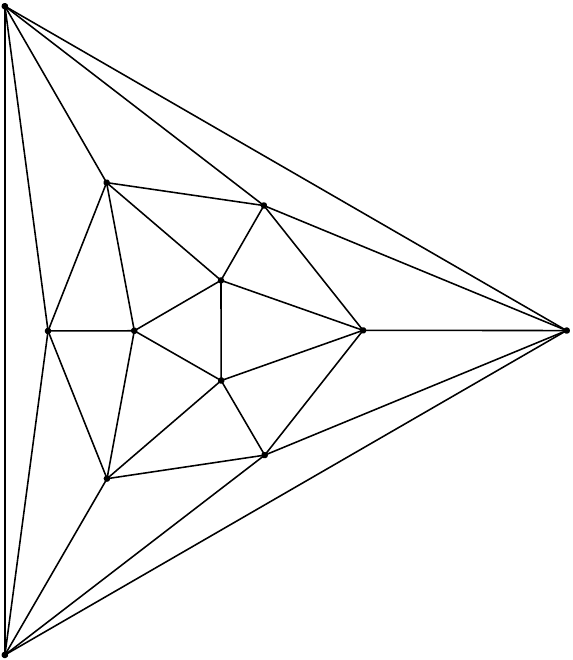}\par
$I_h$, 12, $3R_3$
\end{minipage}
\begin{minipage}[b]{3.0cm}
\centering
\epsfig{height=2.8cm, file=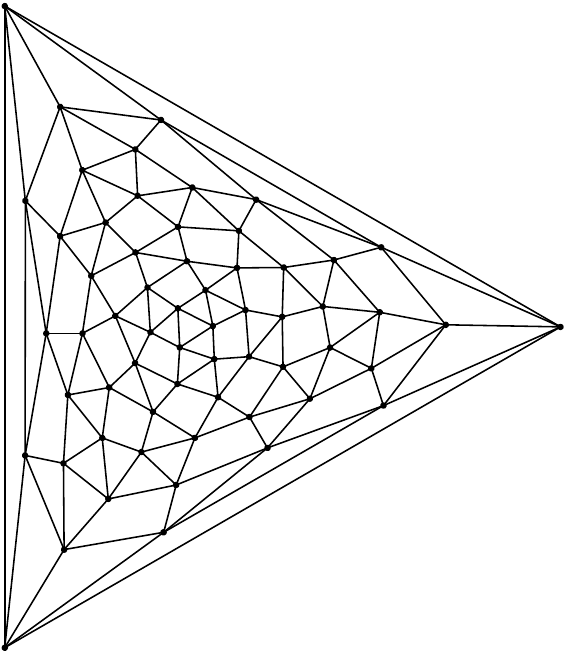}\par
$T$, 60, $4R_0$
\end{minipage}
\begin{minipage}[b]{3.0cm}
\centering
\epsfig{height=2.8cm, file=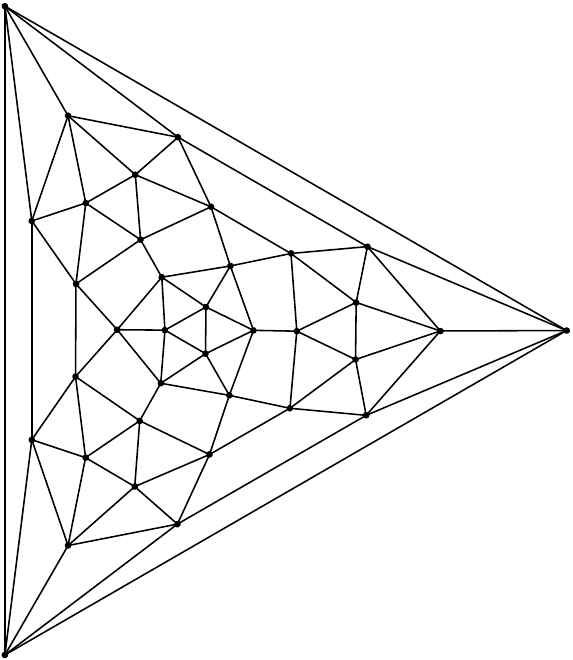}\par
$T_d$, 36, $4R_1$
\end{minipage}
\begin{minipage}[b]{3.0cm}
\centering
\epsfig{height=2.8cm, file=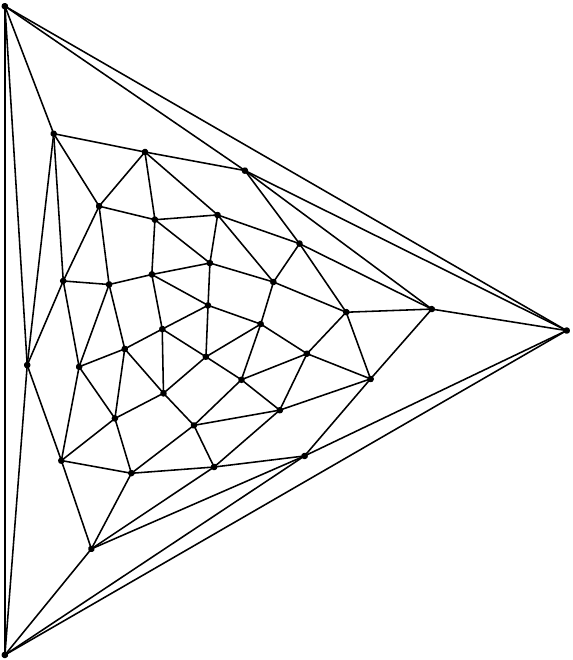}\par
$T_h$, 36, $4R_1$
\end{minipage}
\begin{minipage}[b]{3.0cm}
\centering
\epsfig{height=2.8cm, file=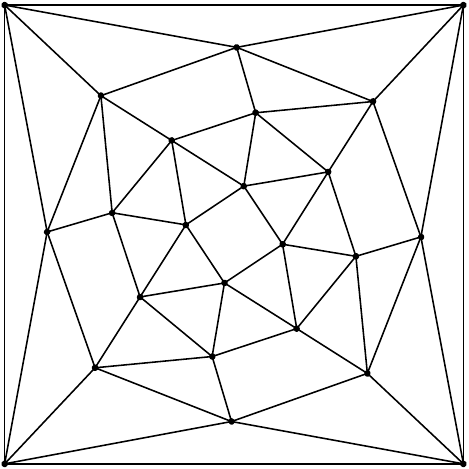}\par
$O$, 24, $4R_0$
\end{minipage}
\begin{minipage}[b]{3.0cm}
\centering
\epsfig{height=2.8cm, file=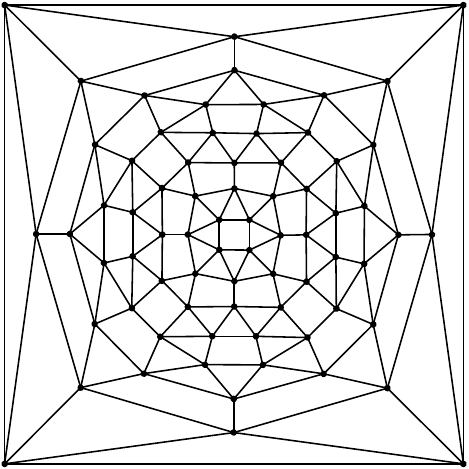}\par
$O_h$, 72
\end{minipage}
\end{center}
\caption{Examples of icosahedrites for all possible symmetry groups (part 2); all $5$ with at most $32$ vertices are minimal ones}
\label{ExampleIcosahedrite_Groups_2}
\end{figure}

The simplest icosahedrite is Icosahedron, which is a $(\{3\},5)$-sphere of
symmetry $I_h$.
One way to obtain icosahedrite is from an {\em octahedrite}, i.e. a
$(\{3,4\},4)$-sphere $G$.
To every vertex of $G$ we associate a square,
to every edge (coherently) a pair of adjacent triangles and faces
are preserved.
Only the rotational symmetries of $G$ are preserved in the final
icosahedrite $C(G)$.
If one applies $C$ to Octahedron, then
one gets the smallest ($24$ vertices) icosahedrite of symmetry $O$ (see Figure \ref{ExampleIcosahedrite_Groups_2}).
Applying it to the infinite regular plane tiling $\{4,4\}$ by squares,
one gets the Archimedean {\em snub square tiling} $(3.3.4.3.4)$.
Note that there is only one other {\em infinite Archimedean icosahedrite},
i.e., a vertex-transitive $5$-valent tiling of the plane by regular $3$-
and $4$-gons only: {\em elongated triangular tiling} $(3.3.3.4.4)$.

For a given icosahedrite, a {\em weak zigzag} $WZ$ is a 
circuit of edges such that one alternate between the left and right way 
but never extreme left or right.
The usual {\em zigzag} is a circuit such that one alternate between the 
extreme left and extreme right way.

A zigzag or weak zigzag is called {\em edge-simple}, respectively
{\em vertex-simple}, if any edge, respectively 
vertex, of it occurs only once.
A vertex-simple zigzag is also edge-simple.
If $WZ$ is vertex-simple
of length $l$, then one can construct another 
icosahedrite with $l$ more vertices.

Clearly, the weak zigzags (as well as the usual zigzags)
doubly cover the edge set. For example,
$30$ edges of  Icosahedron are doubly covered by $10$ weak vertex-simple
zigzags of length $6$, as well as by $6$ usual vertex-simple
zigzags of length $10$.
Clearly, a $3$-gon surrounded by $9$ $3$-gons or a  $4$-gon surrounded 
by $12$ $3$-gons corresponds to weak zigzags of length six or eight,
respectively.
In fact, no other vertex-simple weak zigzags exist.

\begin{figure}
\begin{center}
\begin{minipage}[b]{4.6cm}
\centering
\epsfig{height=3.9cm, file=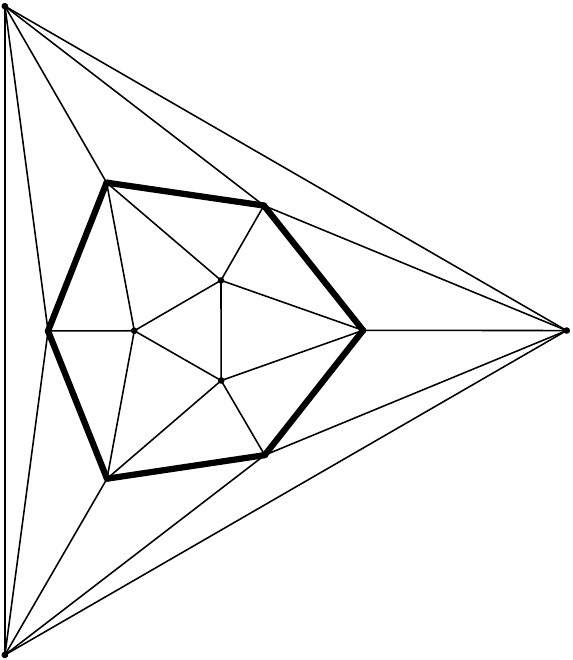}\par
12, $I_h$
\end{minipage}
\begin{minipage}[b]{4.6cm}
\centering
\epsfig{height=3.9cm, file=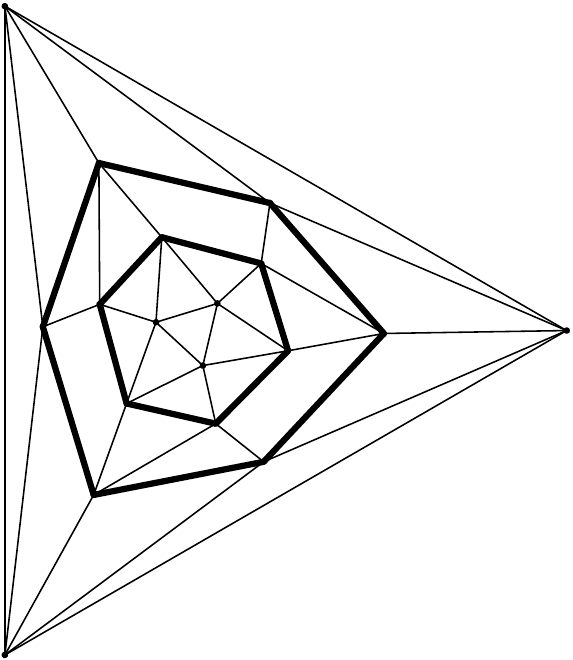}\par
18, $D_3$
\end{minipage}
\end{center}
\caption{A simple weak zigzag in an icosahedrite and the resulting expansion}
\label{WeakZigZagExpansion}
\end{figure}

In Table \ref{NrIcosahedrite} we list the number of $v$-vertex icosahedrites
for $v\le 32$. From this list it appears likely that any icosahedrite
is $3$-connected.

\begin{table}
\caption{The number of $v$-vertex icosahedrites}\label{NrIcosahedrite}
\begin{center}
\begin{tabular}{|c||c|c|c|c|c|c|c|c|c|c|c|}
\hline
v  &12 &14 &16 &18 &20 &22 & 24 & 26  & 28   & 30 & 32\\
\hline
Nr &1  &0  &1  &1  &5  &12 & 63 & 246 & 1395 &7668 & 45460\\
\hline
\end{tabular}
\end{center}
\end{table}

\begin{figure}
\begin{center}
\begin{minipage}[b]{4.6cm}
\centering
\epsfig{height=3.8cm, file=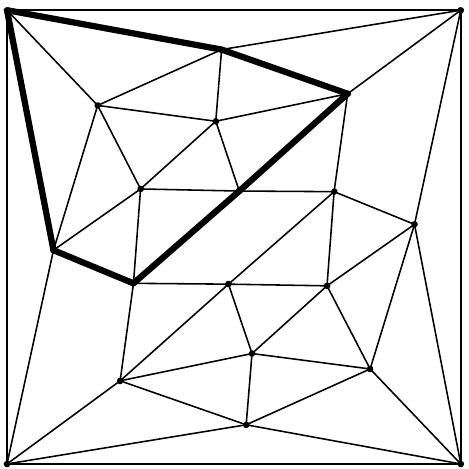}\par
20, $D_{2d}$
\end{minipage}
\begin{minipage}[b]{4.6cm}
\centering
\epsfig{height=3.8cm, file=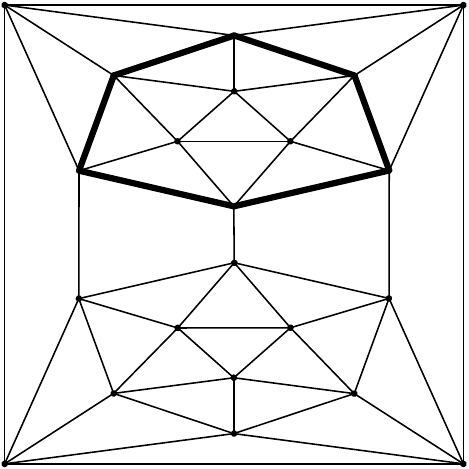}\par
22, $D_{5h}$
\end{minipage}
\end{center}
\caption{Two icosahedrites having a simple weak zigzag of length $6$}
\label{TwoSimpleIcosahedriteWZ6}
\end{figure}

\begin{theorem}
A $(\{3,4\},5)$-sphere exist if $v$ is even, $v\geq 12$ and $v\not= 14$. 
Their number grows at least exponentially with $v$.
\end{theorem}
\proof If there is an weak zigzag of length $6$ in an icosahedrite $G$,
then we can insert a {\em corona} ($6$-ring of three $4$-gons alternated
by three pairs of adjacent $3$-gons) instead of it
and get an icosahedrite with $6$ more vertices.
Since there are such icosahedrite for $n=18$, $20$, $22$
(see Figure \ref{WeakZigZagExpansion}, \ref{TwoSimpleIcosahedriteWZ6}), we can
generate the required graphs.
There is always two options when inserting the corona
and so the number grows exponentially as required. \qed

In the case of fullerenes or other standard $(\{a,b\},k)$-spheres,
the number of $a$-gons is fixed and the structure is made up of
patches of $b$-gons.
The parametrization of graphs, generalizing one in \cite{T}, 
including the case of one complex
parameter - {\em Goldberg-Coxeter construction} (\cite{goldberg}) -
are thus built.
A consequence of this is the polynomial growth of the number of 
such spheres.
This does not happen for the case of icosahedrites and their
parametrization, if any, looks elusive.

\begin{figure}
\begin{center}
\begin{minipage}[b]{6.1cm}
\centering
\epsfig{height=3cm, file=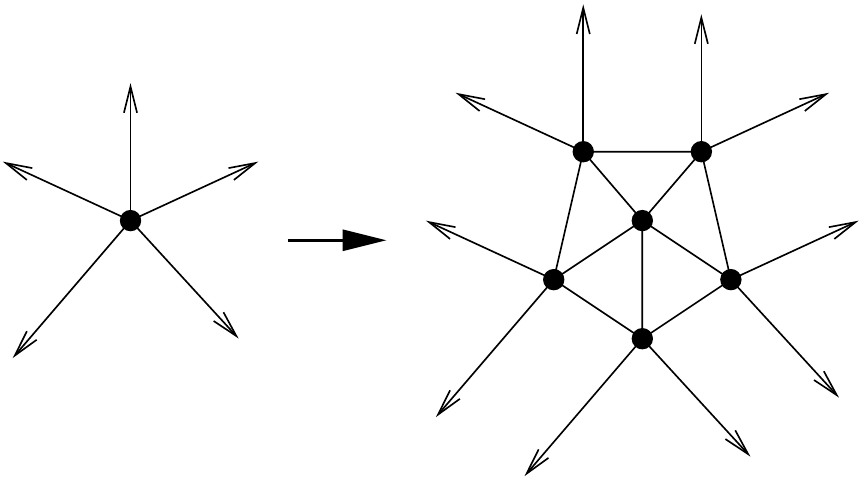}\par
Expansion operation $A$
\end{minipage}
\begin{minipage}[b]{6.1cm}
\centering
\epsfig{height=3cm, file=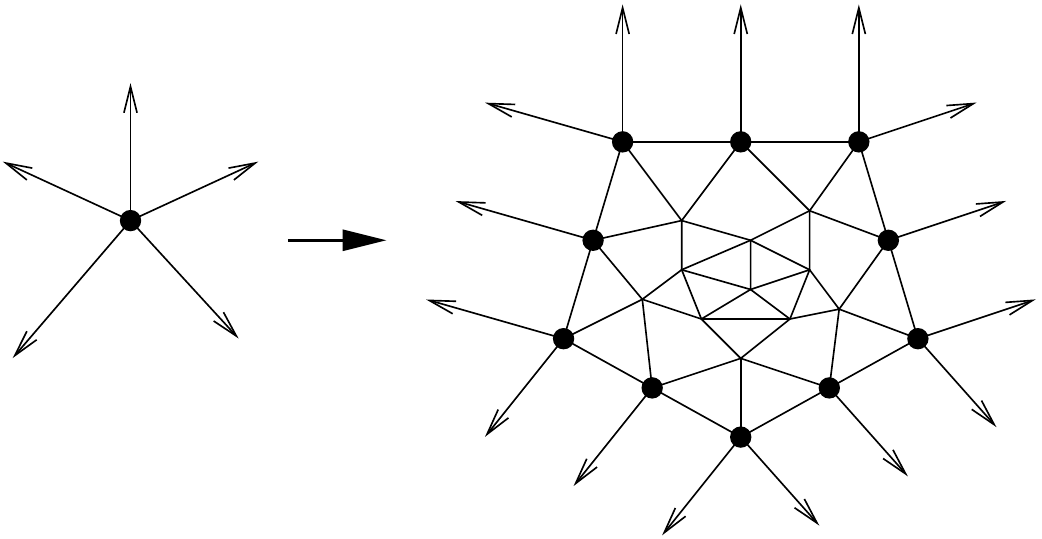}\par
Expansion operation $A'$
\end{minipage}
\begin{minipage}[b]{6.1cm}
\centering
\epsfig{height=3cm, file=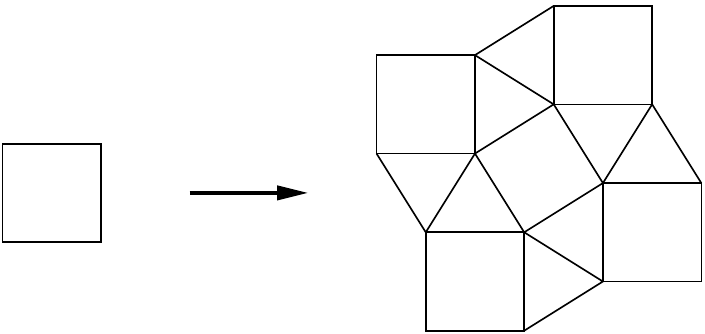}\par
Expansion operation $B_2$
\end{minipage}
\end{center}
\caption{Three expansion operations on icosahedrites}
\label{TwoLocalOperIcosahedrite}
\end{figure}

As a consequence of this increased freedom, we can more easily build
a new icosahedrite from a given one by an expansion operation, while in the
case of standard $(\{a,b\},k)$-spheres, one is essentially restricted to the 
Goldberg-Coxeter construction.
The operation $A$, respectively $A'$, replaces each vertex
of an icosahedrite $G$ by $6$, respectively $21$, vertices and gives
icosahedrites $A(G)$, $A'(G)$ with the same symmetries as $G$
(see Figure \ref{TwoLocalOperIcosahedrite}).
Moreover, for any $m\geq 2$ we can define an operation that replaces every
face $F$ by a patch and add $m-1$ rings of squares and pairs of triangles.
The operation $B_2$ on a $4$-gon is shown in
Figure \ref{TwoLocalOperIcosahedrite}.
The resulting map $B_m(G)$ has only the rotational symmetries of $G$
and associate to every vertex of $G$ $1+5m(m-1)$ vertices in $B_m(G)$.

\begin{theorem}\label{ListOfGroups}
Any symmetry group of icosahedrites if one of following $38$:

$C_1$, $C_i$, $C_s$, $S_4$, $S_6$, $S_8$, $S_{10}$.
$C_2$, $C_{2h}$, $C_{2v}$, 
$C_3$, $C_{3h}$, $C_{3v}$,
$C_4$, $C_{4h}$, $C_{4v}$,
$C_5$, $C_{5h}$, $C_{5v}$, 
$D_2$, $D_{2h}$, $D_{2d}$, 
$D_3$, $D_{3h}$, $D_{3d}$, 
$D_4$, $D_{4h}$, $D_{4d}$, 
$D_5$, $D_{5h}$, $D_{5d}$,
$O$, $O_h$, $T$, $T_d$, $T_h$, $I$, $I_h$.

\end{theorem}
\proof By the face sizes and vertex sizes, the list of possibilities
is the one indicated. We used the enumeration up to $32$ vertices
to find many groups and their minimal representative.
If a sphere has a $3$-fold axis, then it necessarily passes through two
$3$-gons. Those two $3$-gons can be replaced by $4$-gons and ipso facto
we get examples of a sphere with $4$-fold axis. From this we got the
icosahedrites with groups $C_{4}$, $C_{4h}$, $C_{4v}$, $D_{4}$, $D_{4h}$,
$D_{4d}$ and $S_8$.
All but $40$-vertex one $C_{4}$, $C_{4v}$, $D_{4h}$ are minimal.

\begin{figure}

\begin{center}
\begin{minipage}[b]{3.0cm}
\centering
\epsfig{height=2.8cm, file=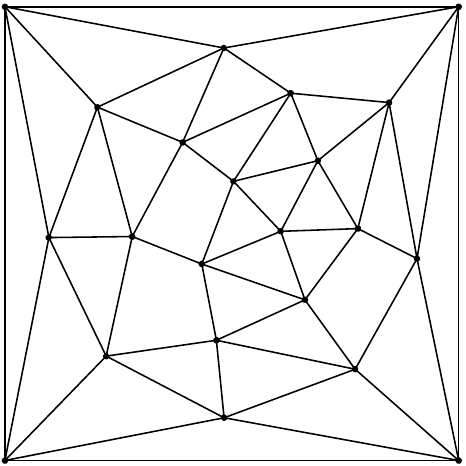}\par
$C_1$, 22
\end{minipage}
\begin{minipage}[b]{3.0cm}
\centering
\epsfig{height=2.8cm, file=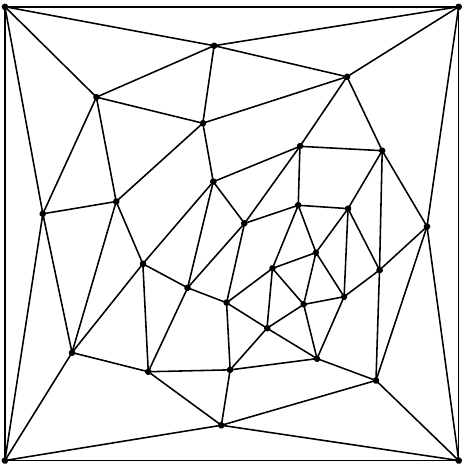}\par
$C_i$, 32
\end{minipage}
\begin{minipage}[b]{3.0cm}
\centering
\epsfig{height=2.8cm, file=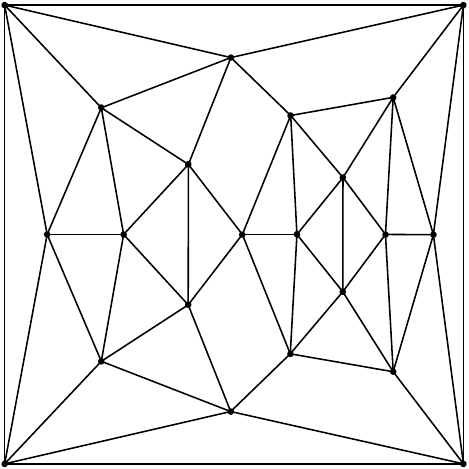}\par
$C_s$, 22
\end{minipage}
\begin{minipage}[b]{3.0cm}
\centering
\epsfig{height=2.8cm, file=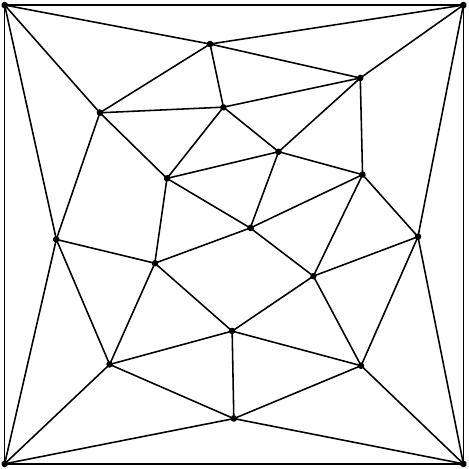}\par
$C_2$, 20, $4R_0$
\end{minipage}
\begin{minipage}[b]{3.0cm}
\centering
\epsfig{height=2.8cm, file=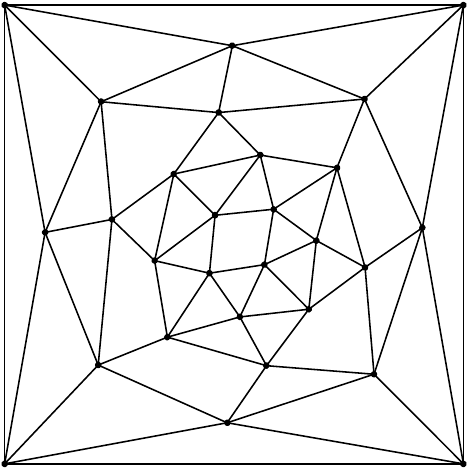}\par
$C_{2h}$, 28
\end{minipage}
\begin{minipage}[b]{3.0cm}
\centering
\epsfig{height=2.8cm, file=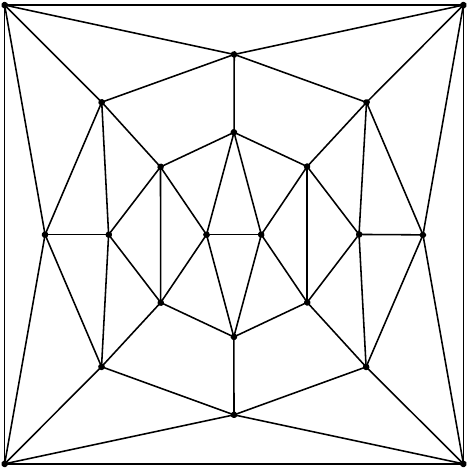}\par
$C_{2v}$, 22
\end{minipage}
\begin{minipage}[b]{3.0cm}
\centering
\epsfig{height=2.8cm, file=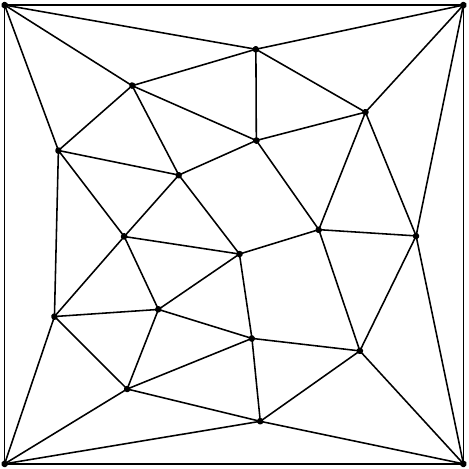}\par
$D_2$, 20, $4R_1$
\end{minipage}
\begin{minipage}[b]{3.0cm}
\centering
\epsfig{height=2.8cm, file=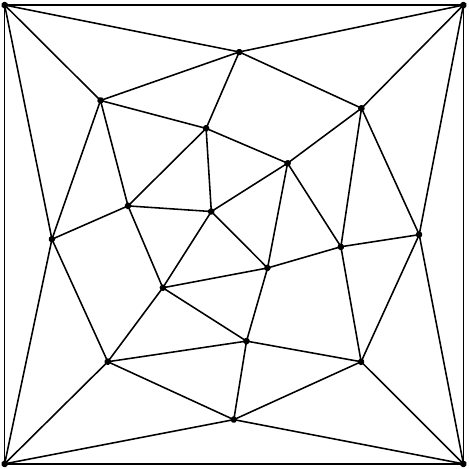}\par
$D_{2d}$, 20, $4R_0$
\end{minipage}
\begin{minipage}[b]{3.0cm}
\centering
\epsfig{height=2.8cm, file=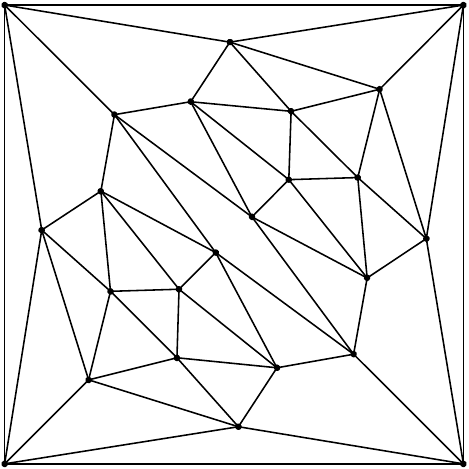}\par
$D_{2h}$, 24
\end{minipage}
\begin{minipage}[b]{3.0cm}
\centering
\epsfig{height=2.8cm, file=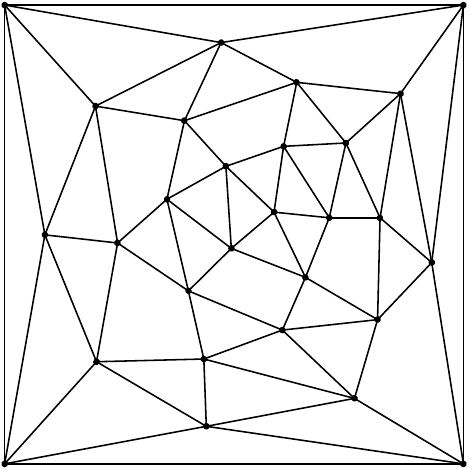}\par
$S_4$, 28, $4R_0$
\end{minipage}
\begin{minipage}[b]{3.0cm}
\centering
\epsfig{height=2.8cm, file=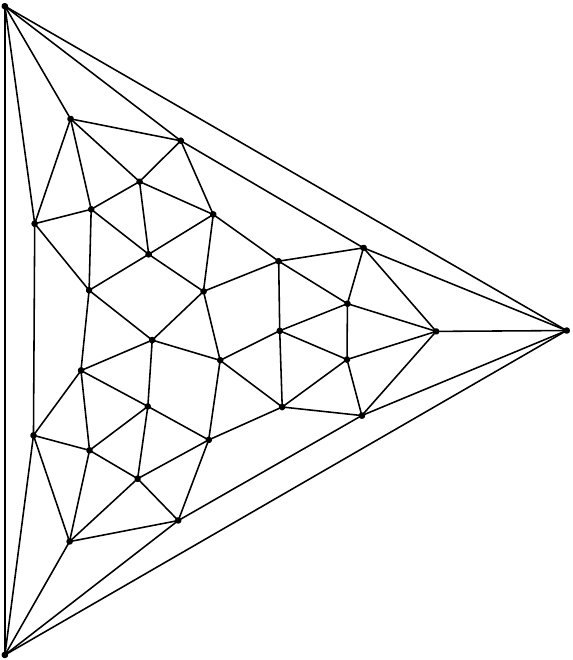}\par
$C_3$, 30
\end{minipage}
\begin{minipage}[b]{3.0cm}
\centering
\epsfig{height=2.8cm, file=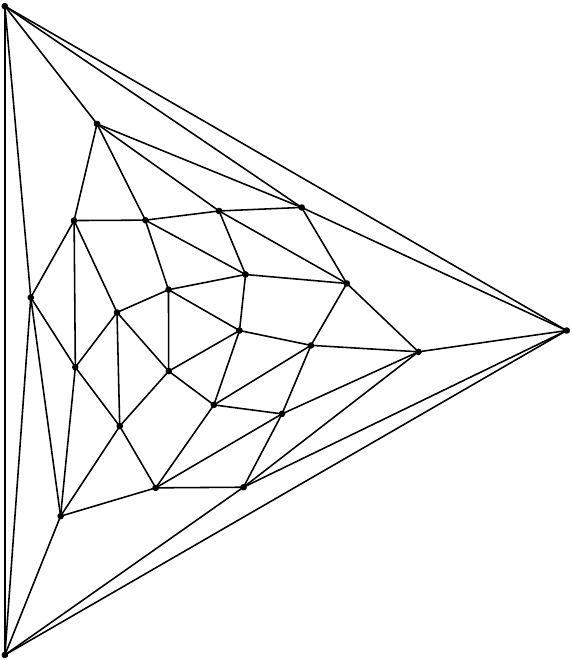}\par
$C_{3h}$, 24, $4R_0$
\end{minipage}
\begin{minipage}[b]{3.0cm}
\centering
\epsfig{height=2.8cm, file=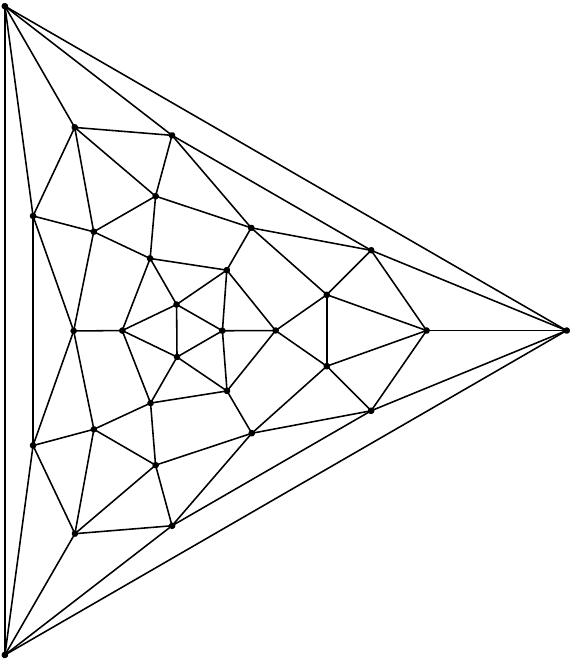}\par
$C_{3v}$, 30
\end{minipage}
\begin{minipage}[b]{3.0cm}
\centering
\epsfig{height=2.8cm, file=PLminD3_18sec.pdf}\par
$D_3$, 18, $4R_0$
\end{minipage}
\begin{minipage}[b]{3.0cm}
\centering
\epsfig{height=2.8cm, file=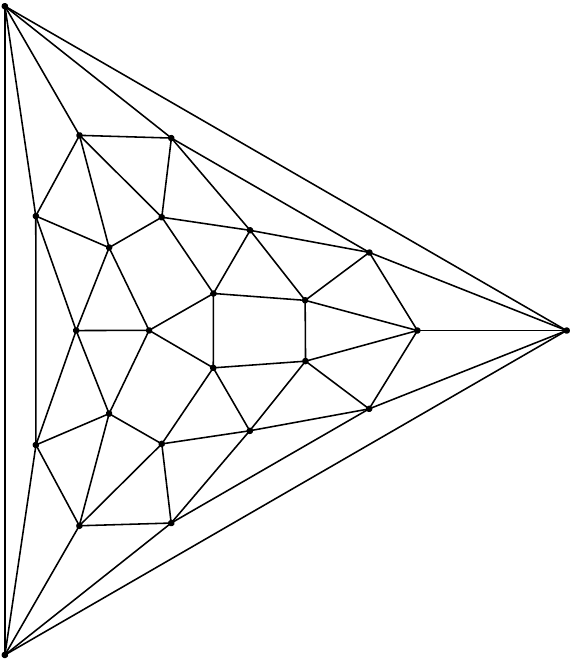}\par
$D_{3d}$ 24, $4R_0$
\end{minipage}
\begin{minipage}[b]{3.0cm}
\centering
\epsfig{height=2.8cm, file=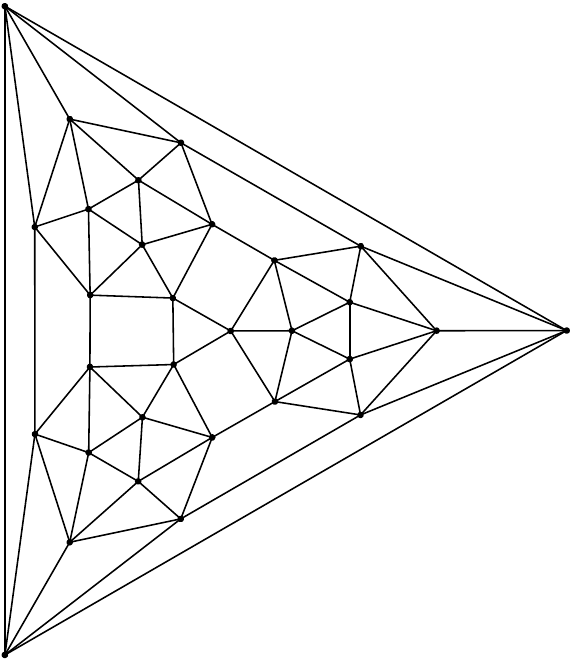}\par
$D_{3h}$, 30
\end{minipage}
\begin{minipage}[b]{3.0cm}
\centering
\epsfig{height=2.8cm, file=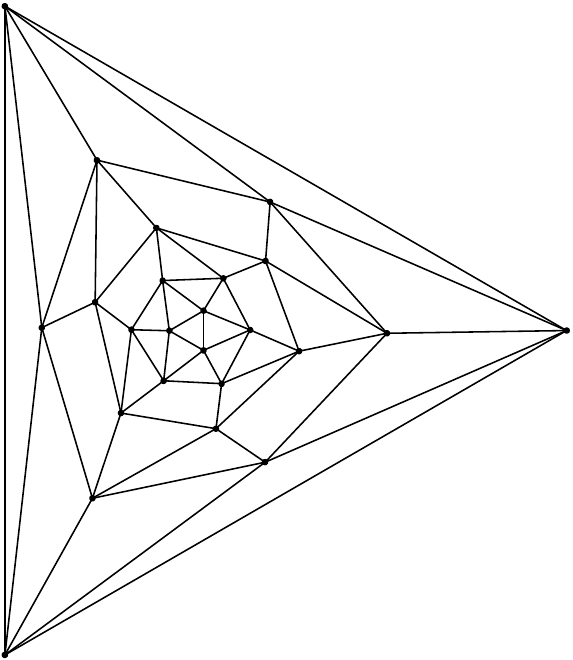}\par
$S_6$, 24, $4R_0$
\end{minipage}
\begin{minipage}[b]{3.0cm}
\centering
\epsfig{height=2.8cm, file=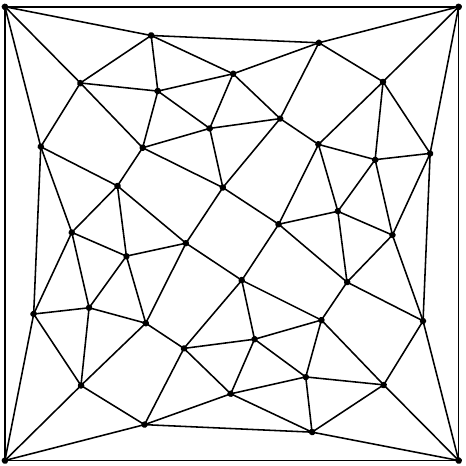}\par
$C_4$, 40
\end{minipage}
\begin{minipage}[b]{3.0cm}
\centering
\epsfig{height=2.8cm, file=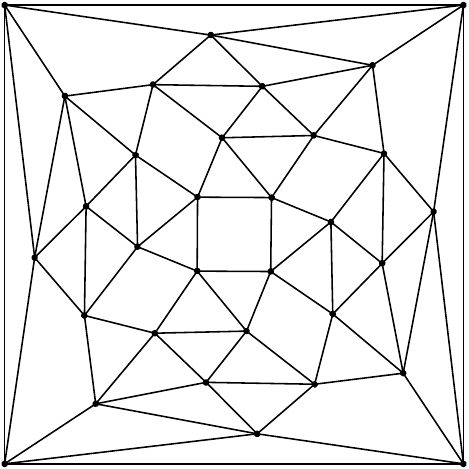}\par
$C_{4h}$, 32, $4R_0$
\end{minipage}
\begin{minipage}[b]{3.0cm}
\centering
\epsfig{height=2.8cm, file=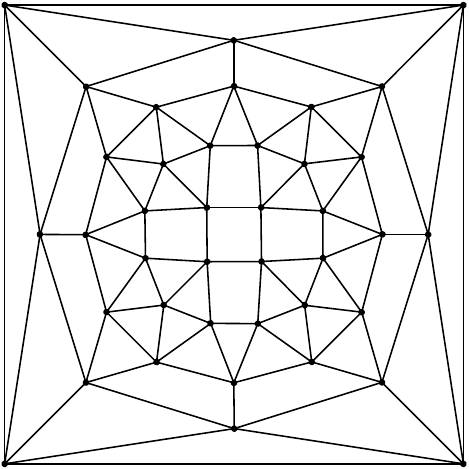}\par
$C_{4v}$, 40
\end{minipage}
\begin{minipage}[b]{3.0cm}
\centering
\epsfig{height=2.8cm, file=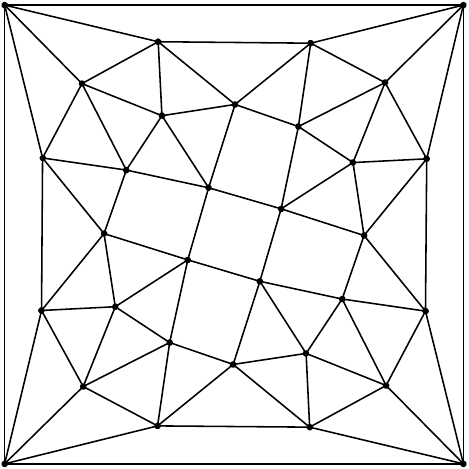}\par
$D_{4}$, 32
\end{minipage}
\begin{minipage}[b]{3.0cm}
\centering
\epsfig{height=2.8cm, file=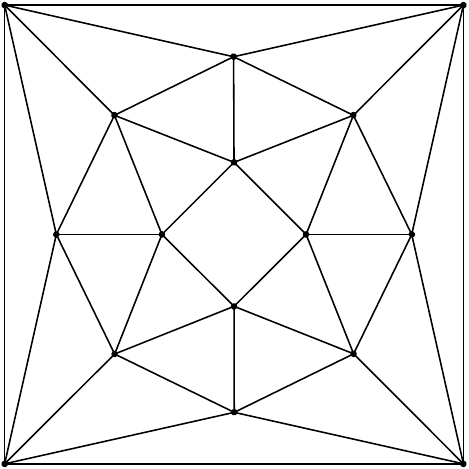}\par
$D_{4d}$, 16, $4R_0$
\end{minipage}
\begin{minipage}[b]{3.0cm}
\centering
\epsfig{height=2.8cm, file=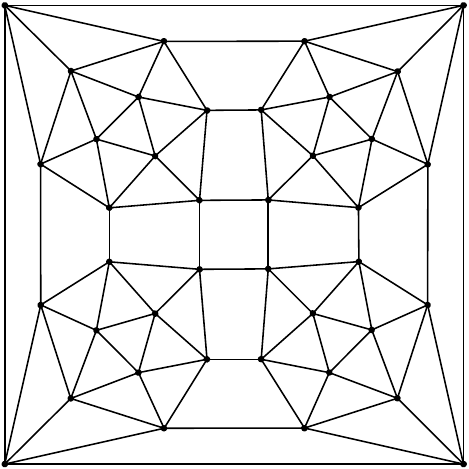}\par
$D_{4h}$, 40
\end{minipage}
\begin{minipage}[b]{3.0cm}
\centering
\epsfig{height=2.8cm, file=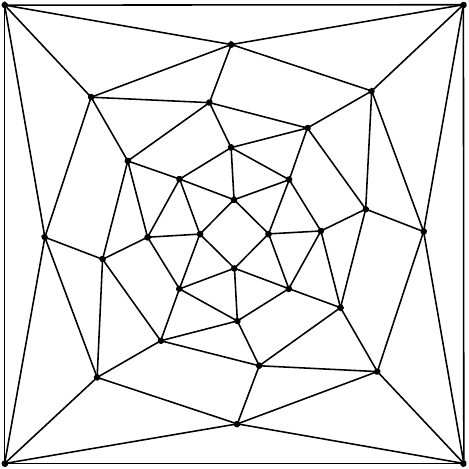}\par
$S_{8}$, 32, $4R_0$
\end{minipage}
\end{center}
\caption{Examples of icosahedrites for all possible symmetry groups (part 1); all $21$ with at most $32$ vertices are minimal ones}
\label{ExampleIcosahedrite_Groups_1}
\end{figure}

Now, suppose that a $2$-fold axis of rotation passes by two edges, 
which are both contained in two triangles.
Then we insert a vertex on those two edges
and replace the $2$-fold symmetry
by a $5$-fold symmetry. By iterating over all known isocahedrites, and
all such $2$-fold axis, we get the symmetries $C_{5}$, $C_{5h}$, $C_{5v}$,
$S_{10}$ and  $D_{5h}$, $D_{5d}$, $D_5$; last three are minimal.
For the cases of $O_h$, $T$, $T_d$ and $T_h$, we obtained examples
by hand drawing.
The $132$ vertices icosahedrite of symmetry $I$ in Figure \ref{ExampleIcosahedrite_Groups_2} is obtained from Isocahedron by operation $B_2$.
The minimal known (actually, examples are minimal whenever $v\le 32$) are
given in Figure \ref{ExampleIcosahedrite_Groups_1}, \ref{ExampleIcosahedrite_Groups_2}.

Aggregating groups ${\bf C_1}$=$\{C_1,C_s,C_i\}$,
${\bf C_m}$=$\{C_m$, $C_{mv}$, $C_{mh}$, $S_{2m}\}$,
${\bf D_m}$=$\{D_m$, $D_{mh}$, $D_{md}\}$, ${\bf T}$=$\{T$, $T_d$, $T_h\}$, 
${\bf O}$=$\{O,O_h\}$, ${\bf I}$=$\{I,I_h\}$,
all $38$ symmetries of $(\{3,4\},5)$-spheres are:
${\bf C_1}$,  ${\bf C_m}$,  ${\bf D_m}$ for $2\le m\le 5$ and  
${\bf T}$, ${\bf O}$, ${\bf I}$.
$5$-, $4$- and $3$-fold symmetry  exists if and only if, respectively,  
$p_4\equiv 0 \pmod 5$ (i.e.,  $v=2p_4+12\equiv 2 \pmod {10}$), 
$p_4\equiv 2 \pmod 4$ (i.e.,  $v=2p_4+12\equiv 0 \pmod 8$) and 
$p_4\equiv 0 \pmod 3$ (i.e.,  $v=2p_4+12\equiv 0 \pmod 6$).

Any group appear an infinite number of times since, for
example, one gets an infinity by applying construction $A$ iteratively.

From the above result it appears that the only limitations
for the group are coming from the rotation axis.
It seems possible that this also occurs for all $(\{a,b\},k)$-spheres
with $b$-gons being of negative curvature.

A map is said to be {\em face-regular}, or, specifically,
$pR_i$ if every face of size $p$ is adjacent to exactly $i$ faces of
the same size $p$.
\begin{theorem}\label{FaceRegularIcosahedrite}
(i) The only icosahedrite which is $3R_i$ is Icosahedron which is $3R_3$

(ii) For $i=0$, $1$ and $2$ there is an infinity of icosahedrites
that are $4R_i$.
\end{theorem}
\proof Let $N_{ij}$ denote the number of $i$-$j$ edges, 
i.e., those which are common to an $i$-gon
and an $j$-gon; so, $e=N_{33}+N_{34}+N_{44}$.
But $N_{34}+N_{44} \le 4p_4$ with equality if and only if our icosahedrite
is $4R_0$.
So, 
\begin{equation*}
2e=3p_3+4p_4=p_3+2(20+2p_4)+4p_4 \le 2N_{33}+8p_4,
\end{equation*}
i.e., $2N_{33}\ge p_3+40$.
It excludes the cases  $3R_0$ and $3R_1$, since then $2N_{33}=0$ and $p_3$,
respectively.

If our icosahedrite is $3R_2$, then $N_{33}+N_{34}= p_3+p_3$, implying 
$N_{44}=p_4-10$. Any $4R_0$ icosahedrite with $32$ vertices (i.e., 
with $p_4=10$) has
$\frac{N_{33}}{p_3}=2$, i.e., it is  {\em $3R_2$ in average}.
Now, $3R_2$ means that the $3$-gons are organized in rings  
separated by $4$-gons.
Such rings can be either five $3$-gons with common vertex, or $12$ 
$3$-gons around a $4$-gon, or $APrism_m$, $m>4$.
In each case, such ring should be completed to an icosahedrite by
extruding edge to keep $5$-regularity.
In order to be isolated from other $3$-gons, the ring should have
twice longer ring  of $4$-gons around it.
The faces touching the ring of $3$-gons in a vertex only, could not 
be $3$-gonal since they have two $4$-gonal neighbors.
So, the isolating ring of $4$-gons consists of $4$-gons adjacent (to the 
ring of $3$-gons) alternated by $4$-gons touching it only in a vertex.
Easy to see that such  process cannot be closed.

\begin{figure}
\begin{center}
\begin{minipage}[b]{4.1cm}
\centering
\epsfig{height=2.8cm, file=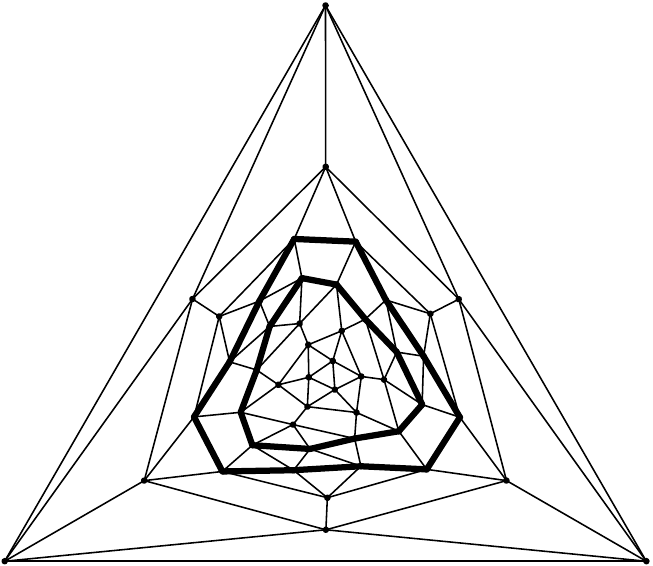}\par
48, $D_3$
\end{minipage}
      
\end{center}
\caption{The icosahedrite used in proof of Theorem \ref{FaceRegularIcosahedrite} as the first one in an infinite series with $4R_1$}
\label{Exmp4R1}
\end{figure}

The number of icosahedrites $4R_0$ is infinite; such series can be 
obtained by the operation $A$ (see Figure \ref{TwoLocalOperIcosahedrite})
from, say,  unique $16$-vertex icosahedrite.
An infinity of $4R_1$ icosahedrites can be obtained in the following way.
We take the $4R_1$ icosahedrite of Figure \ref{Exmp4R1} and the ring
bounded by two overlined circuits can be transformed into any number
of concentrated rings.
An infinity of icosahedrites $4R_2$ can be obtained from the $22$-vertices
icosahedrite of symmetry $D_{5h}$ in Figure \ref{ExampleIcosahedrite_Groups_2}.
It suffices to add layers of five $4$-gons alternated by layers of ten $3$-gons
(as in Archimedean elongated triangular plane tiling $(3.3.3.4.4)$. \qed

Among all $54,851$ icosahedrites with $v\le 32$, there are no $4R_3$ 
and only four specimens $4R_2$. Those with $4R_3$ should have 
$2N_{33}=3p_3-p_4$ and $v$ divisible by $4$, but we doubt they exist.
Note  that part (i) of Theorem \ref{FaceRegularIcosahedrite} can be
easily generalized on $(\{3,b\},5)$-sphere with any $b\ge 4$.

\begin{figure}
\begin{center}
\begin{minipage}[b]{3.0cm}
\centering
\epsfig{height=2.8cm, file=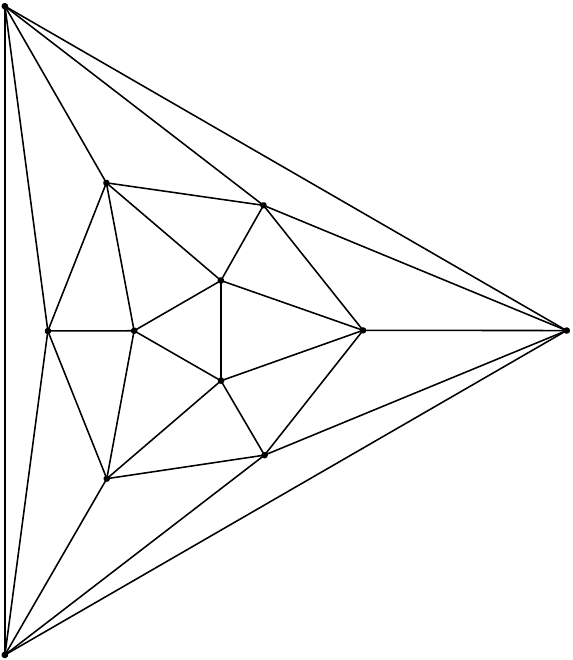}\par
12, $I_h$\par
${6^{10}}$, $3R_3$
\end{minipage}
\begin{minipage}[b]{3.0cm}
\centering
\epsfig{height=2.8cm, file=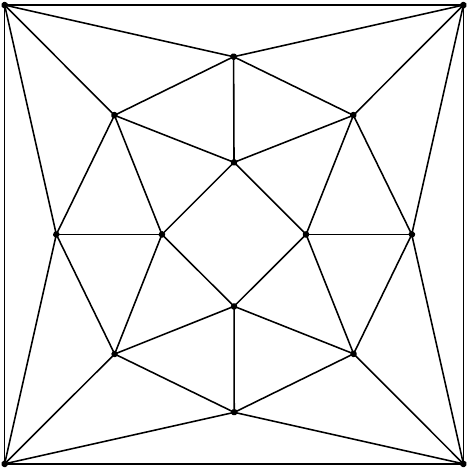}\par
16, $D_{4d}$\par
$(8,24^3)$, $4R_0$
\end{minipage}
\begin{minipage}[b]{3.0cm}
\centering
\epsfig{height=2.8cm, file=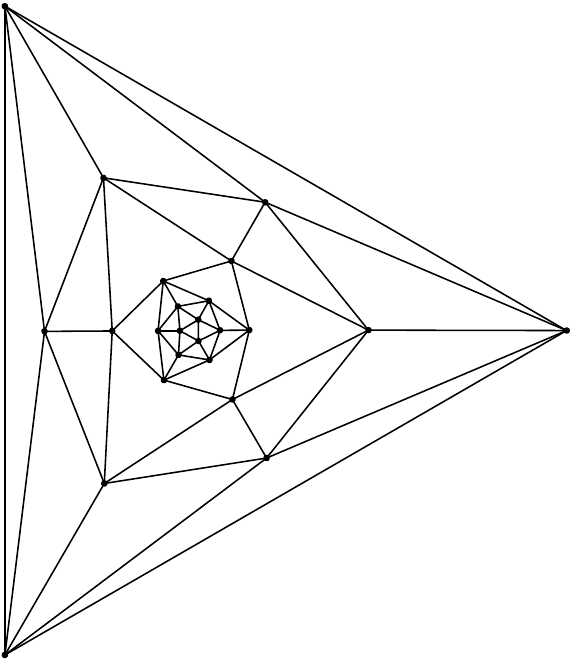}\par
24, $D_{3d}$\par
$(6^8,36^2)$, $4R_2$
\end{minipage}
\begin{minipage}[b]{3.0cm}
\centering
\epsfig{height=2.8cm, file=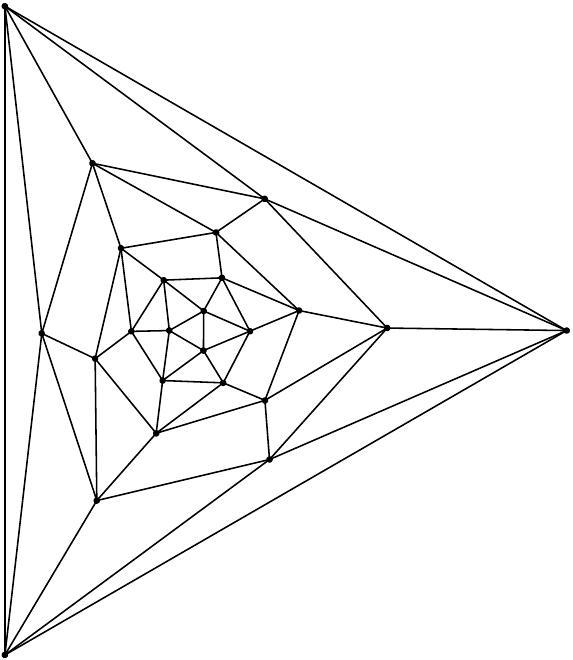}\par
24, $S_{6}$\par
$(6^3,10^2,36^2)$, $R_0$
\end{minipage}
\begin{minipage}[b]{3.4cm}
\centering
\epsfig{height=2.8cm, file=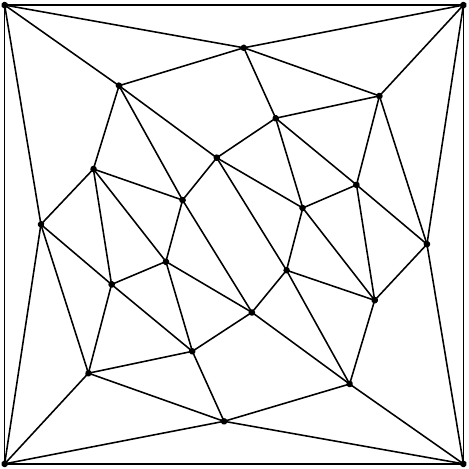}\par
24, $C_{2}$\par
$(10^2,14^2,36^2)$, $4R_0$
\end{minipage}
\begin{minipage}[b]{3.0cm}
\centering
\epsfig{height=2.8cm, file=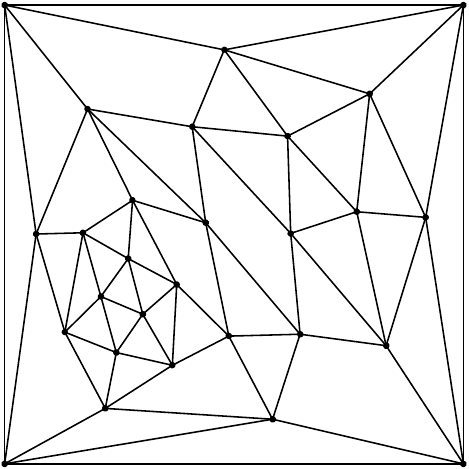}\par
28, $D_{2h}$\par
$(6^6,12^4,14^4)$, $4R_0$
\end{minipage}
\begin{minipage}[b]{3.0cm}
\centering
\epsfig{height=2.8cm, file=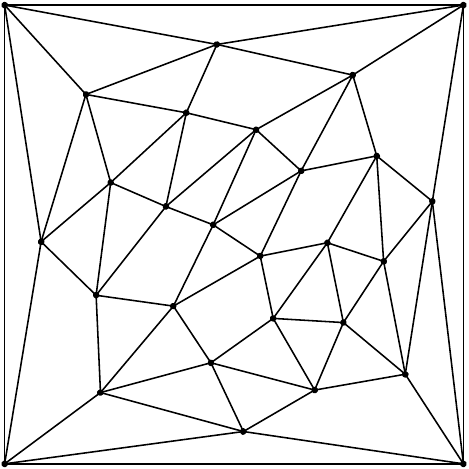}\par
28, $C_{1}$\par
$(8,42^2,48)$, $4R_0$
\end{minipage}

\end{center}
\caption{The icosahedrites with only edge-simple weak zigzags and at most $32$ vertices}
\label{IcosahedriteSimple}
\end{figure}
All icosahedrites with only edge-simple usual
zigzags and $v\le 32$  are three of those seven (see Figure
\ref{IcosahedriteSimple}) having only edge-simple weak
zigzags: $12$, $I_h$ with  $(10^6)$,
$24$, $D_{3d}$ with $(6,10^6,18^3)$ and $28$, $D_{2h}$ with $(10^6,20^4)$.
We expect that  any icosahedrite with only edge-simple usual or weak 
zigzags, if such ones exist for $v\ge 34$, has $v$ divisible by four.
Also, for $v\le 32$, the maximal number of zigzags and weak zigzags are 
realized only by edge-simple ones whenever they exist.

Snub $APrism_b$, Snub Cube and Snub Dodecahedron with $(v,b)=(4b,b)$, $(24,4)$  and $(60,5)$, respectively, are $bR_0$ and only known
$b$-gon-transitive $(\{3,b\},5)$-spheres.
They are also, besides Icosahedron, only known 
$(\{3,b\},5)$-spheres with at most two orbits of $3$-gons.
The Archimedean $(\{3,b\},5)$-plane tilings $(3.3.4.3.4)$, $(3.3.3.4.4)$, 
$(3.3.3.3.6)$  have $b=4$, $4$, $6$, respectively.
They are transitive on $b$-gons and vertices.
They are also, respectively, $(4R_1,3R_1)$, $(4R_2,3R_2)$, $6R_0$ and 
have $1$, $1$, $2$ orbits of $3$-gons.

\section{On $(\{a,b,c\},k)$-spheres with $p_c=1$}

Clearly, an $(\{a,b,c\},k)$-sphere with $p_c=1$ has
\begin{equation*}
v=\frac{2}{k-2}(p_a-1+p_b)=\frac{2}{2k-a(k-2)}(a+c+p_b(b-a))
\end{equation*}
vertices and (setting $b'$=$\frac{2k}{k-2}$)
$p_a=\frac{b'+c}{b'-a}+p_b\frac{b-b'}{b'-a}$ $a$-gons.
So, $p_a=\frac{b+c}{b-a}$ if $b=b'$, i.e., $(\{a,b\},k)$-sphere is standard.

We are especially interested in 
{\em fullerene $c$-disks}, i.e., $(\{5,6,c\},3)$-spheres with $p_c=1$.
It exists for any $c\geq 1$ and has $p_5=c+6$, $v=2(p_6 + c + 5)$.
Clearly, there is an infinity of fullerene $c$-disks for any $c\ge 1$.

The only way to to get fullerene $1$-disk, is to get a
$(\{5,6,r,s\}),3)$-sphere with $p_r=p_s=1$,
$(r,s)=(3,4)$, $(3,3)$, $(2,3)$, $(2,4)$  and add $4$ vertices: $2$-gon on the
$r-s$ edge and then erect an edge with $1$-gon from the middle of $2-r$
edge. In fact, only $(3,4)$ is possible and minimal such graph has $36$
vertices, proving minimality of $40$-vertex $1$-disk.

One can check that all $(\{a,b,c\},3)$-spheres with $p_c=1$, $2\le a,b,c\le 6$
and $c\neq a,b$  are 
four series with $(a,b,c)=(4,6,2)$, $(5,6,2)$, $(5,6,3)$, $(5,6,4)$
and following two $10$-vertex spheres of
symmetry $C_{3v}$ having $p_a=p_b=3$, $p_c=1$: $(\{4,5,3\},3)$-
(Cube truncated on one vertex) 
and  $(\{3,5,6\},3)$- (Tetrahedron truncated on three vertices).
There are also series of $(\{1,4,2\},4)$- and $(\{3,4,2\},4)$-spheres
with $p_2=1$.

Theorem \ref{Case_PB1} implies that $(\{a,b,c\},k)$-sphere with $p_c=1$
and $p_b=0$ has $c=a$, i.e., it is the $k$-regular map $\{a,k\}$
on the sphere.
We conjecture that a $(\{a,b,c\},k)$-sphere with $p_c=p_b=1<a,c$ has $c=b$,
i.e., it is a  $(\{a,b\},k)$-sphere with $p_b=2$.
Note that the $2$-vertex $(\{1,4,2\},4)$-sphere has $p_2=p_4=1=a < c < b$.

Call {\em $(\{a,b,c\},k)$-thimble} any $(\{a,b,c\},k)$-sphere with $p_c=1$
such that the $c$-gon is adjacent only to $a$-gons.
{\em Fullerene $c$-thimble} is the case $(a,b,c;k)=(5,6,c;3)$ of above. 
It exists if and only if $c\ge 5$.

\begin{figure}
\begin{center}
\begin{minipage}[b]{3.0cm}
\centering
\epsfig{height=2.8cm, file=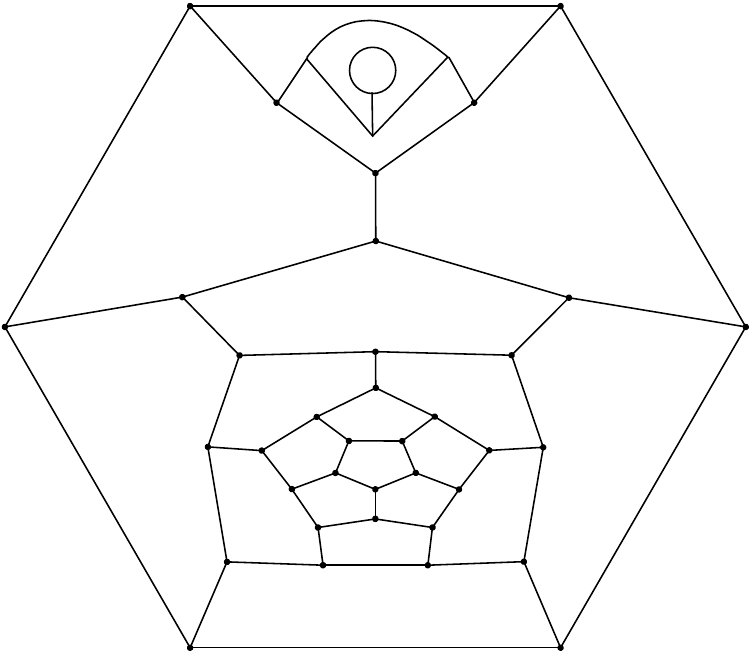}\par
$1$: 40, $C_{s}$
\end{minipage}
\begin{minipage}[b]{3.0cm}
\centering
\epsfig{height=2.4cm, file=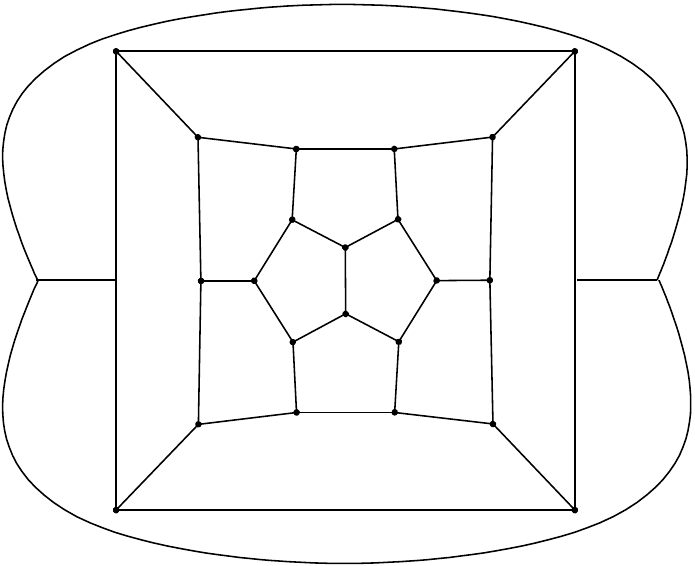}\par
$2$: 26, $C_{2v}$
\end{minipage}
\begin{minipage}[b]{3.0cm}
\centering
\epsfig{height=2.8cm, file=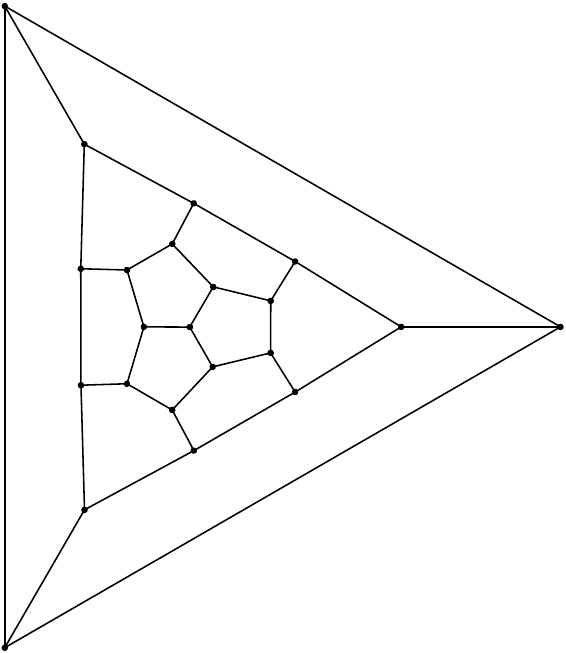}\par
$3$: 22, $C_{3v}$
\end{minipage}
\begin{minipage}[b]{3.0cm}
\centering
\epsfig{height=2.8cm, file=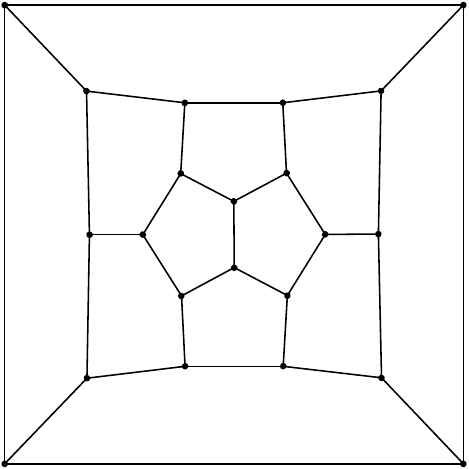}\par
$4$: 22, $C_{2v}$
\end{minipage}
\begin{minipage}[b]{3.0cm}
\centering
\epsfig{height=2.8cm, file=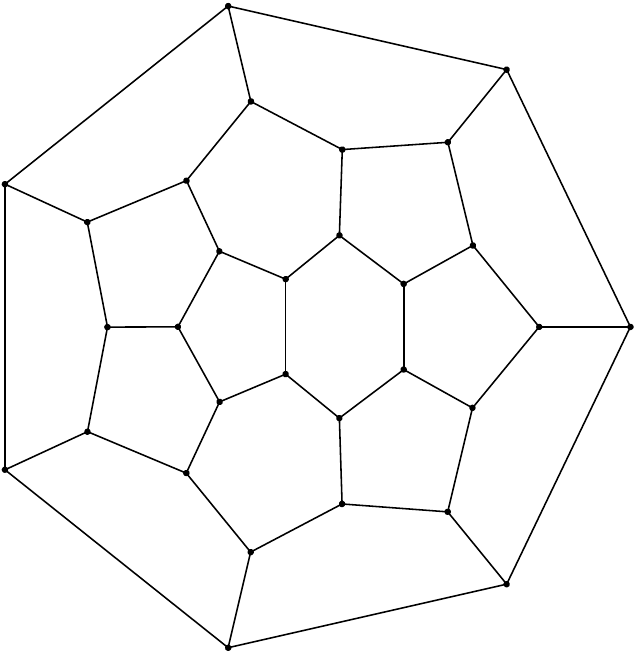}\par
$7$: 30, $C_s$
\end{minipage}
\begin{minipage}[b]{3.0cm}
\centering
\epsfig{height=2.8cm, file=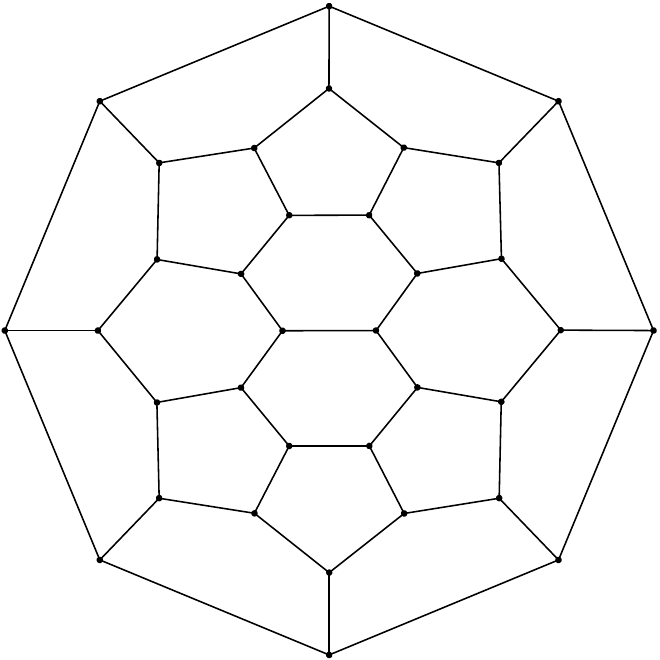}\par
$8$: 34, $C_{2v}$
\end{minipage}
\begin{minipage}[b]{3.0cm}
\centering
\epsfig{height=2.8cm, file=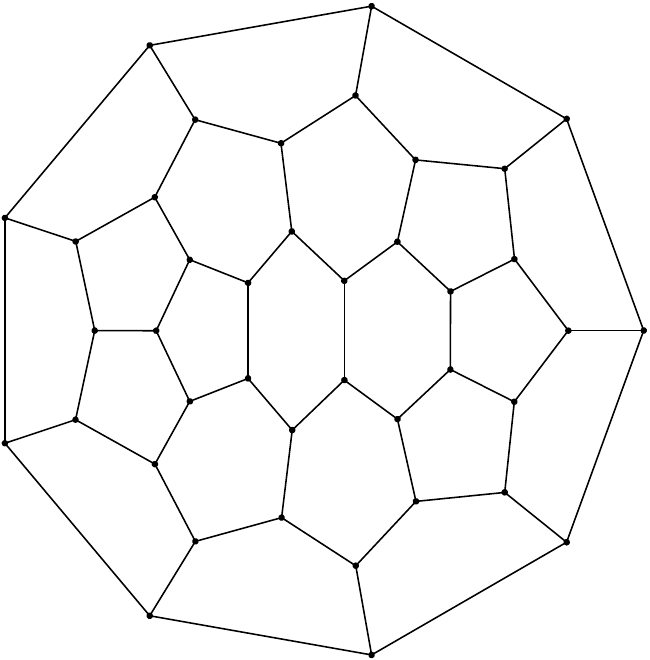}\par
$9-1$: 40, $C_s$
\end{minipage}
\begin{minipage}[b]{3.0cm}
\centering
\epsfig{height=2.8cm, file=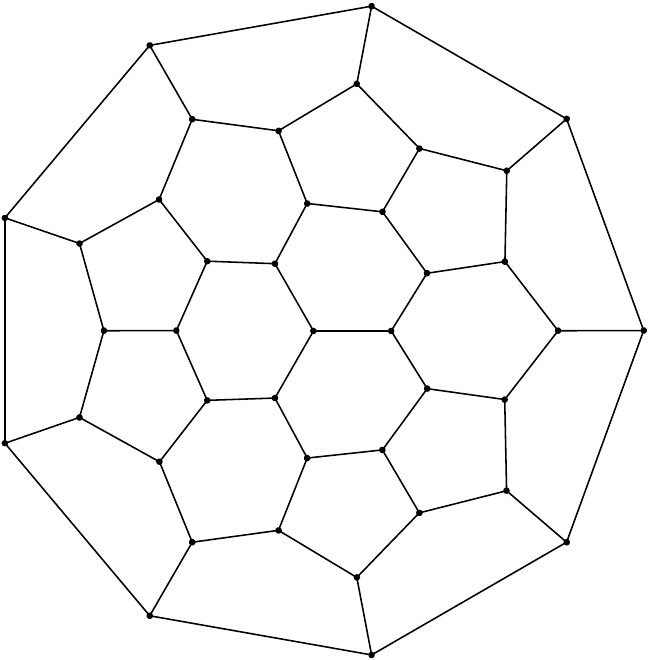}\par
$9-2$: 40, $C_{3v}$
\end{minipage}
\begin{minipage}[b]{3.0cm}
\centering
\epsfig{height=2.8cm, file=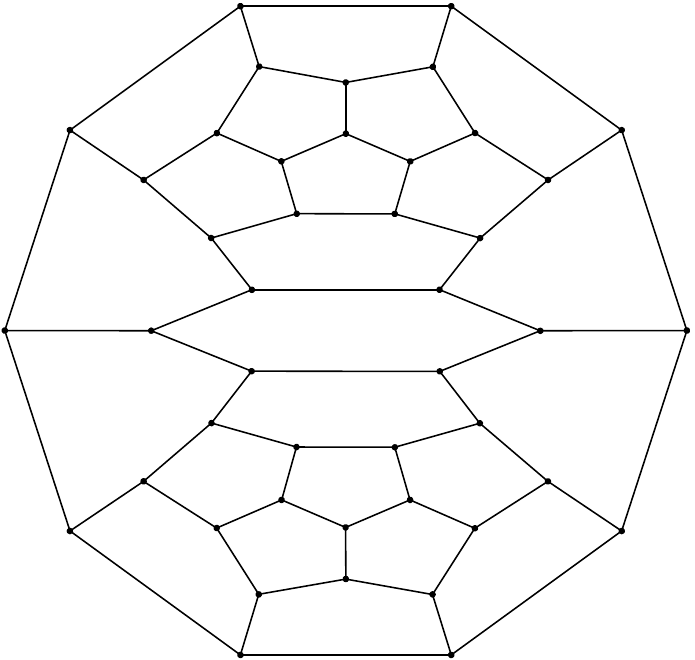}\par
$10-1$: 44, $C_{2v}$
\end{minipage}
\begin{minipage}[b]{3.0cm}
\centering
\epsfig{height=2.8cm, file=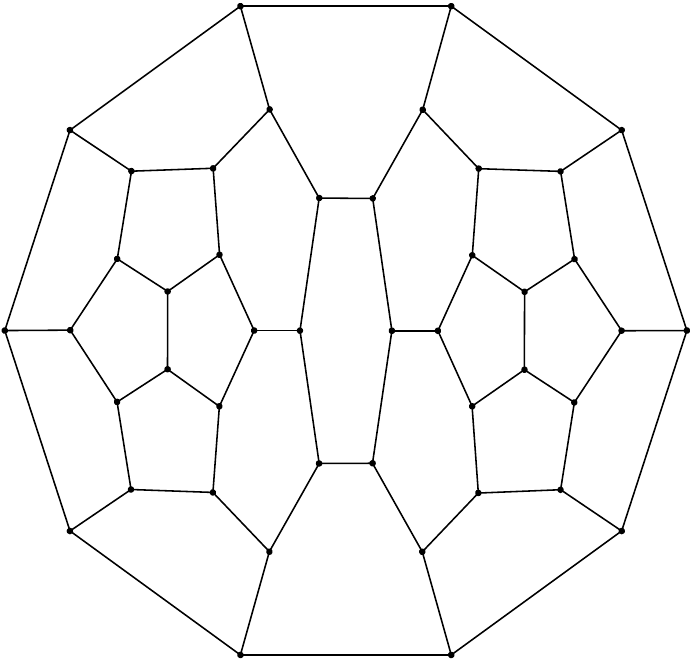}\par
$10-2$: 44, $C_{2v}$
\end{minipage}
\begin{minipage}[b]{3.0cm}
\centering
\epsfig{height=2.8cm, file=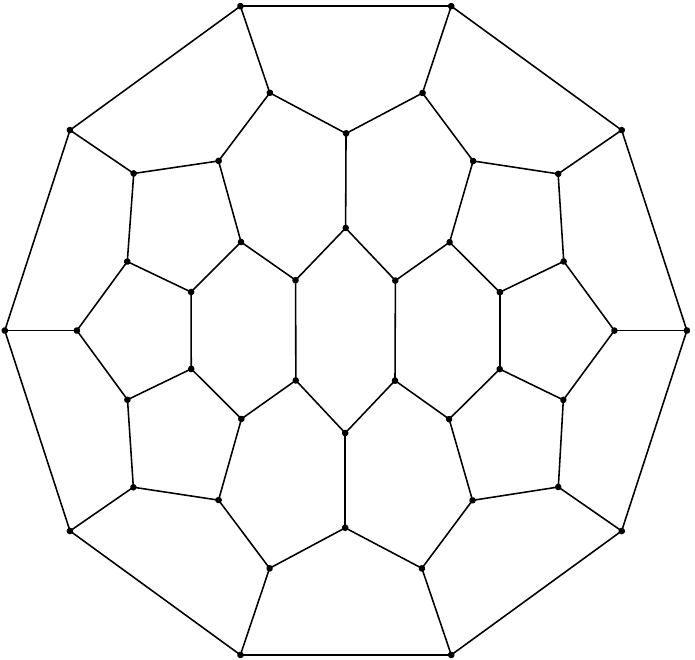}\par
$10-3$: 44, $C_{2v}$
\end{minipage}
\begin{minipage}[b]{3.0cm}
\centering
\epsfig{height=2.8cm, file=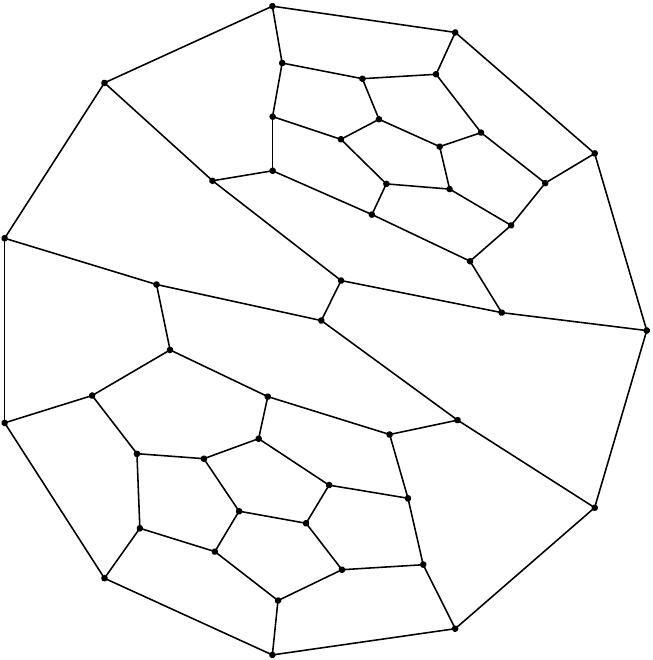}\par
$11-1$: 48, $C_1$
\end{minipage}
\begin{minipage}[b]{3.0cm}
\centering
\epsfig{height=2.8cm, file=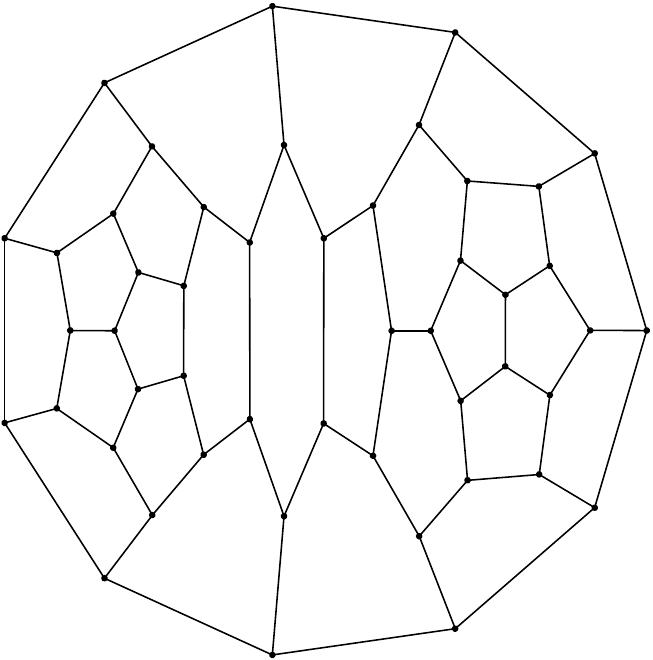}\par
$11-2$: 48, $C_s$
\end{minipage}
\begin{minipage}[b]{3.0cm}
\centering
\epsfig{height=2.8cm, file=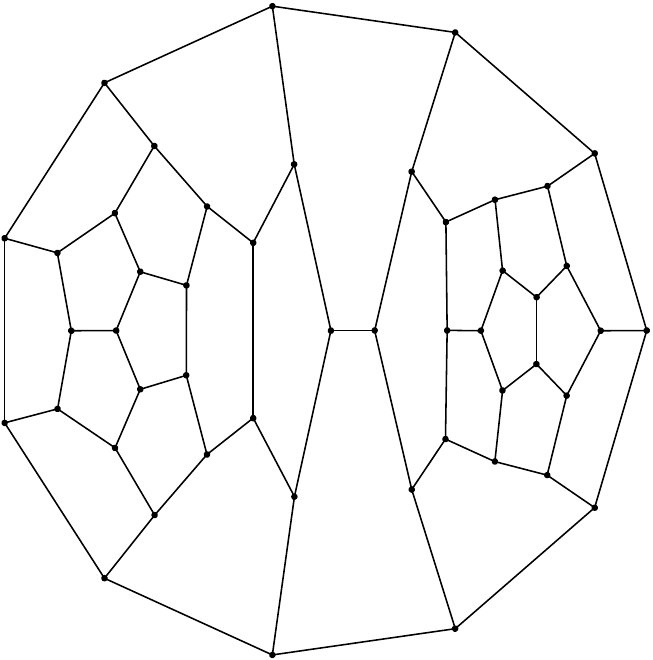}\par
$11-3$: 48, $C_s$
\end{minipage}
\begin{minipage}[b]{3.0cm}
\centering
\epsfig{height=2.8cm, file=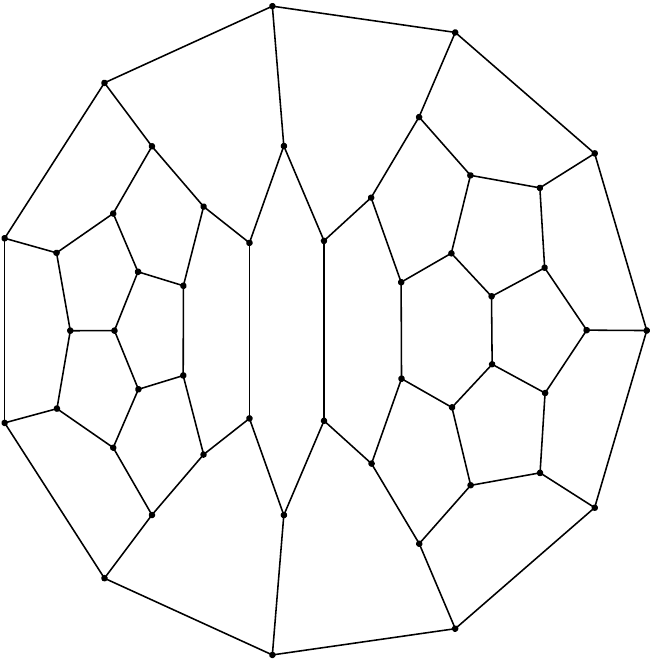}\par
$11-4$: 48, $C_s$
\end{minipage}
\begin{minipage}[b]{3.0cm}
\centering
\epsfig{height=2.8cm, file=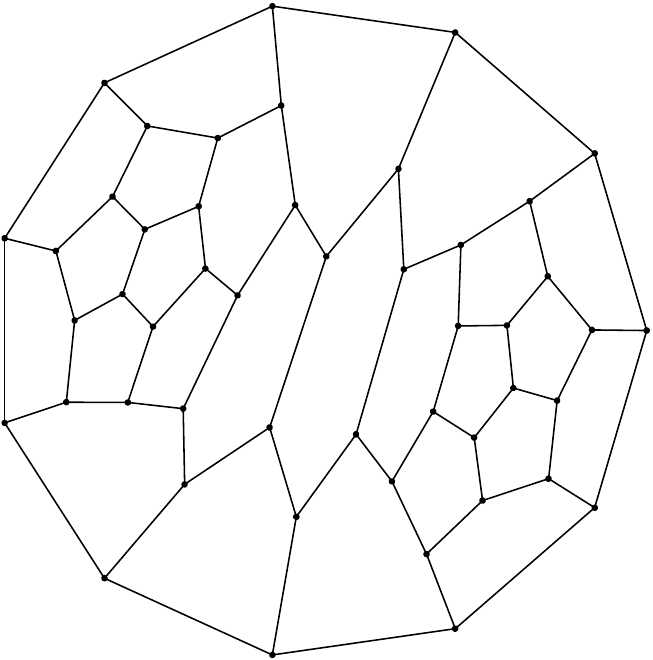}\par
$11-5$: 48, $C_1$
\end{minipage}
\begin{minipage}[b]{3.0cm}
\centering
\epsfig{height=2.8cm, file=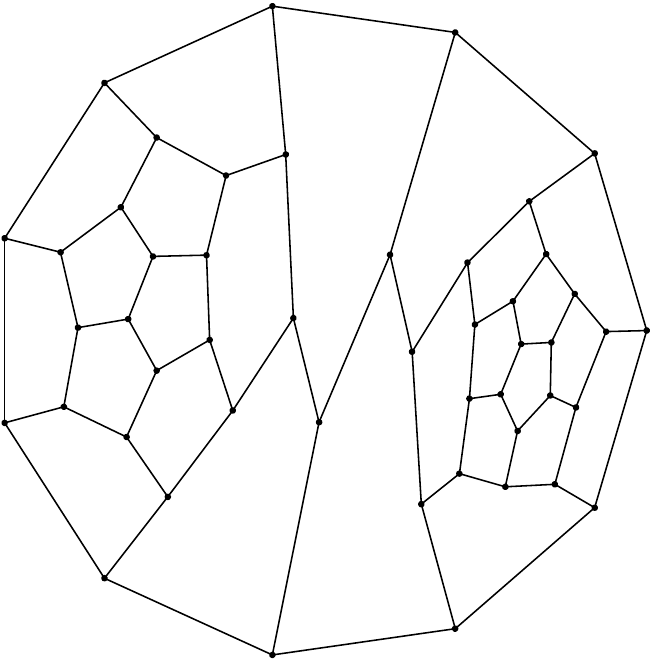}\par
$11-6$: 48, $C_1$
\end{minipage}
\begin{minipage}[b]{3.0cm}
\centering
\epsfig{height=2.8cm, file=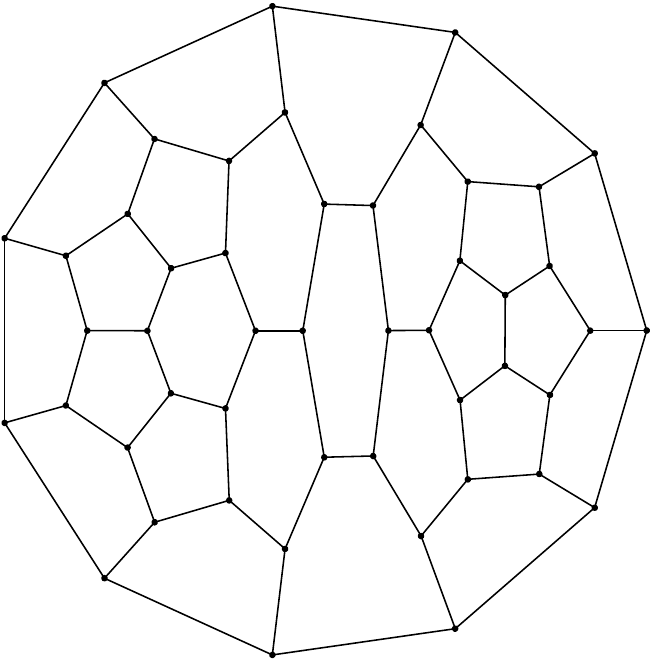}\par
$11-7$: 48, $C_s$
\end{minipage}
\begin{minipage}[b]{3.0cm}
\centering
\epsfig{height=2.8cm, file=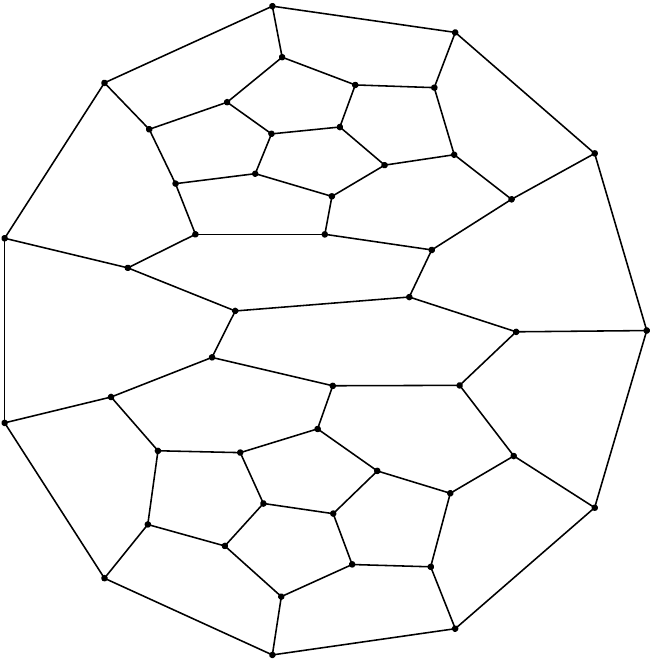}\par
$11-8$: 48, $C_1$
\end{minipage}
\begin{minipage}[b]{3.0cm}
\centering
\epsfig{height=2.8cm, file=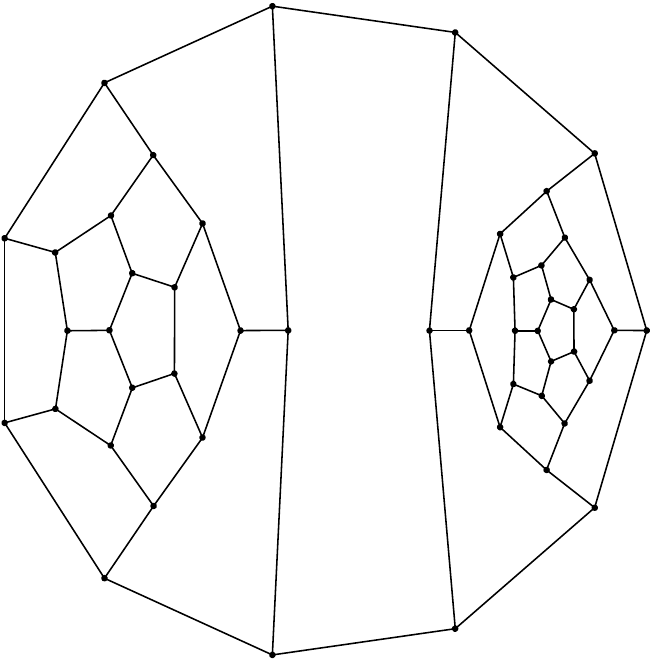}\par
$11-9$: 48, $C_s$
\end{minipage}
\begin{minipage}[b]{3.0cm}
\centering
\epsfig{height=2.8cm, file=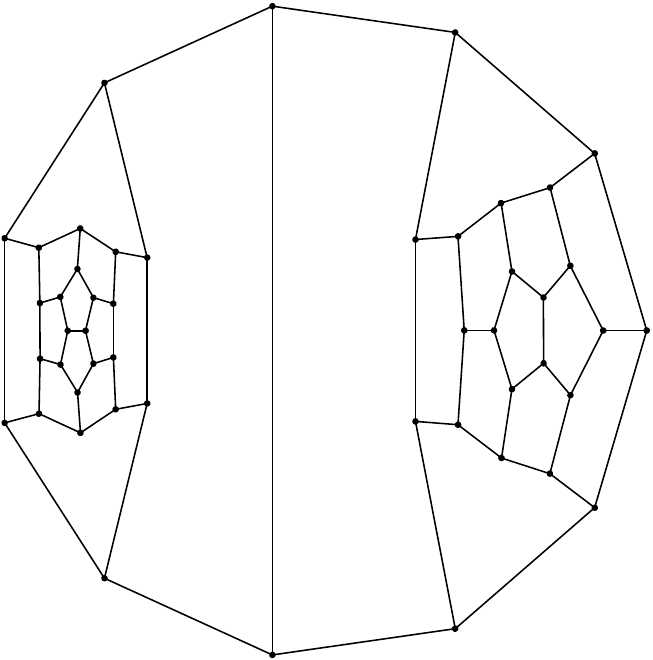}\par
$11-{10}$: 48, $C_s$
\end{minipage}
\begin{minipage}[b]{3.0cm}
\centering
\epsfig{height=2.8cm, file=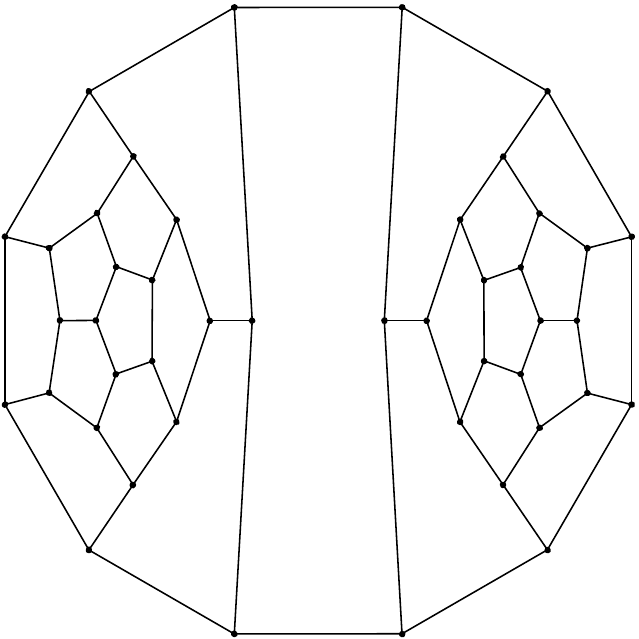}\par
$12$: 44, $C_{2v}$
\end{minipage}
\begin{minipage}[b]{3.0cm}
\centering
\epsfig{height=2.8cm, file=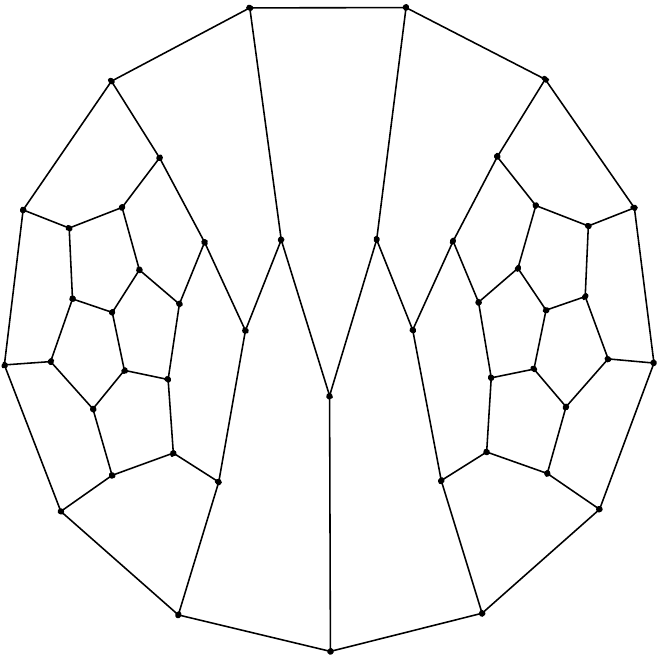}\par
$13$: 48, $C_{s}$
\end{minipage}
\begin{minipage}[b]{3.0cm}
\centering
\epsfig{height=2.8cm, file=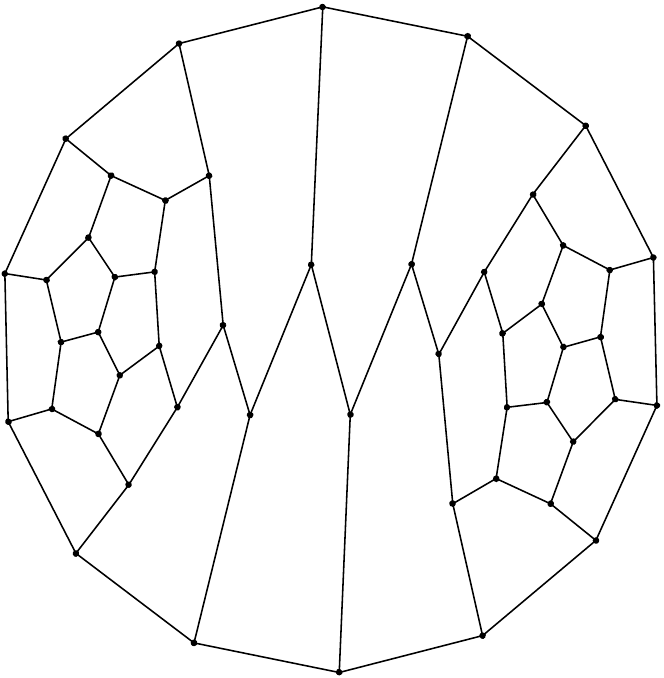}\par
$14$: 50, $C_{2}$
\end{minipage}

\end{center}
\caption{Minimal fullerene $c$-disks with $1 \le c\le 14$, $c\not=5,6$}
\label{MinimalNanodisk}
\end{figure}

Moreover, we conjecture that for odd or even $c\ge 5$, the following
$(5c-5)$- or $(5c-6)$-vertex $c$-thimble is a minimal one; it holds
for $5\le c\le 10$ since this construction generalizes cases $7, 9-1$
and $8$, $10-3$ in Figure \ref{MinimalNanodisk}, as well as
cases $c=5$ and $6$, of minimal $c$-disks.
Take the $c$-ring of $5$-gons, then then put inside a concentric $c$-ring of 
$5$- or $6$-gons: $3$-path of $5$-gons and, on opposite side,
$3$-path or $3$-ring, for even or odd $c$, of $5$-gons.
Remaining $c-6$ or $c-5$ $5$- or $6$-gons
of interior $c$-ring are $6$-gons.
Finally, fill inside of the interior $c$-ring by the $\frac{c-4}{2}$-
or $\frac{c-5}{2}$-path of $6$-gons.
Another generalization of the minimal $6$-disk is, for $c=6t$, 
$\frac{c(c+18)}{6}$-vertex $c$-thimble of symmetry (for $c>6$) $C_6$ or $C_{6v}$.
In fact, take $6t$-ring of $5$-gons, then, inside of it, $6t$-ring, where six equispaced $5$-gons are alternated by $(t-1)$-tuples of $6$-gons.
Finally, put inside, for $i=t-1$, $t-2\dots$, $1$, consecutively
$6i$-rings of $6$-gons.

\begin{figure}
\begin{center}
\begin{minipage}[b]{3.0cm}
\centering
\epsfig{height=2.8cm, file=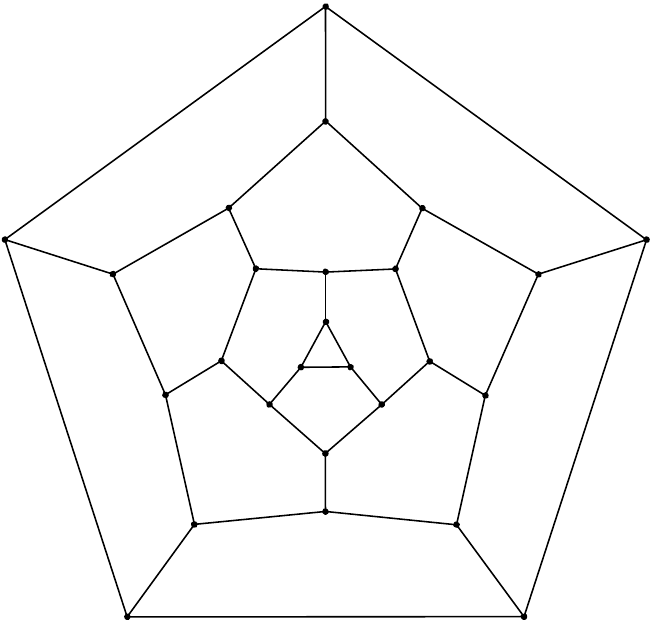}\par
\end{minipage}
\end{center}
\caption{Minimal $3$-disk of type $(6,6,5)$; it is also $C_s$-minimal $3$-disk}
\label{Exmp566}
\end{figure}

Any such $c$-thimble
can be elongated by adding an
outside  ring $5$-gons along the $c$-gon and transforming inside ring
of $5$-gons along it into a ring of $6$-gons.
Let a $(\{5,6,c\},3)$-sphere with $p_c=1$ have a simple {\em zigzag}
(left-right circuit without self-intersections).
A {\em railroad} is a circuit of $6$-gons, each of which is adjacent 
to its neighbors on opposite edges.
Let us elongate above sphere by a railroad  
along zigzag and  then let us cut elongated sphere in the middle of
this ring. We will get two $c$-thimbles.

\begin{theorem}
All simple zigzags of any elongated $c$-thimble are parallel
and its railroads are parallel $c$-rings of $6$-gons, forming a cylinder.
\end{theorem}

In fact, suppose that there is a railroad in an $c$-thimble
not belonging to the cylinder of parallel $c$-rings of
$6$-gons along the boundary $c$-ring of $5$-gons.

Let us cut it in the middle. The thimble is separated into an 
smaller $c$-thimble with $p_5=c$ and a $(\{5,6,c\},3)$-sphere with $p_5=6$.
But this sphere can not contain at least $c+1$ of
$c+6$ original $5$-gons: $c$ boundary $5$-gons
and, at other side of the cylinder, at least one $5$-gon.
Only a railroad belonging to the cylinder can go around this $5$-gon. \qed

Given a $c$-disk, call its {\em type}, the sequence of gonalities of 
its consecutive neighbors. So, any $c$-thimble has type $(5,\dots, 5)$.  

There are bijections between $v$-vertex $1$-disks, $(v-2)$-vertex $2$-disks
of type $(5,6)$ and $(v-4)$-vertex  $(\{5,6,3,4\},3)$-spheres with unique
and adjacent $3$ and $4$-gon. Using it, we found the minimal $1$-disk, 
given in Figure \ref{MinimalNanodisk}.
There is a bijection between $v$-vertex $3$-disks of type $(6,6,6)$ and 
$(v-2)$-vertex fullerenes with $(5,5,5)$-vertex (collapse the $3$-gon).
There is a bijection between $v$-vertex $3$-disks of type $(5,6,6)$ and    
$(v-2)$-vertex fullerenes with $(5,5,6)$-vertex (delete the edge adjacent 
to $5$-gon).
So, the possible types (and minimal number of vertices for examples) are:
$(6)$ ($v=40$) for $1$-disks; $(6,6)$ ($v=26$), $(6,5)$ ($v=40-2$) for
$2$-disks; $(6,6,6)$ ($v=22$), $(6,5,5)$ ($v=40-6$), $(6,6,5)$
($v=26$, see Figure \ref{Exmp566}) for $3$-disks.

One can see unique $(\{5,6\},3)$-spheres with $v=20$, $24$ as minimal
fullerene $c$-disks with $c=5$, $6$.
Clearly, the minimal fullerene $3$- and $4$-disks 
are $22$-vertex spheres obtained from Dodecahedron by truncation on
a vertex or an edge.
As well, $C_1$-minimal $3$-disk and $C_2$-minimal $4$-disk are truncations
of $28$-vertex ($D_2$) or $24$-vertex $(\{5,6\},3)$-sphere on,
respectively, a vertex or an edge.

\begin{conjecture}
An $v$-vertex fullerene $c$-disk with $c\ge 1$, 
except the cases $(c,v)=(1,42)$, $(3,24)$ and $(5,22)$, exists
if and only if $v$ is even and $v \ge 2(p_6(c)+c+5)$.
Here $p_6(c)$ denotes the minimal possible number of $6$-gons in a fullerene 
$c$-disk.
\end{conjecture}
We have $p(1)=14$, $p(2)=6$, $p(3)=3$, $p(4)=2$, $p(5)=0$, $p(6)=1$, $p(7)=3$, $p(8)=4$, $p(9)=6$, $p(10)=7$, $p(11)=8$, $p(12)=5$, and $p(c)=6$ for $c\geq 13$.
For $2\leq c\leq 20$, the conjecture was checked by computation 
and all minimal fullerene $c$-disks ($2,3,10$ for $c=9,10,11$ and
unique, otherwise) are listed; see Figure \ref{MinimalNanodisk}
for $1\le c\le 14$, $c\neq 5,6$.
All $c$-disks there having $c\ge 3$, except three among $11$-disks,
are $3$-connected.

The {\em $c$-pentatube} is
$2(c+11)$-vertex fullerene $c$-disk of symmetry $C_2$, $C_s$
for even, odd $c$, respectively.
Its $5$-gons are organized in two $(5,3)$-polycycles $B_2$ separated
by $Pen_{c-12}$ and its six $6$-gons are organized into two $3$-rings
each shielding a $B_2$ from $Pen_{c-12}$
The $c$-pentatube is unique minimal for $13\le c\le 20$;
we expect that it remains so for any $c\ge 13$.

\begin{theorem}
The possible symmetry groups of a 
$(\{5,6,c\},3)$-sphere with $p_c=1$ and $c\neq 5,6$
are $C_n$, $C_{nv}$ with $n\in \{1,2,3,5,6\}$ and $n$ dividing $c$.
\end{theorem}
In fact, any symmetry of a $(\{5,6,c\},3)$-sphere 
should stabilize unique $c$-gon.
So, the possible groups are only $C_n$ and $C_{nv}$ with $n$
($1\le n\le c$) dividing $c$.
Moreover, $n\in \{1,2,3,5,6\}$ since, on the axis has
to pass by a vertex, edge or face. Remind that $C_s=C_{1v}$. \qed

\begin{figure}
\begin{center}
\begin{minipage}[b]{3.0cm}
\centering
\epsfig{height=2.6cm, file=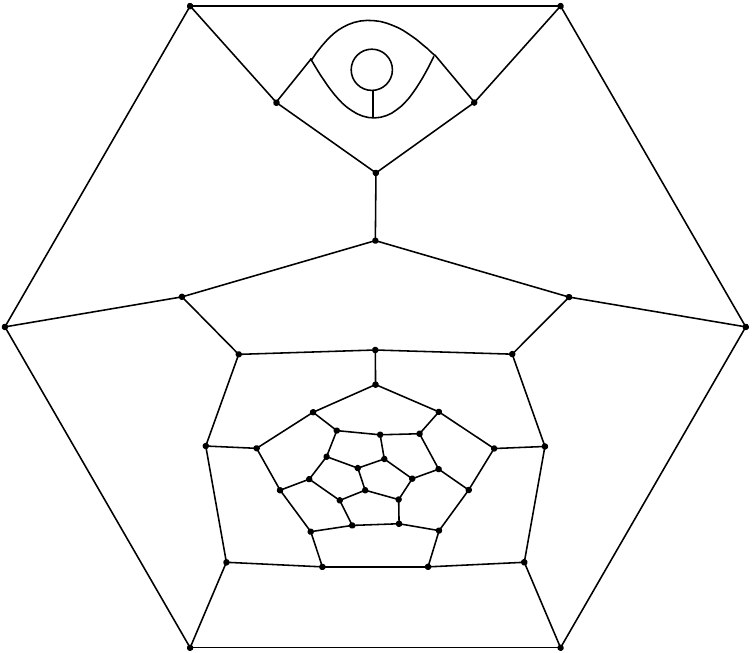}\par
$c=1$: 48, $C_{1}$
\end{minipage}
\begin{minipage}[b]{3.0cm}
\centering
\epsfig{height=2.6cm, file=ExmpB1_sec.pdf}\par
$c=1$: 40, $C_{s}$
\end{minipage}
\begin{minipage}[b]{3.0cm}
\centering
\epsfig{height=2.6cm, file=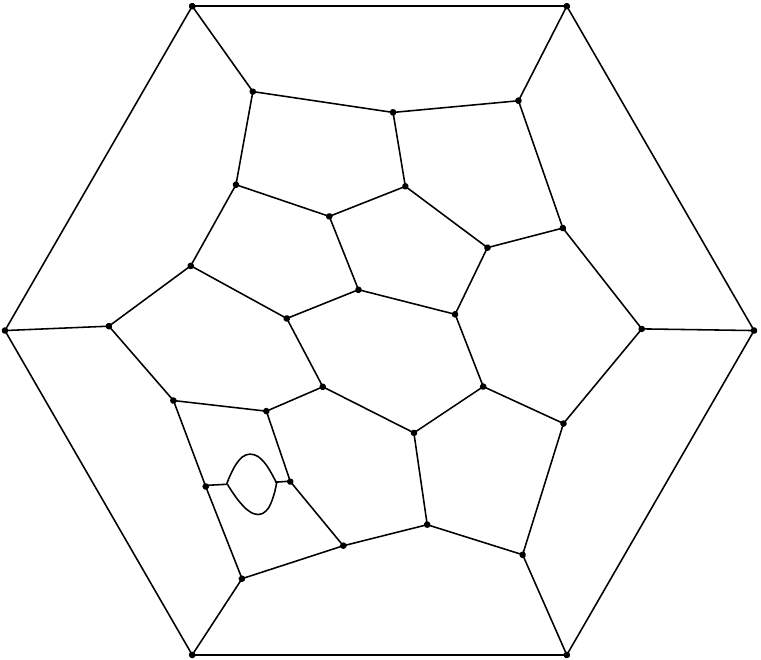}\par
$c=2$: 34, $C_{1}$
\end{minipage}
\begin{minipage}[b]{3.0cm}
\centering
\epsfig{height=2.6cm, file=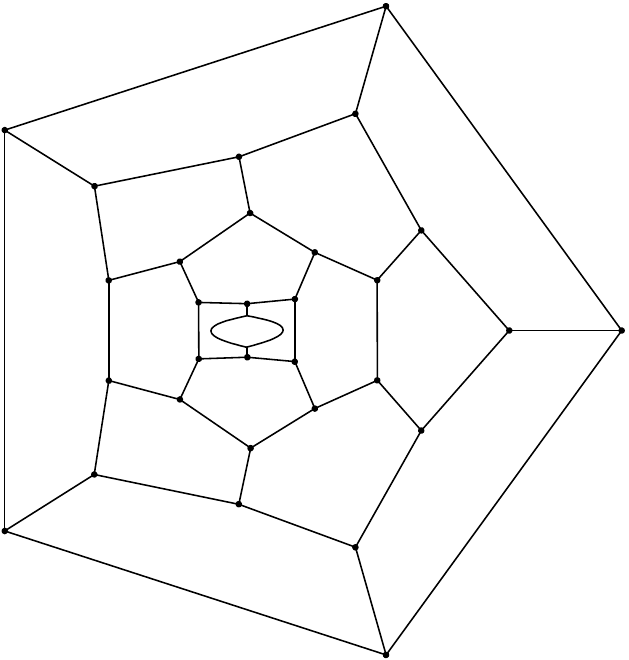}\par
$c=2$: 32, $C_{s}$
\end{minipage}
\begin{minipage}[b]{3.0cm}
\centering
\epsfig{height=2.6cm, file=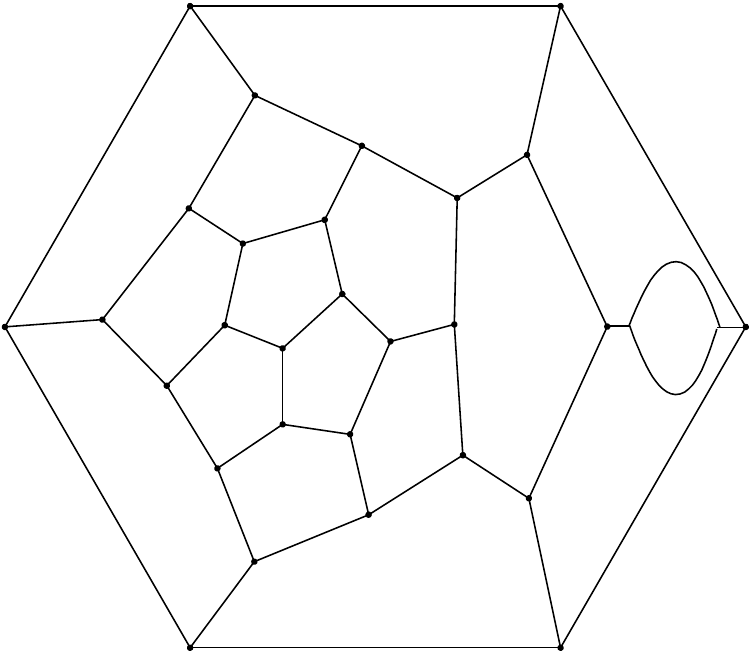}\par
$c=2$: 30, $C_{2}$
\end{minipage}
\begin{minipage}[b]{3.0cm}
\centering
\epsfig{height=2.4cm, file=MinB2_second.pdf}\par
$c=2$: 26, $C_{2v}$
\end{minipage}
\begin{minipage}[b]{3.0cm}
\centering
\epsfig{height=2.8cm, file=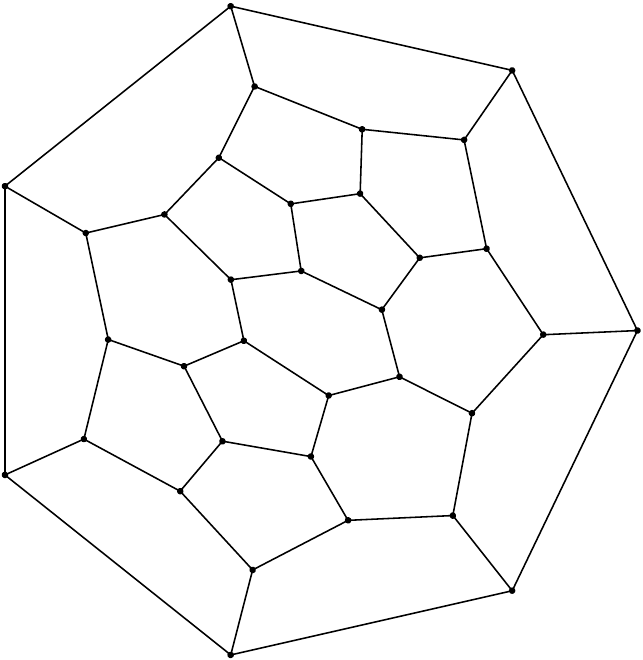}\par
$c=7$: 34, $C_1$
\end{minipage}
\begin{minipage}[b]{3.0cm}
\centering
\epsfig{height=2.8cm, file=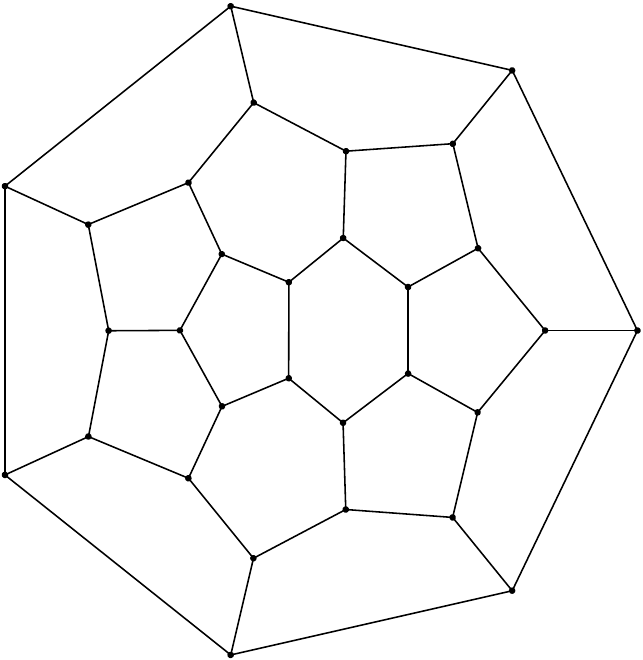}\par
$c=7$: 30, $C_s$
\end{minipage}
\begin{minipage}[b]{3.0cm}
\centering
\epsfig{height=2.8cm, file=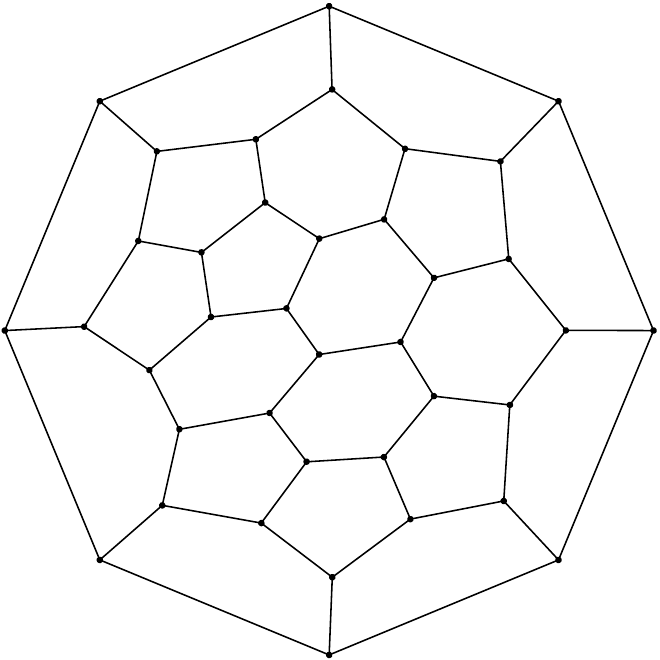}\par
$c=8$: 38, $C_1$
\end{minipage}
\begin{minipage}[b]{3.0cm}
\centering
\epsfig{height=2.8cm, file=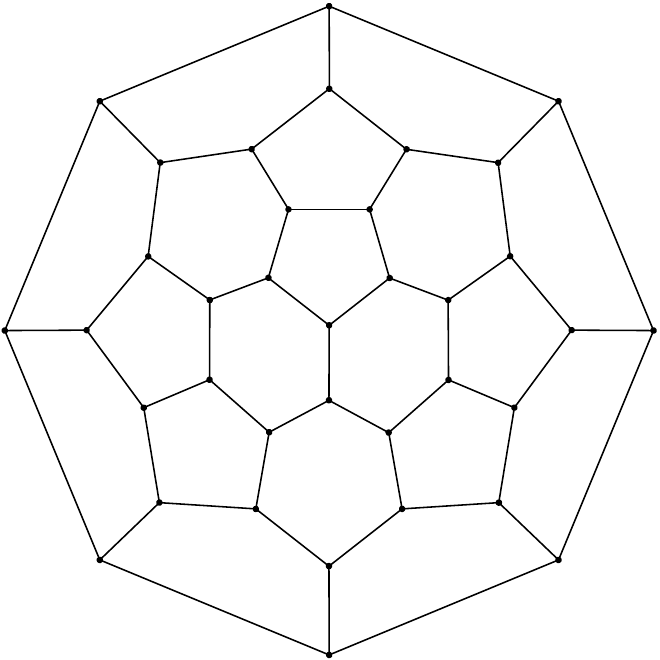}\par
$c=8$: 36, $C_s$
\end{minipage}
\begin{minipage}[b]{3.0cm}
\centering
\epsfig{height=2.8cm, file=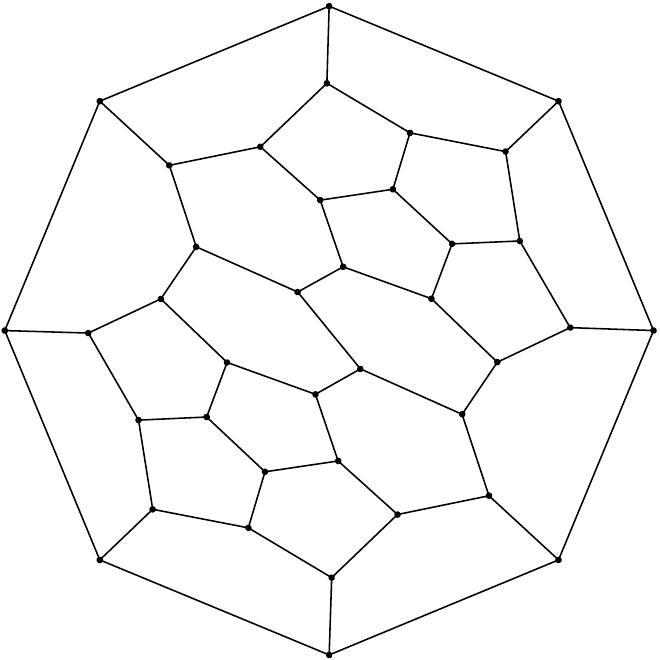}\par
$c=8$: 38, $C_2$
\end{minipage}
\begin{minipage}[b]{3.0cm}
\centering
\epsfig{height=2.8cm, file=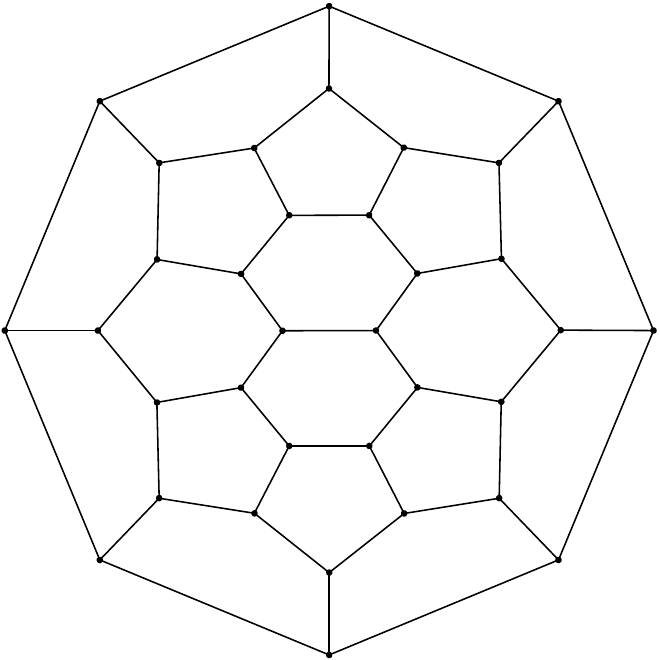}\par
$c=8$: 34, $C_{2v}$
\end{minipage}

\begin{minipage}[b]{3.0cm}
\centering
\epsfig{height=2.8cm, file=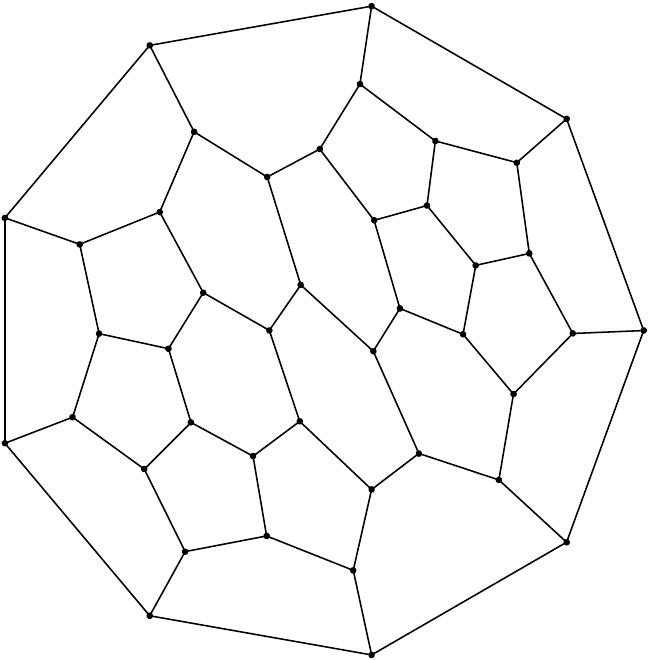}\par
$c=9$: 42, $C_1$
\end{minipage}
\begin{minipage}[b]{3.0cm}
\centering
\epsfig{height=2.8cm, file=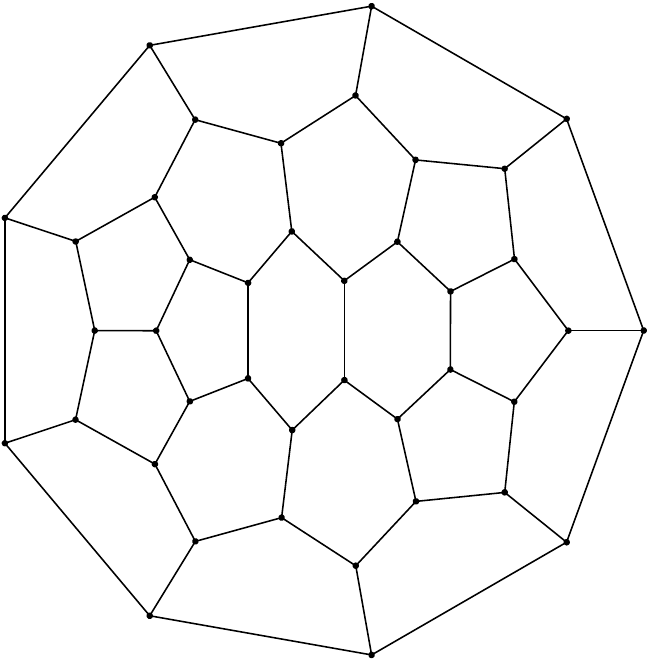}\par
$c=9$: 40, $C_s$
\end{minipage}
\begin{minipage}[b]{3.0cm}
\centering
\epsfig{height=2.8cm, file=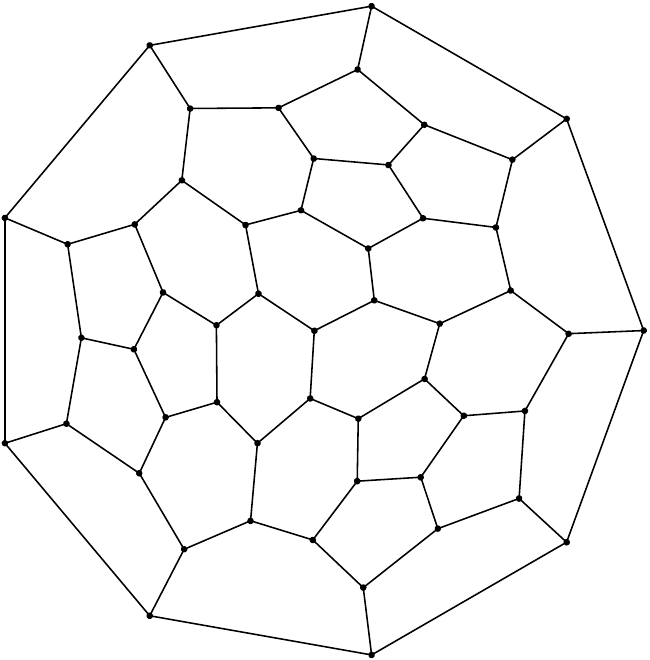}\par
$c=9$: 52, $C_3$
\end{minipage}
\begin{minipage}[b]{3.0cm}
\centering
\epsfig{height=2.8cm, file=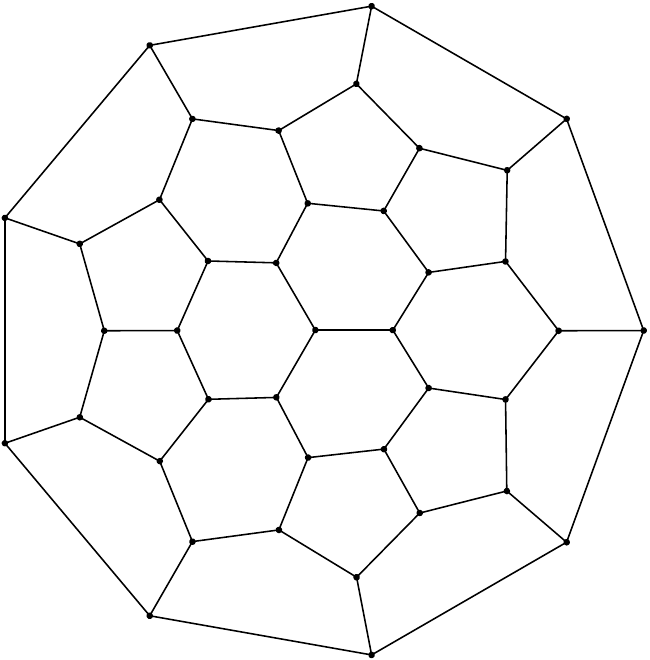}\par
$c=9$: 40, $C_{3v}$
\end{minipage}

\end{center}
\caption{Minimal fullerene $c$-disks for each possible group with
$c=1$, $2$, $7$, $8$, and $9$}
\label{AllKnownGrp_sph789}
\end{figure}

Cases (xi) and (xvi) in \cite{Group345} show that the possible symmetry 
groups (minimal $v$) of a fullerene $c$-disk with $c=3$ and $c=4$
are $C_1\,(30)$, $C_s\,(26)$, $C_{3}\,(34)$, $C_{3v}\,(22)$
and $C_1\,(28)$, $C_s\,(24)$, $C_{2}\,(26)$, $C_{2v}\,(22)$.
Minimal examples are given there. 
For $c=1$, $2$, $7$, $8$, $9$ see minimal examples on Figure
\ref{AllKnownGrp_sph789}.
For $c=2$, those examples are all of type $(6,6)$ and coming
from one-edge truncation of $(\{5,6\},3)$-spheres with, respectively,
$28$ ($D_2$), $26$, $24$ and $20$ vertices.

\end{document}